\documentclass[letterpaper,12pt,leqno,oneside]{article}
\usepackage{color}
\usepackage{bm}
\usepackage{exscale}
\usepackage{amsmath}
\usepackage{amsfonts}
\usepackage{pagecolor,lipsum}
\usepackage{xcolor}
\usepackage{stmaryrd}
\usepackage{amscd}
\usepackage{graphicx}
\usepackage{marvosym}
\usepackage{amsxtra}
\usepackage{amssymb}
\usepackage{theorem}
\usepackage[final]{epsfig}
\newtheorem{proposition}{Proposition}[subsection]
\newtheorem{definition}[proposition]{Definition}

\newtheorem{lemma}[proposition]{Lemma}

{\theorembodyfont{\rmfamily}\newtheorem{remark}[proposition]{Remark}}

\newtheorem{corollary}[proposition]{Corollary}

{\theorembodyfont{\rmfamily}\newtheorem{example}[proposition]{Example}}

\newfont{\abc}{cmtt10 scaled 1200}

\def\R{\mathbb{R}}

\def\Z{\mathbb{Z}}

\def\U{\mathbb{U}}
\def\D{\mathbb{D}}
\def\P{\mathbb{P}}
\def\C{\mathbb{C}}

\def\P{\mathbb{P}}
\def\H{\mathbb{H}}

\def\E{\mathbb{E}}
\def\U{\mathbb{U}}
\def\N{\mathbb{N}}

\def\I{\mathbb{I}}
\def\ve{\varepsilon}

\def\ra{\rightarrow}
\def\cs{\symbol{35}}
\def\p{\partial}
\def\qed{\hfill $\Box$ \\}
\def\mm{\mbox}
\def\v{= \emptyset}
\def\n{\neq \emptyset}
\def\D{\mathbf{ID}}

\def\bp{\langle A \rangle}

\def\su{\mathbb{SD}}

\begin{document}

\vspace*{0cm}

\begin{center}\Large{\bf{Hyperbolic Geometry and Potential Theory on Minimal Hypersurfaces}}\\
\medskip
{\small{by}}\\
\medskip
\large{\bf{Joachim Lohkamp}}\\
\end{center}

\vspace{0.2cm}
\noindent Mathematisches Institut, Universit\"at M\"unster, Einsteinstrasse 62, Germany\\
 {\small{\emph{e-mail: j.lohkamp@uni-muenster.de}}}
\vspace{0.2cm}

{\small
 {\center \tableofcontents}
{\contentsline {subsection}{\numberline {}References}{. .}}}

\vspace{0.3cm} \setcounter{section}{1}
\renewcommand{\thesubsection}{\thesection}
\subsection{Introduction} \label{intro}

\bigskip

This is the second in a series of three papers which started in [L1]  where we established basic skin structural concepts and results for area minimizing hypersurfaces. In the present paper we show how this leads to previously unapproachable and largely unexpected geometric and analytic properties of such hypersurfaces. We do this with an emphasis on the case of singular hypersurfaces.\\

We first show that these generally wrinkled and degenerated spaces admit canonical conformal unfoldings to complete Gromov hyperbolic spaces with bounded geometry and we recover their singular set as the Gromov boundary. This new geometric facet of area minimizers helps us to also reveal analytic properties of these spaces, one commonly  merely expects in cases like regularly bounded Euclidean domains.  This includes potential theoretic details, like asymptotic regularity and extension results near or even to the singular set, for classical elliptic operators. \\

We specify our main results below and recommend to also consult [L1],Ch.1.1 for a broad overview of these concept and of this series of papers, including the third one [L2] where we apply these techniques to problems in scalar curvature geometry.\\

\textbf{Skin Structures} \quad For starters, we recall the basic notions of \emph{skin transforms} and \emph{skin uniformity} introduced in [L1].\\

To this end we henceforth consider a connected locally area minimizing hypersurface $H^n \subset M^{n+1}$ without
boundary in some smooth Riemannian $n+1$-manifold $(M,g_M)$, where either $M$ is compact, or equals $\R^{n+1}$ with its Euclidean metric. $\Sigma \subset H$ denotes its singular set.\\

On the class ${\cal{H}}$ of these hypersurfaces we consider a skin transform $\bp$. We recall from [L1] that an assignment $H \mapsto \bp_H$, to any $H \in {\cal{H}}$, that
commutes with the convergence of sequences of area minimizers, is called a \emph{skin transform} provided
\begin{itemize}
    \item $\bp_H \ge |A_H|$ and for any $f \in C^\infty(H \setminus \Sigma,\R)$ compactly supported in $H \setminus \Sigma$
{\small \begin{equation}\label{hh}
 \int_H|\nabla f|^2  + |A_H|^2 \cdot f^2 dA \ge \tau \cdot \int_H \bp_H^2\cdot f^2 dA, \mm{ for some } \tau = \tau(\bp,H) \in (0,1).
\end{equation}}
    \item $\bp_H \equiv 0$, if $H \subset M$ is totally geodesic*, otherwise, $\bp_H$ is  strictly positive.
    \item When $H$ is not totally geodesic, we define $\delta_{\bp}:=1/\bp$, the \emph{${\bp}$-distance}. It is $L_{\bp}$-Lipschitz regular, for some constant $L_{\bp}=L(\bp,n)>0$:
        \[|\delta_{\bp}(p)- \delta_{\bp}(q)|   \le L_{\bp} \cdot d(p,q), \mm{ for } p,q \in  H \setminus \Sigma. \]
\end{itemize}
We occasionally call (\ref{hh}) the \emph{Hardy inequality} for $\bp$, which actually  is a Hardy inequality for the operator
$-\Delta + |A|^2$ relative to the $\bp$-distance $\delta_{\bp}$, cf. [L1],Ch.3.2 in particular Rm.3.7 for a discussion. (*It is easy to check that any totally geodesic $H \in {\cal{H}}$ must be a regular hypersurface.)\\

With this concept we can formulate and prove that $\Sigma$ can be approached in a quantitatively non-tangential way from $H \setminus \Sigma$:\\

$H \setminus \Sigma$  is a $c$-\emph{skin uniform} space, for some $c >0$, i.e. any pair $p,q \in H \setminus \Sigma$ can be joined by a rectifiable path
$\gamma_{p,q}: [a,b] \ra H \setminus \Sigma$, for some $a <b$, with $\gamma_{p,q}(a)=p$, $\gamma_{p,q}(b)=q$, so that for any $z \in \gamma_{p,q}$:
\[l(\gamma)  \le c \cdot  d(p,q) \, \mm { and }  \, l_{min}(\gamma_{p,q}(z)) \le c \cdot \delta_{\bp}(z),\]
where $l_{min}(\gamma_{p,q}(z))$ is the minimum of the lengths of the subarcs of $\gamma_{p,q}$ from $p$ to $z$ and from $q$ to $z$.\\

\textbf{Statement of Results} \, Now we describe the contents of this paper and state our main results.\\

\textbf{Notations} \, Our hypersurfaces are locally mass minimizing, integer multiplicity rectifiable currents of dimension $n$ without boundary.
The partial regularity theory for these minimizers says that $H$ is a smooth hypersurface except for some singular set $\Sigma_H$ of Hausdorff-dimension $\le n-7$, cf.[L1],Appendix, for details.\\

 We consider the following classes of complete area minimizing hypersurfaces.
\begin{description}
    \item[${\cal{H}}^c_n$]:    $H \subset M$ is a compact and connected hypersurface without boundary.
    \item[${\cal{H}}^{\R}_n$]:  $(M,g_M) = (\R^{n+1},g_{\R^{n+1}})$, $H$ is an oriented boundary of some open set $A \subset \R^{n+1}$  and, thus, $H$ is non-compact and
        complete.
            \item[$\mathcal{C}_{n}$]:   $\mathcal{C}_{n} \subset {\cal{H}}^{\R}_n$ is the space of area minimizing $n$-cones in $\R^{n+1}$ with tip in $0$.
            \item[$\mathcal{SC}_{n}$]:    $\mathcal{SC}_n \subset\mathcal{C}_n$ is the subset of cones singular, at least, in $0$
    \item[$\mathcal{K}_{n-1}$]: For any area minimizing cone $C \subset \R^{n+1}$ with tip $0$, we get the non-minimizing minimal hypersurface $S_C$ in the unit sphere
\[S_C:= \p B_1(0) \cap C \subset S^n \subset  \R^{n+1} \mm{ and set } {\cal{K}}_{n-1}:= \{ S_C\,| \, C \in {\mathcal{C}_{n}}\}.\]
We write ${\cal{K}}= \bigcup_{n \ge 1} {\cal{K}}_{n-1}$, for the space of all such hypersurfaces $S_C$.
\end{description}
The main class of hypersurfaces we study in this paper is given by
\[{\cal{H}}_n:= {\cal{H}}^c_n \cup {\cal{H}}^{\R}_n \, \mm{ and }\, {\cal{H}} :=\bigcup_{n \ge 1} {\cal{H}}_n.\]
${\cal{H}}_n$ is closed under blow-ups. That is, the limit of converging subsequences under scaling by a diverging sequence of real numbers belongs to ${\cal{H}}^{\R}_n$.\\

To better visually identify the cases under consideration, we generally write $\Sigma_H$ for the  singular set of any $H \in {\cal{H}}_n$, but we write $\sigma_C$ for the (non-compact) singular set of a cone $C \in \mathcal{SC}_{n}$, when we want to emphasize that $C$ is a tangential object and we use its cone properties. (This is inspired from the upper versus lower case symbols used for Lie groups versus Lie algebras.)\\

 \textbf{Hyperbolic Geometry} \quad The Riemann uniformization theorem shows that any domain $D \subset \R^2 = \C$, $D \neq \C$ or $\C^*$, is conformally
 equivalent to a complete hyperbolic metric. However, for many applications it is more useful to ensure that it is Gromov hyperbolic. This requires some basic boundary regularity expressed in terms of uniformity properties of $D$.\\

Such uniformity concepts can be defined in arbitrary dimensions: $D \subset \R^n$ is a uniform domain when any pair $p,q \in D$ can be joined by a $c$-\emph{uniform curve}
 in $D$, i.e. a rectifiable path $\gamma_{p,q}: [a,b] \ra D$, for some $a <b$, with $\gamma_{p,q}(a)=p$,
    $\gamma_{p,q}(b)=q$, so that
\[ l(\gamma)  \le c \cdot  d(p,q) \, \mm{ and }\, l_{min}(\gamma_{p,q}(z)) \le c \cdot dist(z,\p D)\]
for any $z \in \gamma_{p,q},$  where $l_{min}(\gamma_{p,q}(z))$ is the minimum of the lengths of the subarcs of $\gamma_{p,q}$ from $p$ to $z$ and from $q$ to $z$.\\

Now we consider the \emph{quasi-hyperbolic metric} $k_D$ which keeps track of the distortion relative $\p D$, cf.[GO],[He],[Ko] or [BHK]: for any two points $x,y \in D$ defined by
\[k_D(x,y) := \inf \Bigl  \{\int_\gamma 1/dist(\cdot,\p D)   \, \, \Big| \, \gamma   \subset D \mbox{ rectifiable arc joining }  x \mbox{ and } y  \Bigr \}\]
Now the basic fact is that when $D$ is a \emph{uniform domain}, $(D,k_D)$ is a complete Gromov hyperbolic space of bounded geometry and
$\p D$ equals the Gromov boundary. This result is due to Gehring, Osgood resp. Bonk, Heinonen, Koskela, cf.[GO], [BHK].\\

The concept of uniformity and quasi-hyperbolic metrics can be extended to  more general \emph{non-complete} metric spaces. Here we view the points added in the metric completion as the boundary.  The quasi-hyperbolic metric for uniform spaces is complete, Gromov hyperbolic and the boundary of the metric completion equals the Gromov boundary, cf.[BHK] and [He]. \\

 Since we have seen in [L1] that $H \setminus \Sigma$ is a uniform space, we can appeal to the latter theory and infer a first coarse hyperbolization result.
We refer to Ch.\ref{discrete} and \ref{gh} for details concerning notions and statements which appear in Theorem 1 and 2.\\

 \textbf{Theorem 1} \textbf{(Quasi-Hyperbolic Geometry)} {\itshape \,  For any singular hypersurface $H \in {\cal{H}}$, that is, with $\Sigma_H \n$, the quasi-hyperbolic metric $k_{H \setminus \Sigma}$ on $H \setminus \Sigma$,
 \[k_{H \setminus \Sigma} := \inf \Bigl \{\int_\gamma 1/dist(\cdot,\Sigma_H)  \, \, \Big| \, \gamma  \subset H \setminus \Sigma \mbox{ rectifiable arc joining }  x \mbox{ and } y  \Bigr \}\]
is defined and has the following properties:
\begin{itemize}
\item   $(H \setminus \Sigma, k_{H \setminus \Sigma})$ is a complete, Gromov hyperbolic and visual metric space.
\end{itemize}}
\bigskip
Towards applications of uniformity and hyperbolizations, we note that smooth but highly curved portions of $H$ have a considerable impact on the elliptic analysis on $H$ not reflected from the uniformity of $H \setminus \Sigma$. For the same reason, it is difficult to understand this analysis relative $k_{H \setminus \Sigma}$.\\

For instance, wrinkles of $H \setminus \Sigma$ can cause degenerations of $(H \setminus \Sigma, k_{H \setminus \Sigma})$ near infinity. In other words, in general, this space does \emph{not} have a bounded geometry and this makes standard analytic  tools, like uniform Harnack inequalities on equally sized domains, unavailable.  Moreover, $k_{H \setminus \Sigma}$ may drastically alter under deformations of $H$. That is, it does not commute with the convergence of area minimizers and this limits its usability in blow-up and  compactness arguments. \\

However, there is a more versatile and natural hyperbolization of $H \setminus \Sigma$, in terms of \emph{skin metrics} $d_{\bp}$, resolving such issues and seamlessly extending to the regular case, with $\Sigma \v$. For their definition, we replace $1/dist(\cdot,\Sigma)$ for $\bp = 1/\delta_{\bp}$ and set, for any two points $x,y \in H \setminus \Sigma$:
\[d_{\bp}(x,y) := \inf \Bigl  \{\int_\gamma  \bp \, \, \Big| \, \gamma   \subset  H \setminus \Sigma\mbox{ rectifiable arc joining }  x \mbox{ and } y  \Bigr \}.\]

The skin metrics play a central role in our treatment of elliptic problems on $H \setminus \Sigma$, in the Theorems 4 - 12 below. This makes the following result our main hyperbolization theorem for $H \setminus \Sigma$. We derive it from the \emph{skin uniformity} of $H \setminus \Sigma$:\\

 \textbf{Theorem 2} \textbf{(Conformal Hyperbolic Unfoldings)} {\itshape \,  For any hypersurface $H \in {\cal{H}}$, the skin metric $d_{\bp}$ has the following properties:
\begin{itemize}
\item The metric space $(H \setminus \Sigma, d_{\bp})$ and its quasi-isometric Whitney smoothing, the smooth Riemannian manifold $(H \setminus \Sigma, d_{\bp^*}) = (H \setminus \Sigma,
    1/\delta_{\bp^*}^2 \cdot g_H)$, are \

\textbf{complete, Gromov hyperbolic} spaces with \textbf{bounded geometry}.
\item $d_{\bp}$ commutes with the convergence of the underlying area minimizers.

 \end{itemize}
 We call the spaces $(H \setminus \Sigma, d_{\bp})$ and $(H \setminus \Sigma, d_{\bp^*})$,  \textbf{hyperbolic unfoldings}  of the conformally equivalent original space $(H \setminus \Sigma, g_H)$.}\\

\bigskip
Next we consider the Gromov boundaries of these hyperbolic spaces.
 We denote the one-point compactification of the space $A$ by $\widehat{A}$ and the point of infinity by $\infty_A$. Then we have, cf.Ch.\ref{grch}\\

 \textbf{Theorem 3} \textbf{(Gromov Boundary of  $\mathbf{H \setminus \Sigma}$)} {\itshape \,  For any singular $H \in {\cal{H}}$ the identity map on
 $H \setminus \Sigma$ extends to homeomorphisms between the on-point compactification $\widehat{H}$ and the Gromov compactifications:
\[\widehat{H}\cong\overline{(H \setminus \Sigma,d_{\bp})}_G \cong \overline{(H \setminus \Sigma,d_{\bp^*})}_G \cong \overline{(H \setminus \Sigma, k_{H \setminus \Sigma})}_G.\]
where $ \cong$ means homeomorphic. In particular, we have:
\[\widehat{\Sigma} \cong\p_G(H \setminus \Sigma,d_{\bp}) \cong \p_G(H \setminus \Sigma,d_{\bp^*}) \cong \p_G(H \setminus \Sigma, k_{H \setminus \Sigma}).\]
For these identifications, one assigns to (equivalence classes of) geodesic rays in $(H \setminus \Sigma,d_{\bp})$, $(H \setminus \Sigma,d_{\bp^*})$ and
$(H \setminus \Sigma, k_{H \setminus \Sigma})$, from some fixed base point, their end points in $\Sigma \subset H$ and $\infty_\Sigma = \infty_H$, if the rays diverge.}\\

For an area minimizing cone $C$ there are two distinguished and symmetrically positioned points $[0]$ and $[1]$ in the Gromov boundary  of $(C \setminus \sigma, d_{\bp})$:
\begin{equation}\label{cgr}
 \p_G (C \setminus \sigma, d_{\bp}) \cong [0,1]\times \Sigma_{\p B_1(0) \cap C}/\sim ,
\end{equation}
  where  $\Sigma_{\p B_1(0) \cap
        C}=\p B_1(0) \cap  \sigma_{C}$ is the singular set of $\p B_1(0) \cap C$, with  $x\sim y$,  if $x , y \in \{0\} \times \Sigma_{\p B_1(0)
        \cap C}$\, or \, $x , y \in\{1\} \times \Sigma_{\p B_1(0) \cap C}$.\\\\

\textbf{Integral Representations} \quad We employ the conformal hyperbolic unfoldings of $H \setminus \Sigma$ to analyze the potential theory and the asymptotic behavior of solutions of elliptic equations on the original space $(H \setminus \Sigma, g_H)$ near $\Sigma$. The outcome is a transparent Martin theory for a large class of elliptic operators on $H \setminus \Sigma$.\\

 To explain this program we recall the classical integral representations of harmonic functions on the unit disc viewed either as a Euclidean domain or alternatively as a version of the hyperbolic plane. The Herglotz theorem, cf.[BJ],1.7.2, shows that a function $f >0$ on the Euclidean unit disk $D^2$ is harmonic, $\Delta_{Eucl} \, f = 0$, if and only if there is a Radon measure $\mu_f$ on $S^1$ such that
\begin{equation}\label{herg}
 f(x)=2 \pi  \cdot \int_{S^1} \frac{1-|x|^2}{|x-y|^2} d\mu_f(y).
\end{equation}
On the other hand, we can equip $D$ with its conformal Poincar\'{e} metric $g_{hyp}$, which is complete and hyperbolic. For the Laplace-Beltrami operator $\Delta_{hyp}$,  we have:\,  $\Delta_{hyp}= (1-|x|)^2/4 \cdot \Delta_{Eucl}$ and, thus, the notion of (super/sub)harmonic and of the Green's functions $G$ does not change. Now it is an elementary fact that \[\frac{1-|x|^2}{|x-y|^2} = \lim_{z \ra y} G(x,z)/G(0,z)=:k(y;z), \mm{ for any } y \in S^1.\]
where the latter limit $k(y;z)$ is the so-called Martin kernel, cf. Ch.\ref{grb}.  This gives us a way to interpret (\ref{herg}) as a Martin integral: $\mu_f$ is now understood as a measure on ideal boundaries  of
the complete space $(D,g_{hyp})$, the Gromov boundary  and the Martin boundary, which in this case equal $S^1$.\\

An extension of these two viewpoints to arbitrary bounded Euclidean domains in dimensions $\ge 3$ applies to the case of a \emph{uniform} domain $D \subset \R^n$  equipped with either its Euclidean or its quasi-hyperbolic metric $k_D$. Again we find two ways to describe a representation of positive harmonic functions on $D$: \\

We can either use Aikawa's work in [Ai2] for $(D,g_{Eucl})$ to show that the Martin boundary equals $\p D$ and to formulate a counterpart of (\ref{herg}).\\

Alternatively, we can use that $(D,k_D)$ is a complete Gromov hyperbolic space of bounded geometry and $\p D$ equals the Gromov boundary. Then we employ Ancona's work, we discuss at length in Ch.3, which applies to such hyperbolic spaces to argue (under some coercivity assumptions on $\Delta$ relative $D$ we omit here) that the Gromov boundary equals the Martin boundary of the Laplacian on $(D,k_D)$ to recover the integral representation of positive harmonic functions.\\

\textbf{Potential Theory on $\mathbf{H \setminus \Sigma}$} \quad  Passing from uniform domains to  $(H \setminus \Sigma, g_H)$, a natural guess is that its (skin) uniformity could help us to develop a comparable Martin theory. But, different from the case of plane domains in $\C$ or uniform domains in $\R^n$, the wrinkled and degenerating geometry of $(H \setminus \Sigma, g_H)$ makes it a too delicate task to work on this space directly.\\

Instead, the idea is to use the skin uniformity to follow the hyperbolic alternative path indicated above and to make the conformal hyperbolic unfolding $(H \setminus \Sigma, d_{\bp^*})$ a workbench to derive analytic results on $(H \setminus \Sigma, g_H)$. \\

  To describe the elliptic problems on $(H \setminus \Sigma, g_H)$ approachable this way, we use the bounded geometry to define \emph{skin adapted charts} on $H \setminus \Sigma$: for given $K
>1$, for some $\gamma(H,K,\bp) >0$ and any $p \in H \setminus \Sigma$, we choose a $K$-bi-Lipschitz chart   $\psi_{p}: B_{\gamma / \bp(p)}(p) \ra \R^n, \, \psi_{p}(p)=0.$\\

\textbf{Definition 1} \, \emph{For any $H \in \cal{H}$, we call  a second order elliptic operator $L$ on $H
\setminus \Sigma$   \textbf{skin adapted}, supposed the following two conditions hold:}\\

 \textbf{$\mathbf{\bp}$-Adaptedness}  \, \emph{$L$ satisfies skin weighted uniformity
conditions relative to the charts $\psi_p$: {\small \[-L(u) = \sum_{i,j}  a_{ij} \cdot \frac{\p^2 u}{\p x_i \p x_j} + \sum_i b_i \cdot \frac{\p u}{\p x_i} + c \cdot u,\]}  \emph{for some locally} $\beta$-H\"{o}lder continuous coefficients $a_{ij}$, $\beta \in (0,1]$, measurable functions $b_i$, $c$, and some $k \ge 1$, so that for any $x \in H\setminus \Sigma$:}
{\small \[k^{-1} \cdot\sum_i \xi_i^2 \le \sum_{i,j} a_{ij} \cdot \xi_i \xi_j \le k \cdot \sum_i \xi_i^2, \,\, \delta^{\beta}_{\bp} \cdot  |a_{ij}|_{C^\beta(B_{\theta(p)}(p))} \le k,
\delta_{\bp} \cdot b_i \le k \mm{ and }\delta^2_{\bp} \cdot c \le k.\]}

\textbf{$\mathbf{\bp}$-Weak Coercivity} \, \emph{There exists a positive supersolution $u$ of the equation $L \, f = 0$ so that} \[L \, u \ge \ve \cdot \bp^2 \cdot u,
\mm{ for  some } \ve  >0.\]

\textbf{Remark 2} \, The class of skin adapted operators is large, but it is only the \emph{Hardy inequality} (\ref{hh}) for $\bp$ that shows that many classical operators on $H \setminus \Sigma$ belong to this class, independent of the chosen of $\bp$, cf. Theorem 10 for  basic examples.  In turn, this coupling to relevant problems is the exclusive use of the Hardy inequality. It is not employed in the proof of the general Theorems 4 - 9, but only to study the more explicit operators of Theorems 10 - 12.\qed

Now we state the basic results for skin adapted operators on $H \setminus \Sigma$. The central machinery to derive these results are boundary Harnack inequalities on
$H \setminus \Sigma$ relative to $\Sigma$, viewed as a boundary.\\

 \textbf{Theorem 4}  \textbf{(Boundary Harnack Inequalities)}\label{mbhsq} {\itshape \, Let $H \in {\cal{H}}$ and $L$ any skin adapted operator on $H \setminus \Sigma$ be given. Then, for any couple of open subsets $U \subset \! \subset V \subset \widehat{H}$, with $U \cap \widehat{\Sigma} \n$, there is a constant $C(L,U,V) >1$, so that for any two  solutions $u,v >0$ of $L\, w = 0$ on $H \setminus \Sigma$, both  L-vanishing* along $V \cap \widehat{\Sigma}$:
\begin{equation}\label{fhpq} u(x)/v(x) \le C \cdot  u(y)/v(y), \mm{ for any two points } x, y \in U \setminus \widehat{\Sigma}.
\end{equation}}
(*A solution $u >0$ \emph{L}-vanishes in some point  $p \in \widehat{\Sigma}$ if there is a supersolution $w >0$,   such that $u/w(x) \ra 0$, for $x \ra p$, $x \in H \setminus \Sigma$. This can also be paraphrased as a \emph{minimal growth} condition. In the case, where constant positive functions are supersolutions, \emph{L}-vanishing matches the classical meaning that $u\ra 0$.) \\

We prove Theorem 4 along with stronger versions of such inequalities. Their formulation requires a more technical language. Therefore, we only stated them in Ch.\ref{agrm}. \\

We derive these boundary Harnack inequalities on $H \setminus \Sigma$, via Gromov hyperbolizations, as a characteristic implication of the skin uniformity of $H \setminus \Sigma$. Indeed, for Euclidean domains $D \subset \R^n$ one knows that uniformity is largely equivalent to the validity of  boundary Harnack principles, for the Laplacian $\Delta$, relative $\p D$, cf.[Ai3].\\

These boundary Harnack inequalities allow us to characterize what can be thought as prime elements: the minimal solutions. Here we call $u >0$ \emph{minimal} if for any other solution $v >0$,  $v \le u$, we have $v \equiv c \cdot u$, for some constant $c >0$.\\

The space of minimal solutions (normalized to $1$ in some basepoint) is the (minimal) Martin boundary. In view of the subtle structure of $H$ near and in $\Sigma$ it is an unexpected outcome that the Martin theory of these operators is easy to describe:\\

 \textbf{Theorem 5} \textbf{(Martin Theory  on  $\mathbf{H \setminus \Sigma}$)} {\itshape \,   Let $H \in {\cal{H}}$ be a singular hypersurface and $L$ some
 skin adapted operator on $H \setminus \Sigma$. Then, we have:
\begin{itemize}
    \item The identity map on $H \setminus \Sigma$ extends to a homeomorphism  between $\widehat{H}$ and the Martin compactification $\overline{H \setminus
        \Sigma}_M$.
         \item All Martin boundary points are extremal points:  $\p^0_M (H \setminus \Sigma,L) \equiv \p_M(H \setminus \Sigma,L)$.
\end{itemize} In particular, $\widehat{\Sigma}$ and the minimal Martin boundary $\p^0_M (H \setminus \Sigma,L)$ are homeomorphic.\\}

Thus, as a counterpart to  (\ref{herg}),  for any function $u >0$ on $H \setminus \Sigma$ we have: $u$ solves $L \, v = 0$ if and only if:  there is a (unique) finite Radon measure $\mu$ on $\widehat{\Sigma}$
so that
\begin{equation}\label{muu} u(x) = u_{\mu}(x) =\int_{\widehat{\Sigma}} k(x;y) \, d \mu(y).
\end{equation}
In this representation formula, $k(x;y)$ denotes the Martin kernel of $L$ on $H \setminus \Sigma$. It is, up to multiples, the unique positive solution of $L \, v = 0$  on $H \setminus \Sigma$ which  $L$-vanishes in all points of $\widehat{\Sigma}$ except for $y$. The $k(x;y)$ are just the minimal solutions.\\

 For details concerning these notions from Martin theory we refer to  Ch.\ref{grb}.\\

\textbf{Extension Results to $\mathbf{\Sigma}$} \quad The next three extension theorems are our counterparts of classical results for suitably regular e.g. uniform Euclidean domains equipped with the Laplacian, cf.[AG],[Ai2] or [JK], for reference.\\

We start with a Fatou theorem for skin adapted operators on $H \setminus \Sigma$. For this, we quantify the non-tangentiality by means of a non-tangential twisted cone or  pencil \[\P(z,\rho):= \{x \in H \setminus \Sigma \,|\, \delta_{\bp} (x) > \rho \cdot d_{g_H}(x,z)\}, \mm{ pointing to } z \in \Sigma.\]
The angle $arctan(\rho^{-1})$ can be thought as the aperture of $\P(z,\rho)$ relative $z$.\\

 \textbf{Theorem 6} \textbf{(Relative Fatou Theorem on  $\mathbf{H \setminus \Sigma}$)} \, {\itshape For any $H \in {\cal{H}}$ and any skin adapted operator $L$ on $H \setminus \Sigma$.
 For any two finite Radon measures $\mu$ and $\nu$ on $\widehat{\Sigma}$  we have: For $\nu$-almost any $z \in \widehat{\Sigma}$ and any fixed $\rho >0$: \[u_{\mu}/u_{\nu}(x) \ra d \mu/d \nu(z),\,\mm{ for } x \ra z, \mm{ with } x \in \P(z,\rho),\] where $u_{\mu}$ and $u_{\nu}$ are the solutions of $L \, v=0$ associated to $\mu$ and $\nu$, according to (\ref{muu}), and
$d \mu/d \nu$ is the Radon-Nikodym derivative of $\mu$ with respect to $\nu$.}\\

In the case of an open subset $A \subset \Sigma$ with $u_{\mu}(A)=u_{\nu}(A)=0$ the Fatou theorem does not provide any insight into the behavior of
$u_{\mu}/u_{\nu}(x)$, for $x \ra z \in A$. On the other hand, in this case $u_{\mu}$ and $u_{\nu}$ both L-vanish along $A$. This makes them amenable to an analysis
using boundary Harnack inequalities. Indeed, we get a versatile complementary result which equally applies tangentially.\\

 \textbf{Theorem 7} \textbf{(Continuous Extensions to $\mathbf{\Sigma}$)}  \,  {\itshape  For any $H \in {\cal{H}}$,  any skin adapted operator $L$ on $H \setminus \Sigma$
and any two solutions  $u,v >0$ of $L\, w = 0$ on $H \setminus \Sigma$ both  L-vanishing along  some common open set $A \subset\widehat{\Sigma}$.  The quotient $u/v$ on $H \setminus \Sigma$ admits a continuous extension to $(H \setminus \Sigma) \cup A$.}\\

Different from these quotients, the individual solutions usually diverge towards $\Sigma$.
In turn, if, in addition to skin adaptedness, we have some a priori control over certain solutions, we may also solve boundary value problems:\\

\textbf{Theorem 8} \textbf{(Dirichlet Problem for Skin Adapted Operators)}  \emph{\, For any $H \in {\cal{H}}$ and let $L$ be a skin adapted operator  on $H \setminus \Sigma$, so
that constant functions solve $L \, v=0$ and $G(x,p)\ra 0$,  for $x \ra \widehat{\Sigma}$, and given  $p \in H \setminus \Sigma$.\\}

\emph{ Then, for any continuous function $f$ on $\widehat{\Sigma}$, there is a uniquely determined continuous function $F$ on $H$ so that}
 \[ L \, F =0 \mm{ with } F|_{\widehat{\Sigma}} \equiv f.\]\

 \textbf{Symmetric Operators} \quad A frequently considered type of elliptic problems
 is the analysis of eigenvalue equations associated to linear elliptic operators. For this type of problems it becomes convenient to focus on symmetric operators.\\

 For a symmetric operator $L$ on $H \setminus \Sigma$, which is adapted relative $\bp$ according to Definition 1, the weak coercivity condition
asserting the existence of a positive supersolution $u$ of the equation $L \, f = 0$ with
\begin{equation}\label{h7}
 L \, u \ge \ve \cdot \bp^2 \cdot u, \mm{ for  some, } \ve  >0.
\end{equation}
        is \emph{equivalent} to validity of the Hardy inequality
{\small \begin{equation} \label{hadi0}  \int_H  f  \cdot  L f  \,  dA \, \ge \, \tau \cdot \int_H \bp^2\cdot f^2 dA,\end{equation}}for any smooth $f $ compactly supported in $H
\setminus \Sigma$ and some positive constant $\tau = \tau(L,\bp,H)>0$.\\

For a given symmetric operator $L$, adapted to $\bp$, there is a largest $\lambda^{\bp}_{L,H} \in [-\infty,+\infty)$  so that the Hardy inequality (\ref{h7}) is satisfied. It
can be viewed as an eigenvalue of the operator $\delta_{\bp}^2 \cdot L$, the generalized \emph{principal eigenvalue} of $\delta_{\bp}^2 \cdot L$ cf. [P], Ch.4 for the associated
spectral geometry on unbounded domains.\\

We notice that $L$ is skin adapted when $\lambda^{\bp}_{L,H} >0$. We use this simple observation to extend definition \ref{sao}:\\

\textbf{Definition 2}  \,  \emph{An $\bp$-adpated  symmetric operator $L$ on $H \setminus \Sigma$ is called \textbf{shifted skin adapted}, when  the
 principal eigenvalue of $\delta_{\bp}^2 \cdot L$ is finite, $\lambda^{\bp}_{L,H} >-\infty$.} \\

For these operators we can still employ the theory for properly skin adapted ones when we consider eigenvalue problems. Since many naturally defined operators can be shown to be shifted skin adapted, cf. Theorem 10 for sample cases, this substantially widens the scope of applications.\\

The role of this generalized notion is explained in the following trichotomy which resembles that well-known from the spectral theory of operators on unbounded Euclidean domains, cf.[P],Ch.4.\\

 \textbf{Theorem 9} \textbf{(Criticality)}  \,  \emph{For any singular $H \in \cal{H}$ and any shifted skin adapted operator $L$ on $H \setminus \Sigma$ with H\"older continuous coefficients we set
\[L_\lambda:= L - \lambda \cdot \bp^2 \cdot Id, \mm{ for }\lambda \in\R.\] Then we have the following trichotomy.
\begin{itemize}
    \item \emph{\textbf{Subcritical}} when $\lambda < \lambda^{\bp}_{L,H}$,  $L_\lambda$ is skin adapted. The minimal solutions of $L_\lambda \, v=0$  L-vanish in  all but the one point in $\widehat{\Sigma}$, which represents this solution as a Martin boundary point.
    \item \emph{\textbf{Critical}} when $\lambda = \lambda^{\bp}_{L,H}$,  there is an, up to multiples,
    \[\mm{ unique positive solution, the \textbf{ground state} } \phi  \mm{ of L} _{\lambda^{\bp}_{L,H}} \, \phi = 0\]
$\phi$  L-vanishes along $\widehat{\Sigma}$ and can be described as the limit of first Dirichlet eigenfunctions for the operator $\delta_{\bp}^2 \cdot L$ on a
        sequence of smoothly bounded domains $\overline{D_m} \subset D_{m+1} \subset H \setminus \Sigma$, $m \ge 0$, with $\bigcup_m D_m = H \setminus \Sigma$.
    \item \emph{\textbf{Supercritical}} when $\lambda > \lambda^{\bp}_{L,H}$, $L_\lambda \, u = 0$ has no positive solution.\\
\end{itemize}}

Due to the boundedness of geometry, $(H \setminus \Sigma, d_{\bp^*})$ is not only complete but also  \emph{stochastically complete} making the geometric analysis on $H \setminus \Sigma$ amenable to many stochastic analysis techniques. An application, in the context of Theorem 9, is that the ground state and subcritical solutions of minimal growth (towards parts of $\Sigma$), where the Martin integral does not carry information,  admit presentations in terms of \emph{Feynman-Kac formulas}, cf.[P],Ch.7.3 and [El],Ch.IX.\\

\textbf{Geometric Operators}\quad  Up to this point we have not considered explicit examples of (shifted) skin adapted operators. Indeed, the coupling of the per se abstract theory above to classical operators
is an important  step. As already mentioned, this is owing to the one $\bp$-axiom we have not used so far, the Hardy inequality (\ref{hh}).\\

From this, many operators we can extract from the Euler-Lagrange equations of natural variational integrals on $H$ can be shown to be (shifted) skin adapted. We consider some basic examples.\\

 \textbf{Theorem 10} \textbf{(Curvature Constraints)} {\itshape \,  For any singular $H \in \cal{H}$, $H^n \subset M^{n+1}$, we have
\begin{enumerate}
 \item If $scal_M \ge 0$, then the \textbf{conformal Laplacian} on $H$ \[L_H: = -\Delta_H +\frac{n-2}{4 (n-1)} \cdot scal_H\] is skin adapted. In general, $L_H$ is shifted skin adapted.
     \item More generally, let $S$ be any smooth function on $M$ and  $scal_M \ge S$, then the  \textbf{S-conformal Laplacian} on $H$
     \[L_{H,S}: = -\Delta_H +\frac{n-2}{4 (n-1)} \cdot (scal_H - S|_H)\] is skin adapted. Again, in general, $L_{H,S}$ is shifted skin adapted.
     \item The \textbf{Laplacian} $-\Delta_H$ is shifted skin adapted. When $H$ is compact, the principal eigenvalue $\lambda^{\bp}_{-\Delta,H}$ vanishes and the ground state is that of a constant function. In particular, $H \setminus \Sigma$ has the Liouville property saying that all bounded harmonic functions are constant.
     \item The \textbf{Jacobi field operator} $J_H=-\Delta_H - |A|^2-Ric_M(\nu,\nu)$ is shifted skin adapted with  principal eigenvalue $\ge 0$.
     \item The \textbf{base operator} $L_{\bot}:= -\Delta_H + |A|^2 = -\Delta_H +  scal_M - scal_H - 2 \cdot Ric_M(\nu,\nu)$ is always skin adapted.
\end{enumerate} $Ric_M(\nu,\nu)$ denotes the Ricci curvature of $M$ for a normal vector $\nu$ to $H$, $scal_H$ and $scal_M$  the scalar
curvature of $H$ and $M$.\\\\}

\textbf{Natural Schr\"odinger Operators}\quad The sample operators we considered in Theorem 10 have two properties in common. Firstly, they are \emph{naturally} associated to any $H \in {\cal{H}}$, that is, they commute with the convergence of sequences of area minimizers. In the case of area minimizing cones, these operators, their Green's function and the set of minimal solutions of $L \, u = 0$ are reproduced under composition with the scaling map
\[ S_\eta: C \ra C, \mm{ given by } x \mapsto \eta \cdot x, \mm{ for } \eta \in (0,\infty).\]
That is we consider the map $u \mapsto u \circ S_\eta$ and regauge the values of the resulting functions $u \circ S_\eta$ in some base point $p \in C \setminus \sigma$ to $1$.
In this fashion we define a scaling action $S^*_\eta$ on the Martin boundary.\\

Secondly, all operators mentioned in Theorem 10 are Schr\"odinger operators and, on area minimizers, their principal eigenvalue remains finite. We merge these properties into one concept, cf. the discussion in \ref{spli} - \ref{ivn}.  \\

\textbf{Definition 3}  \,  \emph{A natural and shifted skin adapted operator $L$ is called a \textbf{natural Schr\"odinger operator}, when $L(H)$ has the form \[L(H)(u)=
-\Delta_H \, u + V_H(x) \cdot \, u \mm{ on } H \setminus \Sigma_H,\] for  some H\"older continuous function $V_H(x)$, for any given $H \in \cal{H}$.}\\

In the case of an area minimizing cone $C \subset \cal{H}$ the naturality of $L$ means $V_C(t \cdot x)=t^{-2} \cdot V_C(x)$, for any $x \in C \setminus \sigma_C$, $t >0$. That is, we can write
$V(x)=r^{-2} \cdot V^\times(\omega)$, for $x=(\omega,r) \in S_C \setminus \sigma \times \R^{> 0} = C \setminus \sigma$. \\

In other words, for cones, there is a separation of variables for $L$. This also gives the two distinguished minimal solutions in the Martin boundary, corresponding to the tip and infinity, a useful product structure. In turn, this inductive asymptotic splitting effect makes these operators on a general area minimizers $H$ amenable to more explicit methods than general shifted skin adapted operators.\\

\textbf{Theorem 11} \textbf{(Separation of Variables)} \, \emph{Let $C$ be a singular area minimizing cone, $L$ a natural Schr\"odinger operator. Then we have for
 the  skin adapted operator $L_\lambda=L - \lambda \cdot \bp^2 \cdot Id$, for  $\lambda < \lambda^{\bp}_{L,C}$:
\begin{enumerate}
\item The scaling action $S^*_\eta$ on $\p_M (C,L_\lambda)$ has exactly two fixed points, namely, the two distinguished Martin boundary points $[0]$ and $[1] \in \p_M (C,L_\lambda) \cong [0,1]\times \Sigma_{\p B_1(0) \cap C}/\sim$. Viewed as functions $\Psi_-=[0]$ and
    $\Psi_+=[1]$  on $C \setminus \sigma$ we have
\[\Psi_\pm(\omega,r) = \psi(\omega) \cdot r^{\alpha_\pm},\]
for $(\omega,r) \in S_C \setminus \sigma \times \R$, with {\small $\alpha_\pm = - \frac{n-2}{2} \pm \sqrt{ \Big( \frac{n-2}{2} \Big)^2 + \mu_{C,L^\times_\lambda}}$}.
 \item $\mu_{C,L^\times_\lambda} > - (\frac{n-2}{2})^2$ is the non-weighted principal eigenvalue and  $\psi(\omega) >0$  the ground state of an associated natural
     Schr\"odinger operator $L_\lambda^\times$:
 \[L_\lambda^\times(v)(\omega)= - \Delta_{S_C} v(\omega) + \big(V^\times(\omega) - \lambda \cdot (\bp^\times)^2(\omega) \big)\cdot v(\omega),\]
    defined on $S_C \setminus \Sigma_{S_C}$, where $\bp(x)=r^{-1} \cdot \bp^\times(\omega)$, for $x=(\omega,r) \in C \setminus \sigma$. That is, we have
    $L_\lambda^\times \, \psi =  \mu_{C,L^\times_\lambda} \cdot  \psi.$\\
\end{enumerate}}

For applications in scalar curvature geometry we derive some more detailed results for the conformal Laplacian $L_H$ on $H$ respectively $L_C$  its tangent cones $C$:\\

\textbf{Theorem 12}  \textbf{(Conformal Laplacians)}    \,\emph{There are constants $\Lambda_n > \lambda_n >0$ depending only on $n$, so that for  $\lambda \in (0, \lambda_n]$ and any singular area minimizing cone $C$, $(L_C)_\lambda$ is skin adapted and for $\Psi_\pm(\omega,r) = \psi(\omega) \cdot r^{\alpha_\pm}$ we have the estimates:
\begin{itemize}
    \item  \quad  {\small $0 > \alpha_+ \ge - (1- \sqrt{3/4}) \cdot \frac{n-2}{2} > - \frac{n-2}{2} > - (1+ \sqrt{3/4}) \cdot \frac{n-2}{2} \ge \alpha_- > -(n-2).$}
      \item \quad $|\psi|_{L^1(S_C \setminus \Sigma_{S_C})} \le a_{n,\lambda} \cdot \inf_{\omega \in S_C \setminus \Sigma_{S_C}}  \psi(\omega)$,
\end{itemize}
for  some constant $a_{n,\lambda} >0$, depending only on $n,\lambda$}.\\

\bigskip
 \setcounter{section}{2}
\renewcommand{\thesubsection}{\thesection}
\subsection{Hyperbolicity of Area Minimizers}

\bigskip
In this chapter we study the skin metric $d_{\bp}$ and the quasi-conformal metric $k_{H \setminus \Sigma}$ on $H \setminus \Sigma$ which may be regarded as
generalizations of the quasi-conformal metric $k_D$ on uniform domains $D \subset \R^n$.\\

 We establish their hyperbolicity properties and determine the associated ideal boundaries. Also we discuss
the difference between these two types of uniformizations of $H \setminus \Sigma$.
\bigskip
\subsubsection{Basic Properties of Skin Metrics}\label{discrete}
\bigskip

We start with some basic definitions.

\begin{definition}   \label{qh} \,  For a locally compact, \textbf{non-complete}, locally complete,
rectifiably connected metric space  $X$ we set \[\p X:= \overline{X} \setminus X, \mm{ where } \overline{X} \mm{ is the metric completion of } X,\] and define the
\textbf{quasi-hyperbolic metric} $k_X$ on $X$   as {\small
\begin{equation} \label{qx0} k_X(x,y) := \inf \Bigl \{\int_\gamma 1/dist(\cdot,\p X) \, \, \Big| \, \gamma \subset  X  \mbox{ rectifiable arc joining }  x \mbox{ and } y \Bigr \},
\end{equation}}
for any two points $x,y\in X$.
\end{definition}

In [L1],5.6,  we have seen from the isoperimetric inequality that for any connected area minimizing hypersurface $H$, $H \setminus \Sigma$ is rectifiably connected.\\

Thus $X=H \setminus \Sigma$ fits into this framework and we can write $\p X = \overline{X} \setminus X = \Sigma$, for $X = H \subset M$, and define the quasi-hyperbolic metric $k_{H \setminus \Sigma}$.\\

This process only uses the intrinsic metric $g_H$ on $H$ which is induced from the embedding in $M$. We use skin transforms to also capture information from $A_H$. In
the definition of $k_{H \setminus \Sigma}$, we replace $1/dist(\cdot,\Sigma)$ for $1/\delta_{\bp} = \bp$. This new metric is a hybrid that does not only hold information
from the intrinsic metric $g_H$ but also from the second fundamental from $A_H$ of $H$ within $M$.

\begin{definition} \emph{\textbf{(Skin Metrics)}} \label{sg} \, For a given skin transform $\bp$ on $H$, the metric
\[d_{\bp}(x,y) := \inf \Bigl  \{\int_\gamma \bp \, \, \Big| \, \gamma   \subset  H
\setminus \Sigma\mbox{ rectifiable arc joining }  x \mbox{ and } y  \Bigr \}, \]
 is called the \textbf{skin metric} on $H$.
\end{definition}

\begin{remark} \, \label{rer}  In the definitions of $k_{H \setminus \Sigma}$ and $d_{\bp}$ we may allow the rectifiable arcs joining $x,y \in H \setminus \Sigma$ also to reach $\Sigma$.
This follows from the path connectedness of $H \setminus \Sigma$ and the inequality $\bp(x) \ge L  / dist(x,\Sigma)$, since $1/dist(x,\Sigma)$ considered along any path reaching $\Sigma$ becomes non-integrable. Thus the infimum remains unchanged when we merely allow rectifiable paths supported in $H \setminus \Sigma$. \qed
\end{remark}

Step by step, we will see that $d_{\bp}$ is a richer and more useful geometry on $H \setminus \Sigma$ than $k_{H \setminus \Sigma}$. It is this metric that
supports our study of the asymptotic analysis near $\Sigma$. We consider $k_{H \setminus \Sigma}$ primarily for a better understanding of $d_{\bp}$.\\

As a first substantial difference between $d_{\bp}$ and $k_{H \setminus \Sigma}$ we notice: the metric $d_{\bp}$ is well-defined regardless whether $H$ is a singular or  a
regular hypersurface, whereas, when $H$ is  regular, $k_{H \setminus \Sigma}$ becomes literally pointless, since there is no distance towards $\Sigma \v$ we could measure.
Thus, we set $k_{H \setminus \Sigma} \equiv 0$, when  $\Sigma \v$.

\begin{proposition} \emph{\textbf{(Basic Properties of $\mathbf{d_{\bp}}$)}} \label{sgp} \, Both $(H \setminus \Sigma, d_{\bp})$ and
$(H \setminus \Sigma, k_{H \setminus \Sigma})$ are \textbf{complete} geodesic metric spaces and we have: \begin{equation} \label{for} \, d_{\bp}(x,y)\ge L \cdot k_{H
\setminus \Sigma} (x,y)\ge  L/2 \cdot  log\Bigl(\bigl(1+\frac {d(x,y)}{dist(x, \Sigma_H)}\bigr)\cdot\bigl(1+\frac {d(x,y)}{dist(y, \Sigma_H)}\bigr)\Bigr).\end{equation}
for any two $x,y \in H \setminus \Sigma$. $L$ denotes the Lipschitz constant for the $\bp$-distance $\delta_{\bp}$ according to axiom (S4) for $\bp$.

\begin{enumerate}
    \item $(H \setminus \Sigma, d_{\bp})$ has \textbf{bounded geometry}.
    \item For a flat norm converging sequence of minimizing hypersurfaces $H_i \ra H$, the skin metrics $d_{\bp}$ converges compactly on smooth domains (via $\D$-maps, cf. [L1), Ch.2.1). \end{enumerate}
\end{proposition}

\begin{remark} \label{dbg} \, 1. The  bounded geometry condition we need here asserts, that for some fixed constant $l \ge 1$ and radius $\rho > 0$ the following holds:\\

For any point $p \in H \setminus \Sigma$ the ball  $B_{\rho}(p)$ in $(H, d_{\bp})$ there is some chart  to some open set $U_p \subset \R^n$
$\phi_{p}: B_{\rho}(p) \ra U_p$ which is $l$-bi-Lipschitz.  For our convenience, we usually assume that $0 \in U_p$ and $\phi_{p}(p) = 0,$ cf. \ref{high} for versions of higher regularity.\\

2. The fact that  $(H \setminus \Sigma, d_{\bp})$ has bounded geometry is an important prerequisite for many analytic arguments we discuss later on. This condition usually fails
for $k_{H \setminus \Sigma}$, since there is no upper bound for $|A|(x) \cdot dist(x, \Sigma)$, for $x \in H \setminus \Sigma$. That is, there is no positive radius $\varrho
>0$ that puts uniform constraints on the geometry or topology of all balls of radius $\varrho$ in $(H \setminus \Sigma, k_{H \setminus \Sigma})$.\\

3. The precise statement of the convergence in (ii) is a bit more technical and described in [L1],Ch.3.1 \qed
\end{remark}

 {\bf Proof of \ref{sgp}} \quad   We choose $x,y \in H \setminus \Sigma$ and  take any smooth curve $\gamma :[0,1]\rightarrow H \setminus \Sigma$ with $\gamma (0)=x$ and $\gamma
(1)=y$. We observe from $dist(x, \Sigma)+ d(x,\gamma (t)) \ge dist(\gamma (t), \Sigma)$ and $|\nabla \gamma(t)| \ge |\nabla d(x,\gamma(t))|$ that due to
  $\bp(x) \ge L  / dist(x,\Sigma)$.
{\small \begin{equation}\int\limits_\gamma \bp  \ge  L \cdot \int_\gamma \frac {1}{dist(\cdot,\Sigma)}  \ge \end{equation}
\[L \cdot \int\limits_0 ^1\frac {ds(t)}{dist(x, \Sigma)+
d(x,\gamma (t))} \ge L \cdot  log\left(1+\frac {d(x,y)}{dist(x, \Sigma_H)}\right).\]} Similarly, we get the inequalities for $\ge  L \cdot  log\left(1+\frac {d(x,y)}{dist(y,
\Sigma_H)}\right)$ and add them up {\small  $$2 \cdot\int\limits_\gamma \bp \ge 2 \cdot L \cdot \int_\gamma \frac {1}{dist(\cdot,\Sigma)} \ge L \cdot \left( log\bigl(1+\frac
{d(x,y)}{dist(x, \Sigma_H)}\bigr) + log\bigl(1+\frac {d(x,y)}{dist(y, \Sigma_H)}\bigr) \right).$$} This holds for all  connecting curves $\gamma$ and hence for the infima $k_{H
\setminus \Sigma} (x,y)$ and $d_{\bp}(x,y)$. The inequality also shows that both metrics are \emph{complete} on $H \setminus \Sigma$.\\

 The fact that both spaces are geodesic easily follows from the Lipschitz continuity of $1/dist(x, \Sigma_H)$ and $\delta_{\bp}$ either using Helly's selection principle, cf. [SG],4.5
 or the more detailed study of path integrals in [BHK],Ch.10, since the spaces are complete and for any two points we can always find a compact subset so that any arc linking
 the two given points with a given upper bounded length must stay in this bounded set.\\

Now we turn to the properties (i) and (ii) asserted for the skin metric.\\

For property (i),  we denote the exponential map of $s \cdot M$ in $p$, scaled by $s \ge 1$, by $\exp_p[s \cdot M]: T_p M \ra s \cdot M$. $T_p M$ carries the
 flat metric $g_{T_p M}$.  Since $H$ is compact we get, for standard results for the exponential map and the regularity theory of $H$:\\

For each positive $\eta >0$ there is some $\Lambda(\eta) \gg 1$ so that for any $p \in H \setminus \Sigma$, $\exp_p[\Lambda(\eta) \cdot \bp(p) \cdot H]$
 is a  a local diffeomorphism from $B_{10^3/L}(0)$ onto its image in $H$ with
\begin{equation}\label{10} |\exp^*_p[\Lambda(\eta) \cdot \bp(p) \cdot H ](\Lambda(\eta) \cdot \bp(p)^2 \cdot g_H) - g_{T_p H}|_{C^5(B_{10/L}(p))} \le
\eta,\end{equation} where the radius $10/L$ is measured relative $\Lambda(\eta) \cdot \bp(p) \cdot H$, cf. Step 2, in the proof of [L1],Prop.3.1 for the elementary details. Thus the radius is $10/(L \cdot \Lambda(\eta)) \ll 1/L$ relative $\bp(p)
\cdot H$.\\

 When we choose a sufficiently small $\eta >0$, the exponential map is $l$-bi-Lipschitz on $B_{10/(L \cdot \Lambda(\eta))}(p) \subset \bp(p) \cdot H$, for some
$l(\eta) \ra 1$ for $\eta \ra 0$. Note that, different from the $C^5$-norm, the Lipschitz constant is scaling invariant.\\

Property (ii) follows from the naturality of the $\bp$ and the fact that the convergence upgrades to compact $C^k$-convergence, for any $k \ge 0$, cf.[L1],Ch.2.1 and that $H_i,
H$ and also their regular portions are rectifiably connected, cf.[L1],Prop.4.5. \qed

In [L1],Ch.3.3, we have seen how some Whitney type smoothing process can be applied to any skin transform $\bp$. It generates a smooth quasi skin
transform $\bp^*$ on $H \setminus \Sigma$, so that
\begin{description}
    \item[(SW1)] \quad $c_1 \cdot \delta_{\bp}(x) \le \delta_{\bp^*}(x)  \le c_2 \cdot \delta_{\bp}(x)$
    \item[(SW2)] \quad  $|\p^\beta \delta_{\bp^*}  / \p x^\beta |(x) \le c_3(\beta) \cdot \delta_{\bp}^{1-|\beta|}(x)$,
\end{description}
for any $x \in H \setminus \Sigma$.
 $\beta$ denotes the usual multi-index for derivatives, with respect to normal coordinates around $x \in H \setminus \Sigma$,  and $c_i >0$, $i=1,2,$ are constants depending on $H$ and $\bp$, $c_3 > 0$
depends on $H$, $\bp$ and on $\beta$. We set
\[d_{\bp^*}(x,y) := \inf \Bigl  \{\int_\gamma \bp^* \, \, \Big| \, \gamma   \subset  H
\setminus \Sigma\mbox{ rectifiable arc joining }  x \mbox{ and } y  \Bigr \}, \]
and notice that $d_{\bp^*}$ is the distance metric for the smooth Riemannian manifold $(H \setminus \Sigma, (\bp^*)^2 \cdot g_H)$ and we get:

\begin{corollary}\label{smskc} \quad  $(H \setminus \Sigma, d_{\bp^*})$, viewed as $(H \setminus \Sigma, (\bp^*)^2 \cdot g_H)$,
is a complete Riemannian manifold with bounded geometry and it is quasi-isometric to $(H \setminus \Sigma, d_{\bp})$.
\end{corollary}
 {\bf Proof } \quad This is an immediate consequence of \ref{sgp}(i) and (SW1).
 \qed

\begin{remark} \label{high}\quad  Actually, $(H \setminus \Sigma, (\bp^*_\alpha)^2 \cdot g_H)$ has a bounded geometry of higher order of regularity, that is, in the
definition \ref{dbg}.1 we also get $C^k$-estimates in place of Lipschitz estimates using also (SW2). But, in this paper, we do not use this detail. \qed
\end{remark}

\bigskip

\subsubsection{Gromov hyperbolicity of $H \setminus \Sigma$}\label{gh}

\bigskip

Now we prove that both the quasi-hyperbolic metric $k_{H \setminus \Sigma}$ and the skin metric $d_{\bp}$ on $H \setminus \Sigma$ are Gromov hyperbolic.\\

These results for $k_{H \setminus \Sigma}$ respectively $d_{\bp}$ are owing to the uniformity respectively to the skin uniformity of $H \setminus \Sigma$ proved in [L1],Th.4.
We restate them as follows.

\begin{proposition} \emph{\textbf{(Uniformity Properties of  $\mathbf{H \setminus \Sigma}$)}} \label{intsuni}  \, For any hypersurface
$H \in {\cal{H}}^n$  we have
\begin{enumerate}
\item When $\Sigma \n$, $H \setminus \Sigma$  is a  \textbf{c-uniform space}, for some $c >0$. That is, any pair $p,q \in H \setminus \Sigma$ can be joined by a
    \textbf{c-uniform curve} in $H \setminus \Sigma$, i.e. a rectifiable path $\gamma_{p,q}: [a,b] \ra H \setminus \Sigma$, for some $a <b$, with $\gamma_{p,q}(a)=p$,
    $\gamma_{p,q}(b)=q$, so that for any $z \in \gamma_{p,q}$:
    \begin{itemize}
    \item \emph{\textbf{Quasi-Geodesic}:} \quad \quad  $l(\gamma)  \le c \cdot  d(p,q),$
    \item  \emph{\textbf{Twisted Double Cones}:}\quad $l_{min}(\gamma_{p,q}(z)) \le c \cdot dist(z,\Sigma)(z)$.
\end{itemize}
\item $H \setminus \Sigma$  is a  \textbf{c*-skin uniform space}, for some $c^* >0$. That is, any pair $p,q \in H \setminus \Sigma$ can be joined by a \textbf{c*-skin uniform
    curve} in $H \setminus \Sigma$, i.e. a rectifiable path $\gamma_{p,q}: [a,b] \ra H \setminus \Sigma$, for some $a <b$, with $\gamma_{p,q}(a)=p$, $\gamma_{p,q}(b)=q$,
    so that for any $z \in \gamma_{p,q}$:
    \begin{itemize}
    \item \emph{\textbf{Quasi-Geodesic}:} \quad \quad  $l(\gamma)  \le c^* \cdot  d(p,q),$
    \item  \emph{\textbf{Twisted Double Skin Cones}:}\quad $l_{min}(\gamma_{p,q}(z)) \le c^* \cdot \delta_{\bp}(z)$.
\end{itemize}
\end{enumerate}
In both cases, $l_{min}(\gamma_{p,q}(z))$ denotes the minimum of the lengths of the subarcs of $\gamma_{p,q}$ from $p$ to $z$ and from $q$ to $z$.\\
\end{proposition}

Next we recall basic concepts from asymptotic geometry.

\begin{definition}\emph{\textbf{(Gromov hyperbolicity)}}\label{grh} \quad   A metric space is \textbf{geodesic}, when any two points can be joined by a geodesic, i.e. a path which is an isometric embedding of an
interval.\\

A geodesic metric space is \textbf{Gromov hyperbolic} or $\mathbf{\delta}$\textbf{-hyperbolic,} if all its geodesic triangles are $\mathbf{\delta}$\textbf{-thin} for
some $\delta >0$. This means that each point on the edge of any geodesic triangle is within $\delta$-distance of one of the other two edges.\\

Moreover, a complete Gromov hyperbolic space $X$ is called \textbf{visual} or $\mathbf{\beta}$\textbf{-roughly  starlike},  for some $\beta  >  0$,  relative some chosen base point  $p \in X$,
if for any $x \in X$ there is a geodesic ray starting from $p$ whose distance to $x$ is at most $\beta$.\\
\end{definition}

The concept of Gromov hyperbolic spaces embraces  a broad range of spaces, from hyperbolic manifolds up to combinatorial objects like trees. It is designed to study the
asymptotic behavior near infinity.

\begin{example}  \textbf{(Uniformity and Hyperbolicity)} \label{grex} \,  We describe some examples E1-E2, and counterexamples C1-C3, of spaces with some
hyperbolicity properties related to our subject.
\begin{itemize}
    \item \textbf{E1}\, Compact Riemannian manifolds are always Gromov hyperbolic: we just choose $\delta = diameter$ and find that the manifold is ${\delta}$-hyperbolic.
        Similarly, we observe that, for any compact Riemannian manifold $M$, the product $M \times \R$ is again Gromov hyperbolic
    \item \textbf{E2}\, Euclidean uniform domains $D \subset \R^n$, and more generally,  uniform spaces $X$,  equipped with their  quasi-hyperbolic metric $k_X$ on $X$ are
        Gromov hyperbolic, cf. [GO], [BHK]. As for conformal hyperbolizations, in the style of the Riemann uniformization theorem, this is a universal example, since the uniformization theory of Bonk, Heinonen and Koskela in [BHK] says:\\
        There is a bijective  conformal  correspondence between  the  quasiisometry  classes  of  proper geodesic  \emph{roughly  starlike  Gromov  hyperbolic} spaces  and  the  quasisimilarity  classes  of  bounded locally  compact  \emph{uniform}  spaces.
    \end{itemize}
    We also note some classes of often considered spaces, which are not  Gromov hyperbolic. In all cases we easily find large and non-thin triangles.
\begin{itemize}
    \item \textbf{C1}\, Neither asymptotically flat spaces, like the Euclidean space nor products of non-compact complete Gromov hyperbolic spaces are Gromov hyperbolic.
   \item \textbf{C2}\,  Even manifolds with sectional curvature $\equiv -1$, e.g. $\Z^2$-coverings of genus $\ge 2$ Riemann surfaces, may not be Gromov hyperbolic.
   \item \textbf{C3}\,  For a compact Riemannian manifold $(M^n,g_M)$, the product space $\R^{\ge 0} \times M^n$ equipped with the warped product metric $g_{\R} + (1+
       a \cdot r)^2 \cdot g_M$, $a \ge 0$, is Gromov hyperbolic if but only if  $a=0$.
\end{itemize}
In other words, although hyperbolicity is commonly associated with fast growth of length and volumina, further spreading of an already Gromov hyperbolic space towards infinity can even destroy its hyperbolicity. \qed
\end{example}

Now we will see that $H \setminus \Sigma$ admits natural hyperbolic geometries with varying additional properties. We start with the quasi-hyperbolic metric.

\begin{proposition}\label{c22} \, For any $H \in \cal{H}$ we have\begin{itemize}
    \item $(H \setminus \Sigma, k_{H \setminus \Sigma})$ is a complete Gromov hyperbolic space.
    \item When $H$ is compact, $(H \setminus \Sigma, k_{H \setminus \Sigma})$ is roughly starlike.
\end{itemize}
\end{proposition}

\textbf{Proof}  \quad We already know from \ref{sgp}  that $(H \setminus \Sigma, k_{H \setminus \Sigma})$ is a complete geodesic metric space.   Since $H \setminus \Sigma$
is a uniform space we can take advantage of the theory for uniform spaces in [GO] and [BHK](3.6) mentioned in \ref{grex} above to see that
$(H \setminus \Sigma, k_{H \setminus \Sigma})$ is Gromov hyperbolic. Also, for any compact $H \in \cal{H}$, [BHK](3.6) shows that $(H \setminus \Sigma, k_{H \setminus \Sigma})$ is  roughly starlike. \qed

As a critical difference to the quasi-hyperbolic metric on uniform Euclidean domains we observe that $(H \setminus \Sigma, k_{H \setminus \Sigma})$ need not to have
bounded geometry. This makes it an analytically hardly usable space. For instance, we cannot expect uniform Harnack inequalities for elliptic problems.\\

The geometric source for this non-bounded geometry are quickly sharpening wrinkles in $H \setminus \Sigma$ while we approach $\Sigma$. These degenerating wrinkles
correspond to singular rays in approximating tangent cones and we observe that along such wrinkles $\bp$, grows much faster than $1/dist(\cdot, \Sigma)$.\\

In turn, these wrinkles are ironed out in  $(H \setminus \Sigma, d_{\bp})$ and we want to show that this metric is still Gromov hyperbolic. However, near $\Sigma$ the
metric $(H \setminus \Sigma, k_{H \setminus \Sigma})$ roughly resembles a warped product, with fiber $\Sigma$, over a product metric as in E1 of \ref{grex}.\\

Thus, in view of counterexamples like C3 in \ref{grex}, we must take into account that conformally deforming $H \setminus \Sigma$ using $\bp \ge (L \cdot dist(\cdot,
\Sigma))^{-1}$, with $\bp \gg (L \cdot dist(\cdot, \Sigma))^{-1}$ in wrinkled places, potentially damages the hyperbolicity of $k_{H \setminus
\Sigma}$, when we think of $d_{\bp}$ as a further spread version of $k_{H \setminus \Sigma}$.\\

It is the sharper skin uniformity of $H \setminus \Sigma$ that \emph{compensates} this adversarial spreading effect on the hyperbolicity properties of skin metrics.

\begin{proposition}\label{c23} \, $(H \setminus \Sigma, d_{\bp})$ and its smoothing
$(H \setminus \Sigma, d_{\bp^*})$ are both complete Gromov hyperbolic spaces of bounded geometry.
\end{proposition}

\textbf{Proof}  \quad The boundedness of the geometry, its completeness and the fact that these are geodesic spaces have been shown in \ref{sgp} and \ref{smskc}.\\

 As explained above, our proof for the hyperbolicity of  $(H \setminus \Sigma, d_{\bp})$ is based on the skin uniformity of  $H \setminus \Sigma$. Nevertheless, it is modelled on
 strategies to derive the
hyperbolicity of $(D,k_D)$ for uniform domains $D \subset \R^n$ due to Gehring and Osgood [GO] and Bonk, Heinonen, Koskela [BHK],Ch.2-3.\\

For the rest of this proof, we let $H \setminus \Sigma$ be $a$-skin uniform, for some fixed $a \ge 1$. To not inflate the exposition with insignificant constants we
 assume that the Lipschitz constant $ L_{\bp}$ of $\delta_{\bp}$ equals $1$,   that is, $L_{\bp} =1$:
  \begin{equation} \label{lip} |\delta_{\bp}(x)-\delta_{\bp}(y)| \le d(x,y),\mm{ for any } x,y \in H \setminus \Sigma\end{equation}

We divide the proof into three steps we formulate in the next three lemmas.
\begin{lemma} \emph{\textbf{(Relations between $\mathbf{d_H}$ and $\mathbf{d_{\bp}}$)}} \label{ue} \,For any two points $x,y \in (H \setminus \Sigma, d_{\bp})$ we have
     \begin{equation} \label{ri} d_{\bp}(x,y) \le 4 \cdot a^2 \cdot \log \left(1+d_H(x,y) \cdot
     \max \{\bp(x), \bp(y)\}  \right), \end{equation}
     and from the elementary inequality $\log(1+x) \le \sqrt{x}$ for $x \ge 0$
   \begin{equation} \label{risq} d_{\bp}(x,y) \le   4 \cdot a^2 \cdot \sqrt{d_H(x,y) \cdot
     \max \{\bp(x), \bp(y)\}.}  \end{equation}
     Conversely, we get
      \begin{equation} \label{ril}   \log (1+l_H(\gamma(x,y)) \cdot \max \{\bp(x), \bp(y)\})  \le d_{\bp}(x,y) \end{equation}
    where $l_H(\gamma(x,y))$ is the length measured in $(H \setminus \Sigma, g_H)$ of a geodesic arc $\gamma(x,y)$ relative $(H \setminus \Sigma, d_{\bp})$ which links $x,y \in H \setminus \Sigma$.
From this we get
     \begin{equation} \label{riw}   \log (1+d_H(x,y) \cdot \max \{\bp(x), \bp(y)\})  \le d_{\bp}(x,y) \end{equation}
     and then, from the Lipschitz inequality $|\delta_{\bp}(x)-\delta_{\bp}(y)| \le d(x,y)$, we have
     \begin{equation} \label{riww}  \big|\log(\delta_{\bp}(x)) - \log(\delta_{\bp}(y)) \big| \le d_{\bp}(x,y). \end{equation}
     \end{lemma}

     \textbf{Proof} \quad  Let $\gamma \subset H \setminus \Sigma$ be an $a$-skin uniform arc of length $L = l_H(\gamma)$ (measured relative $g_H$) joining two points $x,y \in H \setminus \Sigma$. Choose the midpoint
     $z \in \gamma$ with $\gamma = \gamma_1 \cup \gamma_2$, $\{z\} = \gamma_1 \cap \gamma_2$, $l_H(\gamma_1)=l_H(\gamma_2)$ for two subarcs
     $\gamma_i \subset \gamma$, with $x \in \gamma_1$, $y \in \gamma_2$.\\

Now we claim the following main inequalities: \begin{equation} \label{ri1} l_{\bp}(\gamma_1) \le 2 \cdot a \cdot \log(1+ L \cdot \bp(x)) \mm{ and similarly }  l_{\bp}(\gamma_2) \le 2 \cdot a \cdot \log(1+ L \cdot \bp(y))\end{equation}
     where $l_{\bp}$ denotes the length relative $d_{\bp}$.  We first use the Lipschitz estimate $(\ref{lip})$ for $\delta_{\bp}$
to see that in the case where $l_H(\gamma_1) = L/2 < \delta_{\bp}(x)$:
\begin{equation} \label{ri3} l_{\bp}(\gamma_1) = \int_{\gamma_1} 1/\delta_{\bp} \le \int_0^{L/2}1/(\delta_{\bp}(x) - s) \, ds\end{equation} Namely,
 we parameterize $\gamma_1$ by arc length. Then, when we leave $x$ at time $0$, (\ref{lip}) shows that $\delta_{\bp}(\gamma(s)) \ge \delta_{\bp}(x) - s$ and thus (\ref{ri3}).\\

 We use this estimate when either we actually know that $l_H(\gamma_1) = L/2 < \delta_{\bp}(x)$ (in case A below) or we apply it to a subarc of $\gamma_1$ that satisfies this condition (in case B below). Thus, we distinguish between two  cases:
 \[ \mm{\textbf{A.}} \;\; l_H(\gamma_1) \le \frac{a}{a+1} \cdot \delta_{\bp}(x) \;\; \mm{ and } \;\;  \mm{\textbf{B.}}\;\; l_H(\gamma_1) > \frac{a}{a+1} \cdot \delta_{\bp}(x).\]
\textbf{A.}  \, From (\ref{ri3}): $l_{\bp}(\gamma_1) \le \int_0^{L/2}1/(\delta_{\bp}(x) - s) \, ds = \log\left(\delta_{\bp}(x)/(\delta_{\bp}(x)-L/2)\right)$
\[= \log\left(1/(1-L/\left(2 \cdot  \delta_{\bp}(x)\right))\right) \le 2 \cdot a \cdot \log(1 + L/\delta_{\bp}(x)),\]
here we used the elementary inequality $\log(1/(1-x)) \le 2 \cdot k \cdot \log(1+2 \cdot x)$ for any $k \ge 1$, $x \in [0,k/(k+1)]$.\\

\noindent \textbf{B.}  \, The $a$-skin uniformity shows that for $t \le L/2$: $a/l_H(\gamma_1([0,t])) \ge 1/
\delta_{\bp}(\gamma_1(t)).$ We combine this inequality and (\ref{ri3}), for the subarc of $\gamma_1$ from $x$ to the point where the length reaches the value $\frac{a}{a+1} \cdot \delta_{\bp}(x)$, and get, since $a \ge 1${\small \[l_{\bp}(\gamma_1) \le \int_0^{\frac{a}{a+1} \cdot
\delta_{\bp}(x)}1/(\delta_{\bp}(x) - s) \; ds+  a \cdot \int_{\frac{a}{a+1} \cdot \delta_{\bp}(x)}^{L/2}1/s \; ds \]
\[=\log\Big(\frac{1}{1-\frac{a}{a+1}}\Big) + a \cdot \log\Big(\frac{(a+1) \cdot L}{2 \cdot a \cdot \delta_{\bp}(x)}\Big) \le \log(1+a)+a \cdot \log\Big(\frac{L}{\delta_{\bp}(x)}\Big)\]}
\[\le a \cdot \log 2 + a \cdot \log (1+L/\delta_{\bp}(x)) \le 2 \cdot a \cdot \log (1+L/\delta_{\bp}(x)).\]
where we used that $L > 2 \cdot \frac{a}{a+1} \cdot \delta_{\bp}(x) \ge \delta_{\bp}(x)$ and applied the elementary inequality
     \begin{equation} \label{ri2} \log(1+ k \cdot x) \le k \cdot \log(1+   x), \mm{ for } k \ge 1, x \ge 0.\end{equation}
Thus in both cases, A. and B., we have get the first inequality in (\ref{ri1}), the second one follows similarly.\\

Now the inequality (\ref{ri}) follows from the fact that  $l_H(\gamma) \le a \cdot d_H(x,y)$ and another application of (\ref{ri2}). As already noted above this also implies (\ref{risq}). \\

For (\ref{ril}),  we choose $x,y \in H \setminus \Sigma$ and  consider any rectifiable curve $\gamma$ in $(H \setminus \Sigma, d_{\bp})$,
$\gamma :[0,1]\rightarrow H \setminus \Sigma$ with $\gamma (0)=x$ and $\gamma (1)=y$. Then  the Lipschitz continuity of $\delta_{\bp}$,   $\delta_{\bp}(y)  \le
\delta_{\bp}(x) + d(x,y)$ gives
\[\delta_{\bp}(\gamma (t)) \le \delta_{\bp}(x)+ d(x,\gamma (t)) \mm{  and }\delta_{\bp}(\gamma (t)) \le \delta_{\bp}(x)+ l_H(x,\gamma (t))\]
where $l_H(x,\gamma (t))$ is the length of the subarc $\gamma([0,t])$ measured in $(H \setminus \Sigma, g_H)$. From this, we note the following inequalities
\[\frac {\delta_{\bp}(\gamma (t))}{\delta_{\bp}(x)} \le \frac {\delta_{\bp}(x)+d(x,\gamma (t))}{\delta_{\bp}(x)} \le \frac {\delta_{\bp}(x)+l_H(x,\gamma (t))}{\delta_{\bp}(x)}\]

In the case where $\gamma$ is a geodesic arc  in $(H \setminus \Sigma, d_{\bp})$ we find {\small \[\log\left(\frac {\delta_{\bp}(x)+l_H(\gamma(x,y))}{\delta_{\bp}(x)}\right)
= \int\limits_0 ^1\frac {ds(t)}{\delta_{\bp}(x)+ d(x,\gamma (t))} \le \int\limits_0 ^1\frac {ds(t)}{\delta_{\bp}(\gamma (t))} \le d_{\bp}(x,y).\]}

 When we exchange the roles of $x$ and $y$,  we get completely analogous results. From both sets of inequalities and $d_H(x,y) \le l_H(\gamma(x,y))$, we get  (\ref{ril}), (\ref{riw}) and
(\ref{riww}).\qed

\begin{lemma} \emph{\textbf{(Skin Uniformity of Geodesic Arcs)}} \label{hu} \,  Every geodesic arc $\gamma$ in $(H \setminus \Sigma, d_{\bp})$ is a
\textbf{c-skin uniform} arc in  $(H \setminus \Sigma, g_H)$, for some $c\ge a$, which is independent of $\gamma$.
\end{lemma}

\textbf{Proof} \quad We choose a geodesic arc $\gamma$ in $(H \setminus \Sigma, d_{\bp})$ between $x,y \in H \setminus \Sigma$ and verify the two conditions for skin
uniform curves.\\

\textbf{Twisted Double Cone Conditions} \quad We set $D:= \max_{z \in \gamma} \delta_{\bp}(z)$. Then we can find unique integers $N(x), N(y) \ge 0$ so that
\[D/2^{N(x)+1} < \delta_{\bp}(x) \le D/2^{N(x)} \mm{ and }D/2^{N(y)+1} < \delta_{\bp}(y) \le D/2^{N(y)}.\]
Now we chop $\gamma$ into subarcs. We choose cut points  $x_0,...x_{N(x)}$ and $y_0,...y_{N(y)} \in \gamma$ so that $x_i$ is the first point on $\gamma$ starting from
$x$ where $\delta_{\bp}(x_i) = D/2^i$. It exists since $\delta_{\bp}$ is continuous. Similarly we define the points $y_j$ starting from $y$. This defines geodesic arcs
$\gamma_x(i)$ between $x_i$ and $x_{i+1}$, $\gamma_y(j)$ between $y_j$ and $y_{j+1}$ and an arc $\gamma_0$ between $x_0$ and $y_0$.Thus we have for $\bp = 1/\delta_{\bp}$:
\[\bp(\gamma_x(i)) \subset  [2^i/D,\infty),
 \bp(\gamma_y(j)) \subset  [2^j/D,\infty) \mm{ and } \bp(\gamma_0) \subset  [1/D,\infty).\]
Thus, since $\gamma$ is a
geodesic, (\ref{risq}) of \ref{ue} shows
\begin{itemize}
    \item \quad $l_H(\gamma_0)/D \le  l_{\bp}(\gamma_0) \le 4 \cdot a^2 \cdot \sqrt{2 \cdot l_H(\gamma_0)/D}.$
    \item \quad $l_H(\gamma_x(i))/(D/2^{i}) \le  l_{\bp}(\gamma_x(i)) \le 4 \cdot a^2 \cdot \sqrt{l_H(\gamma_x(i))/(D/2^{i+1})}.$
\end{itemize}
and we also get the analogous estimates for $\gamma_y(j)$. From these inequalities we observe
  \begin{equation} \label{rr} l_H(\gamma_0)/D \le  32 \cdot a^4, \, l_H(\gamma_x(i))/(D/2^{i}) \le 32 \cdot a^4, \,
l_H(\gamma_y(j))/(D/2^{j})  \le  32 \cdot a^4
\end{equation}
and in turn this means
 \[l_{\bp}(\gamma_0) \le  32 \cdot a^4, \, l_{\bp}(\gamma_x(i)) \le 32 \cdot a^4, \mm{ and } l_{\bp}(\gamma_y(j))  \le  32 \cdot a^4.\]

We  use this to estimate $\bp$ on $\gamma_x(i)$ away from its endpoints. To this end, we recall that in the start- and endpoints $x_i$ and $x_{i+1}$ of $\gamma_x(i)$ we
have $\delta_{\bp}(x_i)=D/2^i$ and $\delta_{\bp}(x_{i+1})=D/2^{i+1}$, similarly for $\gamma_y(i)$. Let $z \in \gamma_x(i)$ then we have from \ref{ue} (\ref{riw})
\[\big|\log(D/2^i) - \log(\delta_{\bp}(z)) \big|\le d_{\bp}(x_i,z) \le l_{\bp}(\gamma_x(i)) \le 32 \cdot a^4\] and thus  $\exp(-32 \cdot a^4) \cdot D / 2^i \le \delta_{\bp}(z)$, and
thus using (\ref{rr})   \[l_{min}(\gamma(z))  \le \sum_{k \ge i} 32
\cdot  a^4 \cdot D \cdot 2^{-k}  \le 64 \cdot  a ^4 \cdot D / 2^i \le b(a) \cdot \delta_{\bp}(z), \mm{ for }\]\vspace{-0.7cm}
 \begin{equation} \label{bbb}
b(a):= 64 \cdot a^4 \cdot  \exp(32 \cdot a^4).
 \end{equation}\\

 \textbf{Quasi-Geodesics} \, On $\gamma$ we choose the two points $X$ and $Y$ so that for the curves $\gamma_X$ from $x$ to $X$ and
 $\gamma_Y$ from $y$ to $Y$: \[l_H(\gamma_X)=l_H(\gamma_Y)=d_H(x,y)/2\]
 Then  each of the curves reaches at most the midpoint of $\gamma$ and therefore: \begin{itemize}
    \item The length of the arc $l_H(\gamma(X,Y))$ between $X$ and $Y$ is given by \[l_H(\gamma(X,Y))= l_H(\gamma)-d_H(x,y).\]
    \item $l_H(\gamma_X) \le  b(a) \cdot  \delta_{\bp}(X)$ and $l_H(\gamma_Y) \le  b(a) \cdot \delta_{\bp}(Y)$, for $b(a)$ from (\ref{bbb}).
    \item The triangle inequality shows that $d_H(X,Y) \le 2 \cdot d_H(x,y)$.
 \end{itemize}
Now \ref{ue}  (\ref{riw}) gives
\[\log(1+ (l_H(\gamma) - d_H(x,y)) \cdot  \max\{\bp(X),\bp(Y)\}) \le\]
\[ d_{\bp}(X,Y) \le  4 \cdot a^2 \cdot \log(1+2 \cdot d_H(x,y)\cdot  \max\{\bp(X),\bp(Y)\}) \le \]
\[  4 \cdot a^2 \cdot \log(1+ 4\cdot b(a))\] Thus we have for
\begin{equation}\label{b2b}
b^*(a):= \exp(4 \cdot a^2 \cdot \log(1+ 4\cdot b(a)))-1:
\end{equation}
  \begin{equation} \label{ror} l_H(\gamma) \le
d_H(x,y) + b^*(a) \cdot \min\{\delta_{\bp}(X),\delta_{\bp}(Y)\} \end{equation} Now we distinguish two cases:
\[ \mathbf{A.} \;\min\{\delta_{\bp}(X),\delta_{\bp}(Y)\} \le (4 \cdot a^2+1) \cdot d_H(x,y)\; \mm{ and }\ \mathbf{B.} \mm{ otherwise }\]\
In case \textbf{A}, we combine this with (\ref{ror}) and get $l_H(\gamma) \le (1 + b^*(a) \cdot (4 \cdot a^2+1)) \cdot d_H(x,y)$.\\

In case \textbf{B}, we get, again using the Lipschitz condition for $\delta_{\bp}$: {\small \[(4 \cdot a^2 + 1) \cdot d_H(x,y) \le \min\{\delta_{\bp}(X),\delta_{\bp}(Y)\} \le
\min\{\delta_{\bp}(x),\delta_{\bp}(y)\} +  d_H(x,y)\]} This also means $d_H(x,y)/\min\{\delta_{\bp}(x),\delta_{\bp}(y)\} \le (4 \cdot a^2)^{-1}$ and, with $a \ge1$, we find from
$x \ge \log(1+x),$ for $x > 0$:
\[d_{\bp}(x,y) \le  4 \cdot a^2 \cdot \log \left(1+d_H(x,y) \cdot \max \{\bp(x), \bp(y)\}  \right) \le 1\]
Now \ref{ue} (\ref{ri}) and (\ref{ril}) show
\begin{equation} \label{lll}\log (1+ l_H(\gamma) \cdot \max \{\bp(x), \bp(y)\}) ) \le d_{\bp}(x,y) \le 1\end{equation}
Now we use that
\[x \le 2 \cdot \log(1+x), \mm{ for } x \in [0, 1], \log(2) >1/2\mm{ and } x \ge \log(1+x), \mm{ for } x > 0\]
to see from (\ref{lll}) and  \ref{ue} (\ref{ri}) and (\ref{ril})
\[l_H(\gamma) \cdot \max \{\bp(x), \bp(y)\}) \le 2 \cdot\log (1+ l_H(\gamma) \cdot \max \{\bp(x), \bp(y)\})) \le \]
\[2 \cdot  d_{\bp}(x,y) \le  8 \cdot a^2 \cdot \log \left(1+d_H(x,y) \cdot \max \{\bp(x), \bp(y)\}  \right) \le\]
\[8 \cdot a^2 \cdot  d_H(x,y) \cdot \max \{\bp(x), \bp(y)\}.\]
Thus in case B, we get \[l_H(\gamma) \le 8 \cdot a^2 \cdot  d_H(x,y),\] and to serve both cases, A and B, we can choose $c(a):= 1 + b^*(a) \cdot (4 \cdot a^2+1) +  8 \cdot
a^2$ and conclude that $\gamma$ is a $c$-skin uniform arc.
 \qed

\begin{lemma} \emph{\textbf{(Thinness of Geodesic Triangles)}} \label{tg} \, There is some $\delta(a) >0$, so that every geodesic triangle in
$(H \setminus \Sigma, d_{\bp})$ is $\delta$-thin.
\end{lemma}

\textbf{Proof} \quad  Let $x,y,z \in H \setminus \Sigma$ be the vertices of a geodesic triangle in $(H \setminus \Sigma, d_{\bp})$, and $[x,y],[y,z],[x,z]$ the three geodesic
edges which by \ref{hu} are $c$-skin uniform for some $c \ge a \ge 1$. We claim that there is some $\delta(a) >0$ so that
\[d_{\bp}(p,[y,z] \cup [x,z]) \le \delta \mm{ for any } p \in [x,y].\]
We may assume that $l([x,p]) \le l([y,p])$, then the two condition of skin uniformity say that
\begin{equation}\label{tt}
l([x,p]) \le c \cdot \delta_{\bp}(p) \mm{ and } l([x,y]) \le  c \cdot d(x,y.)
\end{equation}

Now we distinguish between two cases:
\[ \mathbf{A.} \;\; c \cdot l([x,z])  <  l([x,p]) \; \mm{ and }\ \mathbf{B.} \;\; c \cdot l([x,z])  \ge  l([x,p])\]\
In case \textbf{A}, there is some $q \in [y,z]$ with $l([q,z])  =  (2 \cdot c)^{-1} \cdot l([x,p]) \le l([y,q])$ and
\[d_H(p,q) \le l([x,p])+l([x,z])+l([q,z]) \le (1 + c^{-1} + c^{-1}/2) \cdot l([x,p]),\]
and the $c$-skin uniformity shows that $l([x,p])/2 = c \cdot l([q,z]) \le c^2 \cdot \delta_{\bp}(q)$. Also we have from (\ref{tt}): $l([x,p])/2 = c/2 \cdot \delta_{\bp}(p)$.
Thus we get
\begin{equation}\label{10}
d_H(p,q) \le 2 \cdot c^2 \cdot (1 + c^{-1} + c^{-1}/2) \cdot \min \{\delta_{\bp}(p), \delta_{\bp}(q)\}.
\end{equation}

 In case \textbf{B}, there is some $q \in [x,z]$ with $l([q,x])  =  (2 \cdot c)^{-1} \cdot l([x,p]) \le l([z,q])$ and
\[d_H(p,q) \le l([x,p])+l([x,q])) \le (1 + c^{-1}/2) \cdot l([x,p]),\]
and the $c$-skin uniformity shows that $l([x,p])/2 = c \cdot l([q,x]) \le c^2 \cdot \delta_{\bp}(q)$. Thus using again $l([x,p])/2 = c/2 \cdot \delta_{\bp}(p)$ we get
\begin{equation}\label{20}
d_H(p,q) \le 2 \cdot c^2 \cdot (1 +  c^{-1}/2) \cdot \min \{\delta_{\bp}(p), \delta_{\bp}(q)\}.
\end{equation}

Thus for both cases, we get $d_H(p,q) \le 2 \cdot c^2 \cdot (1 + c^{-1} + c^{-1}/2) \cdot \min \{\delta_{\bp}(p), \delta_{\bp}(q)\}$ and we insert this, using again \ref{ue}(\ref{ri}), into
\[d_{\bp}(p,[y,z] \cup [x,z]) \le  d_{\bp}(p,q) \le 4 \cdot a^2 \cdot \log \left(1+d_H(p,q) /\min \{\delta_{\bp}(p), \delta_{\bp}(q)\}  \right)\]
to see that
\begin{equation}\label{30}
d_{\bp}(p,[y,z] \cup [x,z]) \le 4 \cdot a^2 \cdot \log \left(1+  c(a) \cdot ( 2 \cdot c(a) + 3)\right)=:\delta(a).
\end{equation}
That is, every geodesic triangle in $(H \setminus \Sigma, d_{\bp})$ is $\delta(a)$-thin. \qed

This also concludes the proof for the hyperbolicity of $(H \setminus \Sigma, d_{\bp})$. \qed\\

\noindent For the smooth version $(H \setminus \Sigma, d_{\bp^*})$ we recall the inequalities [L1],Th.3
\[c^{-1}  \cdot \delta_{\bp_\alpha}(x) \le \delta_{\bp^*_\alpha}(x)  \le c \cdot \delta_{\bp_\alpha}(x),\]
for any $x \in H \setminus \Sigma$ and where $c \ge 1$ is a constant depending on $H$.\\

 They imply that $(H \setminus \Sigma, d_{\bp^*})$ and $(H \setminus \Sigma, d_{\bp})$ are quasi-isometric. Now one uses the basic fact that Gromov hyperbolicity
 remains invariant under quasi-isometries, cf. [BH], Ch.III.H (1.9), to infer that $(H \setminus \Sigma, d_{\bp^*})$ is also Gromov hyperbolic.\qed

\begin{remark} \textbf{(Skin Metrics on Regular Spaces)}  \label{tei} \,  The metric $d_{\bp}$ is well-defined and non-trivial also when $H$ is  a regular hypersurface. By contrast, $k_{H \setminus \Sigma}$ looses its meaning when $H$ is regular,
since the best we could do, is to set the distance
towards $\Sigma \v$ equal to $+\infty$, and thus for any smooth $H$, $k_{H \setminus \Sigma} \equiv 0$.\\

For $\Sigma_H \v$, the metric space $(H,d_{\bp})$  is homeomorphic to $(H,g_H)$, unless $H$ is totally geodesic where it degenerates to a single point, whereas $(H,k_{H \setminus \Sigma})$ always is the one-point space.\\

We also notice, $d_{\bp}$ remains invariant under scalings of the original metric on $H$. The diameter $diam(H, d_{\bp})$ is a measure for the relative curvedness of $H$. For
instance, it is
$0$ iff $H$ is totally geodesic. For compact $H$,  $diam(H, d_{\bp})= \infty$ iff $H$ is singular.\\

The naturality of $\bp$ shows that the metrics $d_{\bp}$ on singular and on regular hypersurfaces match seamlessly:\\

 For a sequence of compact area minimizing hypersurfaces $H_i$ converging in flat norm to a compact minimizer $H_\infty$,
 $\bp_{H_i}$ converges in $C^\alpha$-norm to $\bp_{H_\infty}$, for $\alpha \in (0,1)$, around any given regular
point of $H_\infty$ and hence
\[(H_i \setminus \Sigma_{H_i},d_{\bp_{H_i}}) \ra (H_\infty \setminus \Sigma_{H_\infty}, d_{\bp_{H_\infty}}), \mm{ for } i \ra \infty,\]
 compactly on $H_\infty \setminus \Sigma_{H_\infty} \times H_\infty \setminus \Sigma_{H_\infty}$ via $\D$-maps, cf. [L1], Ch.2.3.
 The various skin metrics in a given converging sequence are not only individually Gromov hyperbolic, but they are $\delta$\emph{-hyperbolic, for the same} $\delta >0$.\\

 To see this, we note that the $\delta$ depends on the skin uniformity parameter. But the only not smoothly approximated portion of  $H_\infty$ is the region near $\Sigma_{H_\infty}$ and we know from [L1],Th.6 and the proof of [L1],Th.5, that there is
 a uniform skin uniformity parameter for all hypersurfaces in $\R^{n+1}$. The parameter for the $H_i$ can be estimated in terms of this uniform parameter and the constant for portions away of $\Sigma_i$ and this means in terms of  $H_\infty$.\\

 Also $H \setminus \Sigma$ has natural inner approximations by skin uniform domains cf.Ch.\ref{bhua} below, and [L1],Ch.4.4 for details. In Ch.\ref{bhua} we will see that
 these skin uniform domains also admit hyperbolic unfoldings and they compactly converge to the skin metric on $H \setminus \Sigma$ while the domain exhausts $H \setminus \Sigma$.\\

This bears similarities not only with the uniformization but also with moduli spaces for Riemann surfaces of genus $\ge 2$. The hyperbolic metrics on degenerating families of smooth surfaces will develop infinite complete ends in those places where the limit surface carries singular points. And in the smooth places the hyperbolic metrics converge smoothly towards the limit metric.\qed
\end{remark}

\begin{remark} \label{teti} \quad  Complementary to the fact that $(H \setminus \Sigma, k_{H \setminus \Sigma})$ may not have bounded geometry, we observe that
$(H \setminus \Sigma, d_{\bp})$ need not to be roughly starlike.\qed \end{remark}

\subsubsection{$\Sigma \subset H$ as a Gromov Boundary}\label{grch}

\bigskip

We use the hyperbolicity of $(H \setminus \Sigma,d_{\bp})$ and $(H \setminus \Sigma, k_{H \setminus \Sigma})$ to retrieve $\Sigma$ as ideal boundary points for
some particular compactifications of these spaces.\\

\textbf{Basic Conepts} \quad We start with a purely geometric setup which leads to the concept of Gromov boundaries, cf. [BH],Ch.III.H and [KB] for detailed expositions.\\

A \textbf{generalized geodesic ray} $\gamma: I \ra X$ is an isometric embedding of the interval $I \subset \R$ into $X$, where either $I = [0,\infty)$, then $\gamma$ is a proper
geodesic ray, or $I = [0,R]$, for some $R \in (0,\infty)$. Then $\gamma$ is a geodesic arc. When we fix a base point $p \in X$ we can use the hyperbolicity to canonically
identify any $x \in X$ with the (properly) generalized ray $\gamma_x$ with endpoint $\gamma(R) = x$.\\

For the following discussion we extend the definition of such a ray to $I = [0,\infty]$ setting $\gamma(x) = \gamma(R)$, when $x \in [R,\infty]$.

\begin{definition} \emph{\textbf{(Gromov product and boundary)}} \,  For a complete Gromov hyperbolic space
$X$ we introduce the following concepts. \begin{itemize}
      \item For any triple of points $x,y,z \in X$ we set
\begin{equation}\label{grop}
(y \cdot z)_x:=1/2 \cdot (d(x, y) + d(x,z)-d(y,z))
\end{equation}
$(y \cdot z)_x$ called the \textbf{Gromov product} of $y$ and $z$ with respect to $x$.
\item The set $\p_G X$ of equivalence classes $[\gamma]$ of geodesic rays,  from a fixed base point $p \in X$, with two such rays being equivalent if they have finite Hausdorff distance. $\p_G X$ is called the \textbf{Gromov boundary} of $X$.
    \end{itemize}
 \end{definition}
($\p_G X$  does not depend on the choice of the base point $p$.)

\begin{remark} \label{grco} \quad 1. $(y \cdot z)_x$ is a measure of how long the two geodesic rays in $X$ from $x$ to $y$ and $z$ remain close together. A space is $\delta$-hyperbolic, if and only if for any four points $w,x,y,z \in X: \, (y \cdot z)_x \ge \min\{(y \cdot w)_x,(z \cdot w)_x\} - \delta.$, cf.[KB],Ch.2 and [BH], III.H.1.22.\\

2.  To get a topology on $\overline{X}_G = X \cup \p_G X$, we define the notion of a converging sequence: $x_n \in \overline{X}$ converges to $x \in \overline{X}$ if there exist
generalized rays $c_n$ with $c_n(0) = p$ and $c_n(\infty) = x_n$ subconverging (on compact sets) to a generalized ray $c$ with $c(0) = p$ and $c(\infty) = x$.\\

3. The canonical map $X \hookrightarrow \overline{X}_G$ is a homeomorphism onto its image, $\p_G X$ is closed and $\overline{X}_G$  is compact, cf.[BH],
Ch.III.H.(3.7). $\overline{X}_G$ is called the \textbf{Gromov compactification} of $X$. It is a metrizable space.\\

4. The Gromov product can be used to describe a \textbf{neighborhood basis} $U(p,r)$, $r \ge 0$, of any point $p \in \p_G X$ within $\p_G X$ respectively $\U(p,r)$, $r \ge 0$, within $X$. We choose a basepoint $q \in X$ and set:
\begin{itemize}
  \item $U(p,r):= \{ x \in \p_G X \,|\, \mm{ there are geodesic rays } \gamma_1, \gamma_2 \mm{ starting from } q$  with
  $[\gamma_1] = p, [\gamma_2] = x,$ so that $\liminf_{t \ra \infty} (\gamma_1(t),\gamma_2(t))_q \ge r\}$,
  \item $\U(p,r):= \{z \in X \,|\, \mm{ there is a geodesic ray } \gamma \mm{ starting from } q$ with
$[\gamma] = p$,   so that $\liminf_{t \ra \infty} (\gamma(t),z)_q \ge r\}$,
\end{itemize}
cf.[KB],Ch.2 and [BH], Ch.III.H.(3.6).\qed
\end{remark}

The following, in a sense, dual concept to express general hyperbolicity properties is due to Ancona [An3], Def.14.

\begin{definition}  \emph{\textbf{($\Phi $-chains)}}  \label{phih} \quad
Let $M$ be a complete metric space and $\Phi$ be an increasing function $\Phi: \R^{\ge 0} \ra \R^{> 0}$ with $\lim_{t \ra \infty}(\Phi(t)) = \infty$ and set $c_0=\Phi(0)$. Then
a $\mathbf{\Phi}$\textbf{-chain} is a sequence of open subsets $V_i \subset M$, $i\ge 1$ with $V_{i+1} \subset V_i$ and a sequence of basepoints $x_i \in \p V_i$  such that
\[c_0 \le d(x_i,x_{i+1}) \le c_0^{-1}\,\, \mm{ and } \,\, d(x,V_{i+1}) \ge \Phi(d(x,x_i)),\] for any $x \in \p V_i$. We say that the points $x_i $ are linked through a $\Phi $-chain.
\end{definition}

$\Phi$-chains can be applied to control positive harmonic functions along curves. To use them to get Harnack inequalities on an open set $V \subset \overline{V} \subset W$, where $W$ is another open set in $X$, we need to make sure that any two points $p \in \p V$ and $q \in \p W$ can be linked through such a $\Phi$-chain passing through the same basepoint to gauge the estimates. Thus,  we need to find general methods to construct $\Phi$-chains.

\begin{example} \label{phh} \quad  1. The classical case is that of balls $B_{1/2^i}(0) \cap U^n$, $i \ge 1$, where $U^n$ is the upper half space model of the hyperbolic space.
These balls constitute a $\Phi$-chain, when we view the
$B_{1/2^i}(0) \cap U^n$ as halfspaces in $\H^n$ with $\Phi(t)=a_n + b_n \cdot t$, for suitable $a_n, b_n > 0$, depending only on $n$, cf.[An6], 3.2.a.\\

2. The latter construction essentially uses the hyperbolicity of $\H^n$: it is not hard to see that there are no $\Phi $-chains on the Euclidean $\R^n$, cf.[An4], Rm.C,p.93.\\

3. For a general $\delta$-hyperbolic space $X$ we can mimic halfspaces of $\H^n$. We choose any geodesic ray $\gamma: \R^{\ge 0} \ra  X$ and define for any $t \ge 0$ the subset
 \begin{equation}\label{chai}
 W_t:=\{x \in X\,|\, dist(x,\gamma([0,t]))<dist(x,\gamma([t,a]))\}
 \end{equation}

An explicit construction in [BHK],8.9, shows that for any integer $m >1$, the $U_i:=W_{i \cdot  22  \cdot \delta}$, $m \ge i \ge 0$, form a $\Phi_\delta$-chain, for  $\Phi_\delta(t)=\max\{\min\{\delta,1/22 \cdot \delta\}, t-6 \cdot \delta\}$.\\

4. A modification of the previous $\Phi$-chain, more directly associated to the hyperbolicity properties, can be defined using the Gromov product, cf.[An6], 6.9. $p=\gamma(0)$,
\begin{equation}\label{ach}
U^*_i := \{x \in X \,|\, (x \cdot \gamma(4 \cdot  i \cdot \delta))_p \ge 4 \cdot i \cdot \delta - 2 \cdot \delta\}
\end{equation}
The $U^*_i$ form a  $\Phi$-chain, for $\Phi^*_\delta(t)=\max\{t-2(\delta+2),  \delta-2\}$, cf. end [An6], 6.9, p.17. \qed
\end{example}

We now use the $\Phi$-chains of examples 3 and 4 to establish properties of a distinguished neighborhood basis around any point $z \in \p_GX$. This reveals some additional similarities to the
properties of the concentric balls $B_{1/2^i}(0) \cap U^n$ of example 1.

\begin{proposition}\emph{\textbf{($\Phi$-Neighborhood Basis)}}\label{neb} \, For a complete $\delta$-hyperbolic space
$X$ and some given $z \in \p_GX$ we can find a neighborhood basis $\N_i=\N_i(z)$, $i \ge 1$, of $z \in \p_GX$ with $\N_{i+1} \subset  \overline{\N_{i+1}} \subset \N_i \subset \overline{X}_G$, in the Gromov compactification $\overline{X}_G$, so that:\\

There are points $p_i \in X$ with $B_{c_0}(p_i) \subset (\N_{i+1} \setminus \N_i) \cap X$, so that any two points $p \in \p\N_i$ and $q \in \p\N_{i+1}$ can be joined through
   a $\Phi$-chain which passes through $p_i$. We call such a family of sets $\N_i$, a $\Phi$-neighborhood basis of $z$.
\end{proposition}

\textbf{Proof}  \quad We choose a geodesic ray $\gamma: \R^{\ge 0} \ra  X$, with basepoint $p = \gamma(0)$ so that $z=[\gamma]$.  Then, for  $U_i$ of (\ref{chai}) and $U^*_i$ and (\ref{ach}),
 we define the two families
 \begin{equation}\label{nbb}
 V_i:= U_i \cup (\overline{U_i} \cap \p_GX) \, \mm{ and } V^*_i:= U^*_i \cup (\overline{U^*_i} \cap \p_GX), \mm{ for } i \ge 1.
 \end{equation}
 We claim that both of them can be used as a $\Phi$-neighborhood basis of $z$.\\

Firstly, the two $\Phi$-chains $U_i$ and $U^*_i$ are coarsely equivalent in sense that for any $i$, there are $i^*$ and $i^+>i$ so that: \, $U_i \supset U^*_{i^*} \supset U_{i^+}$.
This follows from the explicit form of $\Phi_\delta$ and $\Phi^*_\delta$, cf. [BHK],8.3 for details.  In turn it is easily follows from the definition of the Gromov product, cf. \ref{grco}, 4.,  that the $V^*_j= U^*_j \cup (\overline{U^*_j} \cap \p_GX)$  form a neighborhood basis. Also the further conditions are satisfied. They are verified  in [BHK],8.10. \qed \\

\textbf{Identification of $\mathbf{\p_GX}$} \, For the flat  model of a uniform domain $D \subset \R^n$, the Gromov boundary of the complete space $X=(D,k_D)$ is well-understood: there is a canonical bijection
between $\p_GX$ and $\p D$, that assigns to each geodesic ray in $X$ its end point in $\p D$, cf. [BHK](3.6). It is the uniformity that makes sense of this and ensures the bijectivity.\\

The counterparts for the three complete spaces $X_{\bp}:= (H \setminus \Sigma,d_{\bp})$, $X_{\bp^*}:= (H \setminus \Sigma,d_{\bp^*})$ and $X_{1/dist}:= (H \setminus
\Sigma, k_{H \setminus \Sigma})$ read as follows.

\begin{proposition}\label{xgb} \,  For $H \in {\cal{H}}^c_n$, the identity map on $H \setminus \Sigma$ extends to homeomorphisms between $H$ and the Gromov compactifications
of $X_{\bp}$, $X_{\bp^*}$ and $X_{1/dist}$:
\[H \cong\overline{X_{\bp}}_G \cong \overline{X_{\bp^*}}_G \cong \overline{X_{1/dist}}_G.\]
where $ \cong$ means homeomorphic. In particular, we have:
\[\Sigma \cong\p_GX_{\bp} \cong \p_GX_{\bp^*} \cong \p_GX_{1/dist}.\]
For these identifications, one assigns to (equivalence classes of) geodesic rays in $X_{\bp}$, $X_{\bp^*}$ and
$X_{1/dist}$, from some fixed base point, their end points in $\Sigma \subset H$.\\
\end{proposition}

\textbf{Proof}  \quad For the uniform space $H \setminus \Sigma$ equipped with its quasi-hyperbolic metric $k_{H \setminus \Sigma}$, that is for  $X_{1/dist}$, the result is
covered from the general theory of uniform spaces cf. [BHK](3.6), (3.12) and the definition of the topology for the Gromov compactification.\\

Now we turn to the space $X_{\bp}$ and consider its Gromov boundary. We start with the definition of a \textbf{canonical bijection} $\Psi_{\Sigma}:\p_GX_{\bp} \ra \Sigma$.\\

To this end let $\gamma:[0,L) \ra H \setminus \Sigma$, $L \in (0,\infty]$ be a proper geodesic ray relative $X_{\bp}$ initiating from some basepoint $p = \gamma(0) \in H
\setminus \Sigma$ which, now relative $(H \setminus \Sigma, g_H)$, has length $L$ and is parameterized by arc-length. \\

We know from \ref{hu} that $\gamma$ is a $c$-skin uniform curve, for some $c(H) >0$. Thus, since $H$ is compact, the quasi-geodesic condition for $\gamma$, in the guise of its
consequence $diam X_{\bp} < \infty$, imply that $L < \infty$ and we claim that for $t < L$, $t \ra L$: $\gamma(t) \ra x$, for some point $x \in \Sigma$.\\

Indeed, the definition of $[0,L]$ as the maximal interval of definition shows that there is a sequence $t_i \in (0,L)$ with $t_i \ra L$, for $i \ra \infty$, so that $\gamma(t_i) \ra x$,
for some $x \in \Sigma$. Then the quasi-geodesic condition implies that also for any other such sequence $s_i \in (0,L)$ with $s_i \ra L$ we get $\gamma(s_i) \ra x$.\\

Next consider two such geodesic rays  $\gamma[1]$ and  $\gamma[2]$ with end points $x[k] \in \Sigma$, having finite Hausdorff distance in $X_{\bp}$, that is, they define
the same point in $\p_GX_{\bp}$. Then we find sequences $t_i[k] \in (0,L(\gamma[k]))$ with $t_i[k]
\ra L(\gamma[k])$, $k=1,2$, so that $d_{\bp}(t_i[1],t_i[2]) \le c = const. < \infty$ and we note that $\bp(t_i[k]) \ra \infty$, for $i \ra \infty$.\\

Then we infer from (\ref{riw}), that is, $\log (1+d_H(x,y) \cdot \max \{\bp(x), \bp(y)\})  \le d_{\bp}(x,y)$ that $d_H(t_i[1],t_i[2]) \ra 0$, for $i \ra \infty$ and, hence,
$x[1]=x[2]$. Thus every representative of a point in $\p_GX_{\bp}$ has the same endpoint in $\Sigma$. This way we get a well-defined map $\Psi_{\Sigma}$
from $\p_GX_{\bp}$ to $\Sigma$. \\

We claim that $\Psi_{\Sigma}$ is a \textbf{bijective} map.\\

\textbf{Surjectivity of} $\Psi_{\Sigma}$ \quad Let $x \in \Sigma$, then we choose a basepoint $p \in H \setminus \Sigma$ and a sequence $x_i \in H \setminus \Sigma$ with $x_i
\ra x$, for $i \ra \infty$ and a sequence of geodesic arcs $\gamma_i$ from $p$ to $x_i$. Then, using the Arzel\`{a}-Ascoli theorem, we get a compactly converging subsequence
of the $\gamma_i$ with limit geodesic $\gamma$. From the previous argument we see that $\gamma$ links $p$ with some $y \in \Sigma$. The quasi-geodesic condition of the
$\gamma_i$ shows that $y =x$.\\

\textbf{Injectivity of} $\Psi_{\Sigma}$ \quad For geodesic rays  $\gamma[1]$ and  $\gamma[2]$ with end points $x[1]=x[2] \in \Sigma$, we choose any two sequences $t_i[k] \in
(0,L(\gamma[k]))$, $i \in \Z^{\ge 0}$, with
\begin{equation}\label{di}
t_i[k] \ra L(\gamma[k]), k=1,2\mm{ and } l_i:=L(\gamma[1])-t_i[1] = L(\gamma[2])-t_i[2]
\end{equation}

Due to the skin uniformity we infer that for large $i \gg 1$:  \[\bp(\gamma[k](t_i[k])) \le  2 \cdot c /(L(\gamma[k])-t_i[k]).\] In turn, the triangle inequality
shows that
\[d_{H}(\gamma[1](t_i[1]),\gamma[2](t_i[2])) \le (L(\gamma[1])-t_i[1]) + (L(\gamma[2])-t_i[2])=2 \cdot l_i.\]
Now we use (\ref{risq}) in \ref{ue} and (\ref{di}) to see from
\[d_{\bp}(\gamma[1](t_i[1]),\gamma[2](t_i[2]))  \le   4 \cdot a^2 \cdot \sqrt{2 \cdot l_i \cdot  2 \cdot c /l_i}  \le   8 \cdot a^2 \cdot \sqrt{c}\]
that $d_{\bp}(\gamma[1](t_i[1]),\gamma[2](t_i[2]))$ remains
bounded when $i
\ra \infty$.\\

From this we can infer, $\gamma[1]$ and  $\gamma[2]$ have finite Hausdorff distance in $X_{\bp}$ and thus they determine the same point in $\p_GX_{\bp}$.\\

Next we claim that $\Phi_H: \overline{X_{\bp}}_G  \ra H$, defined as $\Phi_H|_ {H \setminus \Sigma}=id_ {H \setminus \Sigma}$ and $\Phi_H|_
{\p_GX_{\bp}}=\Psi_{\Sigma}$,
 is a \textbf{homoemorphism} that the extends the identity may on $H \setminus \Sigma$. For this it is enough to show that it is continuous since
$\Psi_{\Sigma}$ is bijective and $\p_GX_{\bp}$ and $\Sigma$ are compact, and $\overline{X_{\bp}}_G$ is metrizable. But this follows along the same lines as the proof of the
surjectivity above and the definition
of the topology for the Gromov compactification.\\

Finally, for $X_{\bp^*}$, the claim follows from the result for $X_{\bp}$. Since $(H \setminus \Sigma, d_{\bp^*})$ and $(H \setminus \Sigma, d_{\bp})$ are quasi-isometric their
Gromov compactifications and thus their Gromov boundaries are homeomorphic,  cf.[BH], Ch.III.H (3.9). \qed

The case of complete Euclidean hypersurfaces $H \in {\cal{H}}^{\R}_n$ is largely identical. There is only one additional and, for our purposes, important consideration needed.\\

 Due to the non-compactness of any such hypersurface $H$ we need to topologically compactify $H$. We consider the one point compactification $\widehat{H}$ of $H$. That is, we adjoin one additional
point at infinity, denoted $\infty_H$. We show that the Gromov compactification of $H$ also grows by only this one point. This is not hard but also not purely formal,  since it
uses/shows that $(H \setminus \Sigma, d_{\bp})$ has precisely one end.

\begin{proposition}\label{xgbx} \quad   For any $H \in {\cal{H}}^{\R}_n$ the identity map on $H \setminus \Sigma$
 extends to homeomorphisms between $\widehat{H}$ and the Gromov compactifications of $X_{\bp}$,
$X_{\bp^*}$ and $X_{1/dist}$:
\[\widehat{H} \cong\overline{X_{\bp}}_G \cong \overline{X_{\bp^*}}_G \cong \overline{X_{1/dist}}_G.\]
In particular, we have: $\widehat{\Sigma} \cong\p_GX_{\bp} \cong \p_GX_{\bp^*} \cong \p_GX_{1/dist}.$
\end{proposition}

\textbf{Proof}  \quad  This is an extension of the argument for \ref{xgb} to the case where we also have skin uniform curves of \emph{infinite} length relative  $(H \setminus \Sigma, g_H)$.
Again, for the quasi-hyperbolic metric this is contained in [BHK](3.6), (3.12). Now we explain how this is accomplished in the skin metric case.\\

Let $\gamma:[0,L) \ra H \setminus \Sigma$, $L \in (0,\infty]$ be a proper geodesic ray relative $X_{\bp}$ initiating from some basepoint $p = \gamma(0) \in H \setminus \Sigma$
which has length $L$ and is parameterized by arc-length relative  $(H \setminus \Sigma, g_H)$. This time we have the two options: $L < \infty$ or $L=\infty$.\\

For $L < \infty$ we argue as in  \ref{xgb} and find that a homeomorphism $\Psi^*_{\Sigma}$ from $\p^*_GX_{\bp}$ to $\Sigma$, where $\p^*_G X \subset \p_G X$ is the set
of equivalence classes of geodesic rays with finite length relative  $(H \setminus \Sigma, g_H)$. Indeed, we show that in a given equivalence class of geodesic rays either all
representing arcs have finite or all have infinite length.\\

 To see this we note that geodesic arcs relative $X_{\bp}$ are rather special $c$-skin uniform arcs since each subarc is again $c$-skin uniform. From this we infer from the twisted double skin cone
 condition, which is part of  \ref{hu}, for any geodesic ray $\gamma(t)$ parameterized by arc-length, with infinite length, relative
 $(H \setminus \Sigma, g_H)$: \[t \le c \cdot  \delta_{\bp}(\gamma(t)), \mm{ for any } t >0.\]
Now for any other geodesic ray $\gamma^*$ with infinite length relative $(H \setminus \Sigma, g_H)$ inequality (\ref{risq}) in \ref{ue} says
 \[d_{\bp}(\gamma(t),\gamma^*(t)) \le   4 \cdot a^2 \cdot \sqrt{d_H(\gamma(t),\gamma^*(t)) /
     \min  \{\delta_{\bp}(\gamma(t)), \delta_{\bp}(\gamma^*(t))\}.}\]
Then the triangle inequality shows $d_H(\gamma(t),\gamma^*(t)) \le 2 \cdot t$ and thus we get the equivalence of  $\gamma$ and $\gamma^*$ from:
$d_{\bp}(\gamma(t),\gamma^*(t)) \le  8 \cdot a^2$, for any $t >0$.\\

In turn for a geodesic ray  $\gamma^*$ that determines the same point in the Gromov boundary as $\gamma$, (\ref{riww}) in \ref{ue}
\[\big|\log(\delta_{\bp}(\gamma(t))) - \log(\delta_{\bp}(\gamma^*(t))) \big| \le d_{\bp}(\gamma(t),\gamma^*(t)),\] shows that
$\gamma^*$ has infinite length relative $(H \setminus \Sigma, g_H)$.\\

Consequently, there is precisely one point $z_\infty$ in the Gromov boundary  $\p_GX_{\bp}$ that corresponds to geodesic arcs with  infinite length relative $(H \setminus
\Sigma, g_H)$. Also any of these geodesic rays leaves any bounded set in $H$ since otherwise it would approach some $z \in \Sigma$, but these points are reached by rays of finite
length. Thus they all approach $\infty_H$ and we may identify $z_\infty$ with $\infty_H$.\\

 This way we extend the homeomorphism $\Psi^*_{\Sigma}$  from $\p^*_GX_{\bp}$ to $\Sigma$ to a homeomorphism $\Psi_{\widehat{\Sigma}}$  from
$\p_GX_{\bp}$ to $\widehat{\Sigma}$. The remaining assertions follow as in \ref{xgb} above. \qed\\

\textbf{Gromov Boundary of Cones} \quad Area minimizing cones play a central role as the typical blow-up geometries. Here we collect some consequences of their
homogeneity as scaling invariant spaces. We start with some readily checked observations and write $S_C:=\p B_1(0) \cap C \subset S^n$ for any $C \in {\cal{C}}_n$.

\begin{lemma} \emph{\textbf{(Basic Geometry of $(C \setminus \sigma, d_{\bp})$)}} \label{pro} \, The space $(C \setminus \sigma, d_{\bp})$ is canonically diffeomorphic to that of a cylinder $\R \times S_C$ equipped not with a warped product metric:
\[(\R \times S_C, \bp_C^2(1, x) \cdot g_{\R} + \bp_C^2(1, x) \cdot g_{S_C})\]
For each $x \in \R$, the horizontal space $\{x\}\times S_C \subset \R \times S_C$ is totally geodesic. However, the vertical lines defined by $\R \times \{p\}$, for $p \in S_C $, are usually not geodesic.\\

 Geodesics in this warped product are mapped onto geodesics under the projection on the second factor $\pi_2: \R \times S_C \ra  S_C$.\\
\end{lemma}

Now we have the following reformulation of \ref{xgbx} for  $(C \setminus \sigma, d_{\bp})$

\begin{corollary}\emph{\textbf{(Gromov  boundary of Minimal Cones)}} \label{gbc} \,  Let $C^n \subset \R^{n+1}$ be an area minimizing cone with singular set $\sigma$.
Then the Gromov boundary of $(C \setminus \sigma, d_{\bp})$ is homeomorphic to \[ [0,1]\times \Sigma_{S_C}/\sim ,\] where  $\Sigma_{S_C}=\p B_1(0) \cap \sigma_{C}$ is
the singular set of $S_C$, and \[x\sim y \;\mm{ if both } x \mm{ and } y \mm{ belong to } \{0\} \times \Sigma_{S_C}\mm{ or to }\{1\} \times \Sigma_{S_C}.\] In the special case,
where $\sigma = \{0\}$ this reads: the Gromov boundary of $(C \setminus \sigma, d_{\bp})$ contains precisely two points. \qed
\end{corollary}

\subsubsection{Hyperbolicity of Skin Uniform Domains}\label{bhua}

\bigskip
Skin uniform domains $\su(a)$ are regularized extensions of the sets $H \setminus \I(a)$ defined as  $\mathfrak{B}$-hulls, cf.[L1], Ch.4.4. They are skin uniform substructures within $H \setminus \Sigma$ and assembled as union of subcollections of balls in some fixed QT-skin adapted cover $\cal{A}$ of $H \setminus \Sigma$ we also defined in [L1],Ch.4.4. We briefly recall their main properties.

\begin{proposition}\label{std} \emph{\textbf{(Skin Uniform Domains)}} \quad Let $H \in {\cal{H}}$ be a $c$-skin uniform space for some $c(H) >0$.
Then there are some $\iota(c),\kappa(c) > 1$ so that for sufficiency small $a >0$ the following holds:\\

There is a domain $\su= \su(a) \subset H \setminus \Sigma$ with  $\E(a) \subset \su(a) \subset
\E(\iota \cdot a)$, so that any two points $p,q \in \su$ can be linked by an arc $\gamma_{p,q} \subset \su$ with
\begin{itemize}
    \item \quad $l(\gamma)  \le \kappa \cdot d_H(p,q)$,
    \item \quad $l_{min}(\gamma(z)) \le \kappa \cdot \min\{L_{\bp} \cdot dist(z,\p \su),\delta_{\bp}(z)\}$, for any $z \in \gamma_{p,q}$.
\end{itemize}
$\su$ is the union of suitably selected balls all belonging to one fixed QT-skin adapted cover $\cal{A}$ of $H \setminus \Sigma$.\\
\end{proposition}

We recall from [L1],Ch.4.4 that $\mathfrak{B}$-hulls $\su(a)$ are hybrids. They have some skin uniformity properties and  the
properties of (inner) uniform domains in $\R^n$ and we define the hybrid density $\mathbf{d}(z):= \min\{L_{\bp} \cdot dist(z,\p \su),\delta_{\bp}(z)\}$

\begin{definition} \emph{\textbf{(Skin Metrics on $\su$)}} \label{sg} \, For a given skin transform $\bp$ on $H$,
we define the hybrid density \[\mathbf{d}(z):= \min\{L_{\bp} \cdot dist(z,\p \su),\delta_{\bp}(z)\}\]
and call the metric
\[d_{\bp,\su}(x,y) := \inf \Bigl  \{\int_\gamma 1/\mathbf{d}(\cdot) \, \, \Big| \, \gamma   \subset  \su \mbox{ rectifiable arc joining }  x \mbox{ and } y  \Bigr \}, \]
 the \textbf{skin metric} on $\su$.
\end{definition}

Different from the case of the entire space $H \setminus \Sigma$ we have some uniformly controlled bounded geometry on each $\su$ and also uniform controls for its boundary, even for $H \in  {\cal{H}}^{\R}_n$ where the domains $\su$ are always unbounded. Therefore the counterpart of \ref{c22}, \ref{c23} and \ref{xgb} reads as follows.

\begin{proposition}\emph{\textbf{(Geometry of $\mathfrak{B}$-hulls)}} \label{xgbiu} \, For any $H \in \cal{H}$, and $\su = \su(a) \subset H \setminus \Sigma$, for $a >0$,
 we have
\begin{enumerate}
    \item  $(\su, k_{\su})$ and  $(\su, d_{\bp,\su})$ are complete and they have bounded geometry.
    \item  $(\su, k_{\su})$ and  $(\su, d_{\bp,\su})$ are  Gromov hyperbolic and for the Gromov boundary $\p_G (\su, k_{\su})$   we have
    \[\p_G (\su, k_{\su}) \cong \p_G (\su, d_{\bp,\su}) \cong \widehat{\p \su}.\] \end{enumerate}
     The metric spaces $(\su(a), d_{\bp,\su})$ compactly converge  to $(H \setminus \Sigma, d_{\bp})$, for $a \ra 0$.
\end{proposition}

As explained in \ref{tei} this underlines that the concept of skin metrics on regular spaces matches seamlessly the skin metric  $(H \setminus \Sigma, d_{\bp})$.\\

 {\bf Proof} \quad  Again, for notational convenience, we assume that $L_{\bp}=1$. We only check the case of quasi-hyperbolic metrics. The skin metric case follows completely
 similar lines as the proof of \ref{c23}.  The compact convergence of  $(\su(a), d_{\bp,\su})$ to $(H \setminus \Sigma, d_{\bp})$ follows from $\delta_{\bp} \le L_{\bp}  \cdot dist(x,\Sigma)$.\\

 For (i). The completeness follows from \ref{sgp}. To check  the boundedness of the geometry $(\su, k_{\su})$, we define particular charts for $(\su, k_{\su})$.\\

We first note that $\overline{\su}  \subset (H \setminus \I(\alpha /(2\cdot j)))$. This means that in a neighborhood of $\overline{\su}$, we have $|A| \le \bp \le 2\cdot j/\alpha$.
Since the ambient space is
either a compact manifold or the Euclidean space this gives us uniform estimates for the geometry on balls in $\su$.\\

Concretely, we can find some $e_H \ge 1$, so that  the exponential map in $p \in H \setminus \Sigma$, for the original metrics $g_H$,  $\exp_p: T_p(H \setminus \Sigma)
\ra H \setminus \Sigma$, is a local diffeomorphism on  $B_{dist(p,\p
\su)/(8 \cdot e_H \cdot j)}(0),$ for $p \in \su$, with bi-Lipschitz constant $\le 2$. In the cone case, we can choose $e_H=e_n$ for some $e_n$ depending only on the dimension.\\

Now we consider the original balls $B^*(p):= B_{dist(p,\p \su)/(8 \cdot e_H \cdot j)}(p)$ equipped with $k_{\su}$. To control their  geometry we note from the triangle
inequality $|dist(p,\p \su) - dist(q,\p \su)| \le d(p,q),$ for any two points $p,q \in T_j:$
 \[1/2 \le dist(p,\p \su)/dist(q,\p \su)^{-1} \le 2, \mm{ for } q \in  B^*(p) \subset B_{dist(p,\p \su)/8}(p).\]

This shows that the identity map on $B^*(p)$ is a $K$-bi-Lipschitz map, for $K = 4$, when viewed as a map from the ball  scaled by the constant $dist(p,\p \su)^{-1}$ to the
ball we equipped with $k_{\su}$ written as $dist(x,\p \su)^{-2}\cdot g_H$
 {\small \[I_p: (B^*(p),dist(p,\p \su)^{-2}\cdot g_H) \ra (B^*(p), dist(x,\p
\su)^{-2}\cdot g_H),\mm{ for }x \in B^*(p)\]} Finally, we define a local chart $\phi_p$ as the inverse of the composition map  {\small \[I_p \circ \exp_p (dist(p,\p \su) \cdot x),
\mm{ for } x \in B_{1/(8 \cdot e_H \cdot j)}(0) \subset T_p(H \setminus \Sigma),\]} and observe that $\phi_p$ is $K^*$-bi-Lipschitz map, for $K^* = 8$.\\

For (ii). The assertions $\su$ is bounded $\Leftrightarrow H \in {\cal{H}}^c_n$ resp. $\su$ is unbounded $\Leftrightarrow H \in {\cal{H}}^{\R}_n$,  follow from the fact that
$\I_H^c(a) \subset \su(a)$.\\

 Since any $\su$ is an (inner) uniform domain we can take advantage of the theory for uniform spaces in [GO] and [BHK](3.6)  to see that
$(\su, k_{\su})$ is Gromov hyperbolic with Gromov boundaries equal to $\p \su$ for $\su$ bounded and $\p \su^*$  for $\su$ unbounded. \qed

\setcounter{section}{3}
\renewcommand{\thesubsection}{\thesection}
\subsection{Martin Theory on $H \setminus \Sigma$}

\bigskip

Now we turn to more analytic aspects of  $H \setminus \Sigma$ and develop a potential theory for elliptic operators  on $H \setminus \Sigma$ regarding $\Sigma$ as a boundary. \\

Besides the Gromov boundary and compactification there are other important concepts of ideal boundaries and associated compactifications related to further structural information, like analytic or algebraic structures on
the given spaces, cf.[BL].\\

The structures we are interested in are elliptic operators on $H \setminus \Sigma$. In this case, Martin compactifications are a natural choice and we briefly review the basic
 notions and Ancona's theory to characterize these compactifications in \ref{grb} and \ref{bhai} below.\\

In \ref{agrm} we derive the main results of this chapter. We define the class of skin adapted operators on $H \setminus
\Sigma$ and show that their Martin boundary is homeomorphic to $\Sigma$.

\subsubsection{Martin Compactifications}\label{grb}
\bigskip

We recall some basic notions from Martin theory, cf. [M], [P],7.1 or [BJ],Ch.I.7. for details.

\begin{definition}\emph{\textbf{(Martin Boundary)}}\label{mb} \quad For a non-compact Riemannian manifold $X$ and some linear
 second order elliptic operator $L$ on $X$, with some positive minimal Green's function $G: X\times X \ra (0,\infty]$, we choose some base point $p \in X$ and consider the space $S$ of
 sequences $s$
 of points $p_n \in X, n \ge 1$
 with \begin{itemize}
    \item $p_n \in X$ has no accumulation points in  $X$,
    \item $K(x,p_n):= G(x,p_n)/G(p,p_n) \ra K_s(x)$, in compact convergence, for $n \ra \infty$, for some function $K_s$ on $X$.
 \end{itemize}
 Then, the \textbf{Martin boundary}\, $\p_M (X,L)$ is the quotient of $S$ modulo the relation on $S$: $s \sim s^*$,   if $K_s \equiv K_{s^*}$. \\

\noindent Also we define the \textbf{Martin kernel} $k(x;y)$ on $X \times \p_M (X,L)$ by $k(x;y):= K_s(x)$, for some $K_s$ representing $y \in \p_M (X,L)$.
\end{definition}
(Again, the definitions do not depend on the choice of the base point $p$.)

\begin{remark} \label{mrtt} \quad 1. The Harnack inequality and elliptic theory show that each $K_s \in \p_M (X,L)$ is a positive solution of $L \, u = 0$ on $X$. This also shows
that the convex set $S_L(X)$ of positive solutions of $L u =0$ on $X$ with $u(p) =1$ is compact in the topology of compact convergence. In turn, $\p_M (X,L)$ is a compact
subset of $S_L(X)$. Being a \textbf{minimal Green's function} says that for any $x \in X$ \begin{itemize}
    \item  $G(\cdot,x)$ solves $L \v=0$ on $X \setminus \{x\}$,
    \item $L G(\cdot,x)= - \delta_x$, where $\delta_x$ denotes the Dirac measure in $x$,
    \item  There is no positive solution $w$ on $X$ with $w \le G(\cdot,x)$.
\end{itemize}
 The latter property is also paraphrased saying $G(\cdot,x)$ is an $L$\emph{-potential}.\\

2. The \textbf{Martin topology} on  $\overline{X}_M:= X \cup \p_M (X,L)$ is defined as follows: a sequence $s$, $p_n \in X$, with no accumulation points in  $X$ converges to a
point $y \in \p_M (X,L)$ iff $K(x,p_n)$ converges compactly to $K_s(x)$ representing $y$. And $y_n  \in  \p_M (X,L)$ converges to $y  \in \p_M$ iff there are functions
$K_s(y_n)(x)$ representing
$y_n$ compactly converging to a functions $K_s(y)(x)$  representing $y$.\\

The outcome is that $\p_M (X,L)$ is closed and $\overline{X}_M$ is compact. This is the \textbf{Martin compactification} of $(X,L)$. It is also easy to see that $\overline{X}_M$
is metrizable, cf.[BJ], Ch.I.7. or [H], Ch.12 for further details.  \qed
\end{remark}

To motivate the following ideas, recall the classical result, due to Minkowski, cf. [C],Ch.6, that each point in a convex set  $K \subset \R^n$ is a uniquely determined
convex combination of extremal points of $K$.\\

We consider the convex set  $S_L(X)$. The extremal elements of $S_L(X)$ form a particular subset $\p^0_M (X,L) \subset \p_M (X,L)$ of the Martin boundary.\\

It is easy to check that a positive solution $u$ of $L \, u = 0$ on $X$ with $u(p)=1$ is extremal iff $u$ is a minimal solution. Here we call $u$ \textbf{minimal}
if for any other solution $v >0$,  $v \le u$, we have $v \equiv c \cdot u$, for some constant $c
>0$. Therefore $\p^0_M (X,L) \subset \p_M(X,L)$ is also called the \textbf{minimal Martin boundary}.\\

The Choquet integral representations in [C], Ch.6, give the following general version of the Martin representation theorem, cf.[P],7.1 or [An3], Cor.13 :

\begin{proposition}\emph{\textbf{(Martin Integral)}} \label{mrt} \quad For any positive solution $u$ of $L \, u = 0$ on $X$, there is a
unique finite Radon measure $\mu_{u}$ on $\p^0_M (X,L)$ so that
\[u(x)  =\int_{\p^0_M (X,L)} k(x;y) \, d \mu_u(y).\]
Conversely, for any finite Radon measure $\mu$ on $\p^0_M (X,L)$,
\[h(x)  =\int_{\p^0_M (X,L)} k(x;y) \, d \mu(y),\]
defines a positive solution $h$ of $L \, h = 0$ on $X$.
\end{proposition}

\bigskip
\subsubsection{Boundary Harnack Inequalities}\label{bhai}
\bigskip

In general, the Martin boundary for a given pair $(X,L)$ is hard to determine. However, this changes when we can establish \emph{boundary Harnack inequalities}, for positive solutions of $L\, u =0$, either along the boundary of locally bounded Euclidean domains or along ideal boundaries at infinity.\\

To explain this type of inequalities we first look at the classical situation of a Euclidean domain $D \subset \R^n$. We use the common abbreviation of the assertion about the \emph{validity of the boundary Harnack principle} as BHP:

\begin{proposition}\emph{\textbf{(Euclidean BHP)}} \label{abhp} \,  For any sufficiently* regular domain $D \subset \R^n$, *cf. \ref{rec1} below, we have the following two versions of a BHP.
\begin{itemize}
  \item \emph{\textbf{(Global BHP)}} For any pair of a compact set $K$ and an open $V$ with $K \subset V \subset \R^n$, $K \cap D$ and $K \cap \p D$, there is a constant $C$ depending only on $V,K$  so that for any two harmonic functions $u,v >0$  on $V \cap D$ vanishing along $V \cap \p D$:
\begin{equation}\label{gf} u(x)/v(x) \le C \cdot  u(y)/v(y), \mm{ for any two points } x, y \in K \cap D.
\end{equation}
  \item \emph{\textbf{(Uniform BHP)}}  There are constants $A,C >1$ depending only on $D$, so that for any point $p \in \p D$ and $R>0$ small enough the following holds: for any two harmonic functions $u,v >0$  on $B_{A \cdot R}(p) \cap D$ vanishing along $B_{A \cdot R}(p) \cap \p D$:
\begin{equation}\label{gf2} u(x)/v(x) \le C \cdot  u(y)/v(y), \mm{ for any two points } x, y \in B_R(p) \cap D.
\end{equation}
\end{itemize}
\end{proposition}

\begin{remark} \label{rec1} \textbf{(Regularity of Domains)}\quad  The uniform BHP implies the global BHP. And we note that it is the uniform BHP that is used in Martin theory.\\

To explain the attribute \emph{sufficiently regular} we note that, in consecutive work, starting in the early 20th century (with results like the Phragm\'{e}n-Lindel\"of principle), one became aware of the importance of the shape of the boundary in the control of harmonic functions. Gradually, it became possible to establish the BHP for increasingly general types of domains.\\

 A particularly important and quite general case that even includes some fractal domains, is that of \emph{uniform} Euclidean domains $D \subset \R^n$. Aikawa [Ai2] has shown that uniform domains  satisfy the BHP for the Laplacian $\Delta$. Remarkably, this is an optimal result and characterizes uniform domains, cf.[Ai3]: under some mild additional hypotheses, concerning the exterior
capacity density, uniformity of a Euclidean domain is \emph{equivalent} to the validity of the BHP for $\Delta$. \qed
\end{remark}

Now we turn to the counterpart of this analysis on a complete Gromov hyperbolic manifold. We recall from the introduction that this is classically suggested from the
two ways to derive and interpret the integral representations of harmonic functions through the Poisson integral.\\

In work, in particular, due to Kifer [K], Anderson [A], Sullivan [S] and  Anderson and Schoen [AS]  the authors succeeded to understand the space of
positive harmonic functions on (simply connected) manifolds of pinched negative sectional curvature. Inspired from [AS], Ancona developed
a broad theory to understand potential theory on spaces with some rather general hyperbolicity properties, cf.[An3], [An4].\\

 It is Ancona's theory that
 matches the degree of abstraction and generality we need here. To state  Ancona's BHP
we first describe the range of admissible operators:

\begin{definition} \label{sao0} \quad  For a complete Riemannian manifold $X$ with bounded geometry, we call a second order elliptic operator
$L$ on $X$   \textbf{adapted weakly coercive} provided the following conditions hold:
\begin{itemize}
    \item $L$ is \textbf{adapted}: $L$ is a linear second order elliptic operator so that relative to the charts $\phi_{p}$: $-L(u) =
         \sum_{i,j} a_{ij} \cdot \frac{\p^2 u}{\p x_i \p x_j} + \sum_i b_i \cdot \frac{\p u}{\p x_i} + c \cdot u,$\\
 for $\beta$-H\"{o}lder continuous $a_{ij}$, $\beta \in (0,1]$, measurable functions $b_i$, $c$, with \small{\[k^{-1} \cdot\sum_i \xi_i^2 \le \sum_{i,j} a_{ij} \cdot \xi_i \xi_j
\le k \cdot\sum_i \xi_i^2 \mm{ and } |a_{ij}|_{C^\beta(B_{\rho}(p))}, |b_i|_{L^\infty}, |c|_{L^\infty} \le k,\]}  for some $k \ge 1$ and any $p \in X$.
    \item  $L$ is \textbf{weakly coercive}: There is a positive supersolution $u$ of the equation $L \, f = 0$ with the quantitative estimate \, $L \, u \ge \ve \cdot u ,
        \mm{ for some } \ve  >0.$\\
\end{itemize}
\end{definition}

The associated BHP is formulated relative $\p_GX$. The boundary condition is that of  \textbf{\emph{L}-vanishing}: a solution $u \ge 0$ \emph{L}-vanishes along some open subset $V \in \p_GX$, if there is a supersolution $w >0$  such that $u/w \ra 0$ when we approach $V$ in $X$.\\

It is easy to see that the minimal Green's functions $G(\cdot,p)$ for any basepoint $p \in X$ $L$-vanishes along  $\p_GX$, cf.[An3], lower part of p.509 for an explicit construction.\\

The counterpart to the uniform Euclidean BHP can be stated by means of a $\Phi$-neighborhood basis $\N_i(z)$ of a given $z \in \p_GX$ in $\overline{X}_G$ we described in \ref{neb}.

\begin{proposition}\emph{\textbf{(Hyperbolic BHP)}} \label{abxhp} \, Let $X$ be a complete Riemannian manifold which is $\delta$-hyperbolic with bounded geometry, for some $\delta>0$, and equipped with an adapted weakly coercive operator $L$.\\

 Then, there is a constant $C>1$ depending only on $X$ and $L$, so that for any $z \in \p_GX$ and any two  solutions $u,v >0$ of $L\, w = 0$ on $X$ both  L-vanishing along $\N_i(z) \cap \p_GX$:
\begin{equation}\label{fhepq} u(x)/v(x) \le C \cdot  u(y)/v(y), \mm{ for any two points } x, y \in \N_{i+1}(z) \cap X.
\end{equation}
\end{proposition}

\textbf{Proof} \quad This is a version of [An3], Th.5'. The validity of the assumptions is verified in [BHK],Prop.8.10 (b), the second half of the proof of [BHK],Prop.8.15  and  from \ref{phh} which says that
the $\Phi$, for the $\Phi$-chains used, is determined from the $\delta$ of the $\delta$-hyperbolic space $X$. \qed

We occasionally  need a version of \ref{abxhp} that includes also minimal Green's functions. Implicitly, the version for Green's functions is part of the proof of the usual BHP. But we could not localize an appropriately explicitly stated version. Thus, for completeness, we explain a way to derive it separately. Note that this cannot be accomplished from a simple restriction argument near some point in $\Sigma$, since need to ensure the availability of sufficiently many links between points through $\Phi$-chains. \\

For this, we define the following class $H(X)$ of functions on $X$ possibly with poles:
\begin{equation}\label{gre}
S^*_L(X):= \{u >0 \,|\, u \mm{ solves } L\, w = 0 \mm{ on } X \} \cup \{G(\cdot,p) \,|\, p \in X\}.
\end{equation}

\begin{corollary}\emph{\textbf{(Hyperbolic BHP for Green's Functions)}} \label{abxhpg} \, Let $X$ and $L$ be as in \ref{abxhp}.
 Then, there is a constant $C^*>1$ depending only on $X$ and $L$, so that for any two $u,v  \in S^*_L(X)$ both  L-vanishing along
 $\N_i(z) \cap \p_GX$, there is some $k  >i$, so that:
\begin{equation}\label{fhepq2} u(x)/v(x) \le C^* \cdot  u(y)/v(y), \mm{ for any two points } x, y \in \N_k(z) \cap X.
\end{equation}
\end{corollary}

(Note that \emph{L}-vanishing along $\N_i \cap \p_GX$ is not an extra condition for $G(\cdot,p)$ since it already \emph{L}-vanishes along $\p_GX$.)\\

\textbf{Proof} \quad We reduce this claim to \ref{abxhp}. For this we choose some \emph{minimal} function $u$ \emph{L}-vanishing along $\N_i(z) \cap \p_GX$ and $G(\cdot,p)$, for some $p \in X$. Also we choose some basepoint $q \in X$, $q \neq p$. Then we can find some $k >i$ so that $p,q \notin \N_{k-2}$ and note the following inequality which follows from [An3], Th.5:
\begin{equation}\label{eqq}
C^{-1}_0 \cdot G(x,p)/G(q,p) \le u(x)/u(q)\le C_0 \cdot G(x,p)/G(q,p), \mm{ for any } x \in  \N_k,
\end{equation}
for some $C_0>1$ depending only on $X$ and $L$, cf. [An4],Th.6.2, part A, p.98 (and alternatively the proof of [An3],Th.7) for the details of how to establish this inequality.\\

We combine this inequality with (\ref{fhepq}) to get (\ref{fhepq2}) for any solution $u >0$ of $L\, w = 0$ on $X$ \emph{L}-vanishing along $\N_i(z) \cap \p_GX$. Also applying this argument twice we get the corresponding inequalities when both $u$ and $v$ are Green's functions.\qed

\subsubsection{Martin Boundaries of $(H \setminus \Sigma, d_{\bp^*})$}\label{mss}
\bigskip

Towards Martin theory, we note that the BHP implies, from some by now standardized comparison processes, that the
(minimal) Martin boundary $\p^0_M (X,L)$ equals the topological/ ideal boundary $\p X$ of the underlying space, cf.[Ai1], [AG] or [An3]. We omit here the version
for Euclidean domains which the reader finds broadly discussed in the literature, and directly proceed to the case of Gromov hyperbolic spaces and their Gromov boundary:

\begin{proposition}  \label{abhp} \, Let $X$ be a complete Gromov hyperbolic manifold with bounded geometry equipped with an adapted weakly coercive operator $L$. Then the Gromov compactification $\overline{X}_G$ and the Martin compactification $\overline{X}_M$ are homeomorphic:\[\overline{X}_G \cong \overline{X}_M\] and all Martin boundary points are minimal.\end{proposition}

The proof is given in [An3],Th.2 and Th.8, and [An4], Th.V.6.2, p.97, cf. also [BHK],Ch.8 for further details.\qed

This general theory leads us to the following basic result for the smooth unfolding $X=(H \setminus \Sigma, d_{\bp^*})$ of $H \setminus \Sigma$.

\begin{proposition}\emph{\textbf{(Martin Boundary of Unfoldings)}}\label{gr} \, Let $H \in {\cal{H}}$  and  $L$ be an adapted weakly coercive operator over
$X = (H \setminus \Sigma, d_{\bp^*})$. Then, we have for compact resp. non-compact $H$:  \begin{itemize}
    \item The identity $id_{H \setminus \Sigma}$ extends to homeomorphisms between $H$ resp.$\widehat{H}$, the Gromov compactification
    $\overline{X}_G$ and the Martin compactification $\overline{X}_M$ of $X$.
    \item In particular, $\Sigma \subset H$ resp. $\widehat{\Sigma} \subset \widehat{H}$, the Gromov boundary $\p_GX_{\bp^*}$ and the Martin boundary $\p_M (X,L)$ are homeomorphic.
        \item All boundary points are minimal points, that is,  $\p^0_M (X,L) \equiv \p_M(X,L)$.
\end{itemize}
Thus we have
\[\Sigma \mm{ resp. } \widehat{\Sigma} \cong\p_GX_{\bp} \cong \p_GX_{\bp^*} \cong \p^0_M (X,L) \equiv \p_M (X,L),\] where $ \cong$ means homeomorphic.
\end{proposition}

\textbf{Proof} \quad The properties of $(H \setminus \Sigma, d_{\bp^*})$ we use are, that it is as complete smooth Riemannian manifold, which is Gromov hyperbolic and has
bounded geometry. Also we know from \ref{xgb} that the identity map on $H \setminus \Sigma$ extends to homeomorphisms between $H$ resp. $\widehat{H}$ and
$\overline{X}_G$. Then \ref{abhp} gives us the asserted properties for $\overline{X}_M$ for any adapted weakly coercive operator $L$. \qed

\subsubsection{Skin Adapted Operators on $(H \setminus \Sigma, g_H)$}\label{agrm}
\bigskip

Now we want to see how this Martin theory on the unfolding $X=(H \setminus \Sigma, d_{\bp^*})$ can be transferred to the original space $(H \setminus \Sigma, g_H)$.\\

 To this end, we start on the input side and first formulate a counterpart to the adapted weakly coercivity on $(H \setminus \Sigma, d_{\bp^*})$ for $H \setminus \Sigma$. These are the skin
adapted operators. For the output, that is the Martin theory for these operators,  we establish a simple correspondence between the Martin theories
on $(H \setminus \Sigma, d_{\bp^*})$ and $H \setminus \Sigma$.\\

\textbf{Skin Adapted Operators} \quad We first define \textbf{skin adapted charts} for any $H \in \cal{H}$:  For any $K >1$, there is a radius $\theta(p):= \gamma / \bp(p)$, for
some $\gamma(H,K,\bp) >0$, so that  for any $p \in H \setminus \Sigma$, the exponential map $\exp_p|_{B_{\theta(p)}(0) \subset T_pH}$  is a $K$-bi-Lipschitz
$C^\infty$-diffeomorphism onto its image. Details can be seen e.g. in the proof of \ref{sgp} or [L1], 4.1. This way, we get the smooth charts
\[ \psi_{p}:=\exp^{-1}_p|_{B_{\theta(p)}(0)}: B_{\theta(p)}(p) \ra B_{\theta(p)}(0) \subset \R^n, \, \psi_{p}(p)=0,\] on  $H \setminus \Sigma$.
While these are our default charts, we notice that, in what follows, we may use any fixed collection of charts, for some $K >1$, with these properties.

\begin{definition} \label{sao} \quad  For an area minimizing hypersurface $H$ with singular set $\Sigma$ and given skin transform $\bp$, we call a second order elliptic operator
$L$ on $H \setminus \Sigma$ a \textbf{skin adapted operator}, supposed\\

\emph{\textbf{$\mathbf{\bp}$-Adaptedness}} \, $L$ satisfies skin weighted uniformity conditions relative to the charts $\psi_p$: \small{\begin{equation}
        \label{ada} -L(u) = \sum_{i,j} a_{ij} \cdot \frac{\p^2 u}{\p x_i \p x_j} + \sum_i b_i \cdot \frac{\p u}{\p x_i} + c \cdot u,\end{equation}}
for some locally $\beta$-H\"{o}lder continuous coefficients $a_{ij}$, $\beta \in (0,1]$, measurable functions $b_i$, $c$, and some $k \ge 1$, so that for any $x \in H \setminus
\Sigma$
\[k^{-1} \cdot\sum_i \xi_i^2
\le \sum_{i,j} a_{ij} \cdot \xi_i \xi_j \le k \cdot \sum_i \xi_i^2, \,\,  \delta^{\beta}_{\bp} \cdot  |a_{ij}|_{C^\beta(B_{\theta(p)}(p))} \le k,  \,\,\delta_{\bp}
        \cdot b_i \le k \mm{ and } \delta^2_{\bp} \cdot c \le k.\]

\emph{\textbf{$\mathbf{\bp}$-Weak Coercivity}} \, There exists a positive supersolution $u$ of the equation $L \, f = 0$ so that: \,  $L \, u \ge \ve \cdot \bp^2 \cdot u
, \mm{ for  some, } \ve  >0.$\\
\end{definition}

In chapter 5 below we discuss typical classes of examples and sources of such operators. They most naturally appear in the context of variational problems on $H$.\\

We observe that the notion of skin adaptedness on $(H \setminus \Sigma, g_H)$ the counterpart of adapted weakly coercivity on $(H \setminus \Sigma, d_{\bp^*})$.

\begin{proposition}\emph{\textbf{(Unfolding Correspondence)}}\label{ucor} \, In the canonical  correspondence  \[(H \setminus \Sigma, g_H)  \mm{ equipped with }
L \, \rightleftharpoons \, (H \setminus \Sigma, d_{\bp^*}) \mm{ equipped with } \delta_{\bp^*}^2 \cdot L,\] we have the equivalence
\begin{equation}\label{unfold}
L \mm{ is \textbf{skin adapted}} \,\, \Leftrightarrow \,\, \delta_{\bp^*}^2\cdot L \mm{ is \textbf{adapted weakly coercive},}
\end{equation}
\end{proposition}

\textbf{Proof} \quad  For this we take an atlas  $\{\psi_{p}\,|\,  p \in H \setminus \Sigma \}$ of charts $\psi_{p}: B_{\gamma / \bp(p)}(p) \ra \R^n, \, \psi_{p}(p)=0$ for some $\gamma(H,\bp)
>0$ and a bi-Lipschitz constant $K \ge 1$ independent of $H$ so that the given operator $L$ satisfies the skin adaptedness
conditions in \ref{sao}.\\

Then we define the \emph{locally constantly} scaled version of $\psi_{p}$, cf. \ref{bre} below: \[\psi^{\bp}_{p}:(B_{\gamma}(p),\bp^2(p) \cdot g_H) \ra \R^n \mm{ by }
\psi^{\bp}_{p}(x):= \bp(p) \cdot \psi_{p}(x).\]
 Since both, source and goal, have been scaled by the same constant, $\psi^{\bp}_{p}$ is again a $K$-bi-Lipschitz map. \\

With these particular scalings we find a transparent transformation behavior for the given skin adapted operator $L$:  When the new coordinates relative to the local maps
$\psi^{\bp}$ are denoted by $y_i$, then $-L \, u$ can be recomputed, using the chain rule, as
\[\bp^2(p) \cdot \sum_{i,j} a_{ij} \cdot \frac{\p^2 u}{\p y_i \p y_j} + \bp(p) \cdot \sum_i b_i \cdot \frac{\p u}{\p y_i} + c \cdot u.\]

The adaptedness of $L$ relative $\bp$ and the  Lipschitz continuity of both, $\psi^{\bp}_{p}$ and $\delta_{\bp}$, then shows that $\delta_{\bp}^2\cdot L$ satisfies the adaptedness condition
for adapted weakly coercive operators on $(H \setminus \Sigma, d_{\bp^*})$ relative to the charts $\psi^{\bp}_{p}$.\\

The weak coercivity of $L$ relative $\bp$ implies that there a positive supersolution $u$ of the equation $L \, f = 0$ so that \[\delta_{\bp^*}^2\cdot L \, u \ge \ve \cdot
\delta_{\bp^*}^2\cdot \bp^2 \cdot u, \mm{ for  some, } \ve  >0.\] Thus the approximation property $c_1 \cdot \delta_{\bp}(x) \le \delta_{\bp^*}(x)  \le c_2 \cdot
\delta_{\bp}(x)$, gives the
weak coercivity of $\delta_{\bp^*}^2\cdot L$ on $(H \setminus \Sigma, d_{\bp^*})$.\\

For the reversed direction, from adapted weakly coercive to skin adapted operators, one argues completely similarly.\qed

\begin{remark}\label{bre} \quad Relating $L$ to $\delta_{\bp^*}^2\cdot L$ makes this process formally similar to the case of the Laplacian on two-dimensional  unit disk we mentioned in the introduction.
 \[(D, g_{Eucl})  \mm{ equipped with } \Delta_{Eucl}  \, \rightleftharpoons \, (D, g_{hyp}) \mm{ equipped with }  \Delta_{hyp},\]
In this special case, $\Delta_{hyp}$ is the transformed version of $ \Delta_{Eucl}$ and we actually have $\Delta_{hyp}= (1-|x|)^2/4 \cdot \Delta_{Eucl}$. In general, however, $\delta_{\bp^*}^2\cdot L$ is \emph{not} the transformed operator we get from $L$ under the conformal deformation $g_H \mapsto \delta_{\bp^*}^{-2}\ \cdot g_H$.\\

 The option to work with $\delta_{\bp^*}^2\cdot L$ on $(H \setminus \Sigma, d_{\bp^*})$, to draw conclusions for $L$ on $(H \setminus \Sigma, g_H)$,
 is owing to the flexibility of the adaptedness condition, its quasi-isometry invariance. In turn, this relies on the framework of Brelot's potential theory [Br].\qed
\end{remark}

Now we reach the analytic core results for skin adapted operators. We use again a $\Phi$-neighborhood basis $\N_i(z)$ of a given $z \in \widehat{\Sigma}$ in $\widehat{H}$ of \ref{neb}, we defined in the unfolded setting.

\begin{proposition}\emph{\textbf{(BHP for Skin Adapted Operators)}}\label{mbhsq} \quad Let $H \in {\cal{H}}$ and $L$ some skin adapted operator on $H \setminus \Sigma$. Then, there is a constant $C>1$ depending only on $H$ and $L$, so that for any two  solutions $u,v >0$ of $L\, w = 0$ on $H \setminus \Sigma$ both  L-vanishing along $\N_i(z)\cap \widehat{\Sigma}$:
\begin{equation}\label{fhepq1} u(x)/v(x) \le C \cdot  u(y)/v(y), \mm{ for any two points } x, y \in \N_{i+1}(z) \cap H \setminus \Sigma.
\end{equation}
\end{proposition}

This is our counterpart of the uniform version of the Euclidean BHP with the two differences that we assume that $u,v$ solve $L\, w = 0$ on the entire space $H \setminus \Sigma$. Also the $\N_i$ merely form a topological neighborhood basis of $z$. For later reference we also notice the version for minimal Green's functions. We recall the definition
 \begin{equation}\label{gre1}
S^*_L(H \setminus \Sigma):= \{u >0 \,|\, u \mm{ solves } L\, w = 0 \mm{ on } H \setminus \Sigma \} \cup \{G(\cdot,p) \,|\, p \in H \setminus \Sigma\}.
\end{equation}

\begin{corollary}\emph{\textbf{(BHP for Green's Functions)}} \label{mbhsqg} \, Let $H \in {\cal{H}}$ and $L$ some skin adapted operator on $H \setminus \Sigma$.
 Then, there is a constant $C^*>1$ depending only on $H$ and $L$, so that for any two $u,v  \in S^*_L(H \setminus \Sigma)$ both  L-vanishing along $\N_i(z) \cap \widehat{\Sigma}$, there is some $k  >i$, so that:
\begin{equation}\label{fhepq3} u(x)/v(x) \le C^* \cdot  u(y)/v(y), \mm{ for any two points } x, y \in \N_k(z) \cap H \setminus \Sigma.
\end{equation}
\end{corollary}

Also we rewrite the inequality (\ref{eqq}) needed to derive (\ref{fhepq3})  from  (\ref{fhepq1}): for some $k >i$ so that $p,q \notin \N_{k-2}$ we have for any minimal function $u>0$ \emph{L}-vanishing along $\N_i(z) \cap \widehat{\Sigma}$
\begin{equation}\label{eqq2}
C^{-1}_0 \cdot G(x,p)/G(q,p) \le u(x)/u(q)\le C_0 \cdot G(x,p)/G(q,p), \mm{ for any } x \in  \N_k,
\end{equation}
for some $C_0>1$ depending only on $H$ and $L$.\\

An immediate consequence of  \ref{mbhsq} is the following weaker, global variant.

\begin{corollary}\emph{\textbf{(Global BHP)}} \label{mbhsqc}  \, Let $H \in {\cal{H}}$ and $L$ some skin adapted operator on $H \setminus \Sigma$. Then there is a constant $C >1$ depending only on $H$ and $L$ so that: For any point $p \in \widehat{\Sigma}$ and any neighborhood $V \subset \widehat{H}$ of $p$, there is a smaller neighborhood $W$ of $p$
 with $\overline{W} \subset V$, so that for any two  solutions $u,v >0$ of $L\, w = 0$ on $H \setminus \Sigma$ both  L-vanishing along $V \cap \widehat{\Sigma}$.
\begin{equation}  u(x)/v(x) \le C \cdot  u(y)/v(y), \mm{ for any two points } x, y \in W \setminus \widehat{\Sigma}.
\end{equation}
\end{corollary}

Finally, we reach the Martin theory on $(H \setminus \Sigma, g_H)$:

\begin{proposition}\emph{\textbf{(Martin Boundary for Skin Adapted Operators)}}\label{mbhs} \, Let $H \in {\cal{H}}$ and $L$ some skin adapted operator on $H \setminus \Sigma$. Then, we have:
\begin{itemize}
    \item The identity map on $H \setminus \Sigma$ extends to a homeomorphism  between $\widehat{H}$ and the Martin compactification $\overline{H \setminus
        \Sigma}_M$.
         \item All Martin boundary points are minimal:  $\p^0_M (H \setminus \Sigma,L) \equiv \p_M(H \setminus \Sigma,L)$.
\end{itemize} Thus,  $\widehat{\Sigma}$ and the minimal Martin boundary $\p^0_M (H \setminus \Sigma,L)$ are homeomorphic.\\
\end{proposition}

\textbf{Proof} \quad  The results, \ref{mbhsq}, \ref{mbhs} and \ref{mbhsqg}, translate from the results for $(H \setminus \Sigma, d_{\bp^*})$ equipped with $\delta_{\bp^*}^2 \cdot L$ upon applying the unfolding correspondence (\ref{unfold}). The positive factor  $\delta_{\bp^*}^2$ does neither change the notion of positive solutions of $L \, u = 0$ nor that of minimal solutions.\\

In other words, up to the obvious adjustments in their statements, the respective BHP's for positive solutions of $L \, u = 0$ but also two Martin
compactifications and boundaries are simply identical, that is
\[\overline{(H \setminus \Sigma, g_H)}_M \equiv \overline{(H \setminus \Sigma, d_{\bp^*})}_M \]
and, in particular, we have:   $\p_M((H \setminus \Sigma, g_H),L)  \equiv \p_M((H \setminus \Sigma, d_{\bp^*}),\delta_{\bp^*}^2 \cdot L)$. \qed

\begin{remark}\textbf{(Martin Boundary for Domains)}\label{mgr} \,   For some further types of domains $D \subset H$, so that $D$ is a \emph{uniform space},  we can derive completely similar BHP's and Martin theories. Here we consider the following two cases:
\begin{enumerate}
  \item $D \subset H \setminus \Sigma$, where $dist(z,\Sigma) \le a_D \cdot \delta_{\bp}(z)$, for some $a_D >0$.
  \item $\su \subset H \setminus \Sigma$ a skin uniform domain
\end{enumerate}

We choose, in case (i), the density $\mathbf{d}_1(z):=dist(z,\Sigma)$, in case (ii) the merged density $\mathbf{d}_2(z):= \min\{L_{\bp} \cdot dist(z,\p \su),\delta_{\bp}(z)\}$ and define the metrics $d[i]$, $i=1,2$:
\[d[i](x,y)= \inf \Bigl  \{\int_\gamma 1/\mathbf{d}_i(\cdot)  \, \, \Big| \, \gamma   \subset D \setminus \Sigma \mbox{ smooth arc joining }  x \mbox{ and } y  \Bigr \}\]

In both cases $(D \setminus \Sigma, d[i])$ is again a \emph{complete, Gromov hyperbolic spaces with bounded geometry}. Case (ii) was treated in Ch.\ref{bhua}.  The \emph{boundedness of geometry} results,  in case (i), follows from $dist(z,\Sigma) \le a_D \cdot \delta_{\bp}(z)$.\\

Next, on the analytic side we consider a second order elliptic operator $L$ on $D$ with the appropriate adaptedness properties.\\

$\bullet$ \, $L$ satisfies $\mathbf{d}$-weighted uniformity conditions: \small{\begin{equation}
        \label{ada} -L(u) = \sum_{i,j} a_{ij} \cdot \frac{\p^2 u}{\p x_i \p x_j} + \sum_i b_i \cdot \frac{\p u}{\p x_i} + c \cdot u,\end{equation}}
for some locally $\beta$-H\"{o}lder continuous coefficients $a_{ij}$, $\beta \in (0,1]$, measurable functions $b_i$, $c$, and some $k \ge 1$, so that for any
$x \in D$
\[k^{-1} \cdot\sum_i \xi_i^2
\le \sum_{i,j} a_{ij} \cdot \xi_i \xi_j \le k \cdot \sum_i \xi_i^2, \, \mathbf{d}(z)^{\beta} \cdot  |a_{ij}|_{C^\beta(B_{\theta(p)}(p))} \le k,  \,\mathbf{d}(z) \cdot b_i \le k \mm{ and } \mathbf{d}(z)^2 \cdot c \le k.\]

$\bullet$ \, There exists a positive supersolution $u$ of the equation $L \, f = 0$ so that: \[L \, u \ge \ve \cdot \mathbf{d}(z)^{-2} \cdot u
, \mm{ for  some } \ve  >0.\]

Then, as a counterpart to \ref{mbhs}, we have, in case (i), that the identity map on $D$ extends to a homeomorphism  between $\widehat{D}$ and the Martin compactification $\overline{D}_M$, all Martin boundary points are minimal and $\widehat{\p D}$ and the minimal Martin boundary $\p^0_M (D,L)$ are homeomorphic.\\

In case (ii), the minimal Martin boundary is again the Gromov boundary, cf.\ref{xgbiu}.  \qed
\end{remark}

\setcounter{section}{4}
\renewcommand{\thesubsection}{\thesection}
\subsection{Asymptotic Analysis near $\Sigma$}
\bigskip

We use the BHP for skin adapted operators and the associated Martin theory to derive results describing the asymptotic behavior of positive solutions of the related equations.
Remarkably, although in many cases the solutions diverge towards the singular set or admit merely discontinuous extensions, we will observe that the quotient of any two such solutions
remains well-behaved.\\

\subsubsection{Non-Tangentiality and Fatou Theorems}\label{nt00}
\bigskip

To motivate the results in this section, we recall a classical version of so-called Fatou type theorems, cf.[St],Ch.VII,1.:\\

\emph{For any bounded harmonic functions $h$ on the unit disc $B_1(0) \subset \R^2$, there is an extension of $h^*$ of $h$ onto $\overline{B_1(0)}$, so that for
\emph{almost} any $p \in S^1$: $h(x) \ra h^*(p)$,  for $x \ra p$, provided $x \in
C_\theta$, the inner cone with any angle $\theta \in (0,\pi)$ pointing to $p$.}\\

This was extended by Na\"{i}m, Doob and Gowrisankaran, cf.[Db] and [Go], to relative versions of quotients of not necessarily bounded harmonic functions on open subsets of locally compact Hausdorff
spaces in the axiomatic potential theory \`{a} la Brelot. In this generality non-tangential limits are viewed as fine limits and the boundary is the minimal Martin boundary.\\

Here we prove a Fatou theorem on $H \setminus \Sigma$.  We recall the following notion to describe accessing $\Sigma_H \subset H$ non-tangentially in some quantitative way.

\begin{definition} \label{pen} \quad For $H \in {\cal{H}}$ and any $\rho >0$, we define the \textbf{regular pencil} $\P(z,\rho)$, pointing to $z \in \Sigma$
\begin{equation}\label{ii}
 \P(z,\rho):= \{x \in H \setminus \Sigma \,|\, \delta_{\bp} (x) > \rho \cdot d_{g_H}(x,z)\},
\end{equation}
where $arctan(\rho^{-1})$ is a measure for the size of the (solid) opening angle of $\P(z,\rho)$.
\end{definition}

Next we  compare concepts of non-tangentiality in
$(H \setminus \Sigma, d_{\bp})$ and $(H \setminus \Sigma, g_H)$.
To this end, we consider the distance tube $U_r(\gamma_z)$, for $r >0$, around a geodesic ray $\gamma_z \subset (H \setminus \Sigma,
d_{\bp^*})$ (respectively in $(H \setminus \Sigma, d_{\bp})$) from a basepoint $p \in H \setminus \Sigma$ to $z \in \Sigma$, that is, $\gamma_z$ represents $z$.

 \begin{lemma} \emph{\textbf{(Cylinders in $(H \setminus \Sigma, d_{\bp})$ $\rightleftharpoons$ Pencils in $(H \setminus \Sigma, g_H)$)}}  \label{cy} \quad For any $\rho >0$, and sufficiently small $\ve >0$ there is a some large $r
>0$ so that $\P(z,\rho) \cap B_{\ve}(z)$ equipped with
$d_{\bp^*}$ resp. $d_{\bp}$ is contained in $U_r(\gamma_z)$.
\end{lemma}

 {\bf Proof} \quad  Since $d_{\bp^*}$ and $d_{\bp}$ are quasi-isometric, it is enough to consider  tubes around such geodesics relative $d_{\bp}$.\\

We  start from a tube around $\gamma_z$ in $(H \setminus \Sigma, d_{\bp})$ and we recall from \ref{hu} that $\gamma_z$
  is c-skin uniform relative $(H \setminus \Sigma, g_H)$, for some $c > 0$. In particular, we have $l_{min}(\gamma_z(q)) \le c \cdot \delta_{\bp}(q)$ for any point $q \in \gamma_z$.
Since we are only interested in what happens close to $z$, we can assume that $l_{min}(\gamma_z(q))$ equals the length of the subarc to $z$ and, hence, we have
$d_{g_H}(q,z) \le c \cdot \delta_{\bp}(q)$. That is, we can assume that for  $\gamma_z \subset \P(z,\rho)$, for any $\rho $ with $1/c > \rho >0$.\\

In this situation we can benefit from tangent cone approximations, for instance, in the form of the freezing effect of $H$ around the singular point $z$, cf. [L2],Ch.2.1:\\

Let  $\delta > 0$, $1/c > \rho >0$ and $R \gg 1 \gg r >0$ be given.
Then, for some $\tau_{\delta ,\rho, R , r,p} \gg 1$ we have: for every $\tau \ge \tau_{\delta ,\rho, R , r,p}$, there is a tangent cone $C_p^\tau$ of $H$ at $p$ so that
    \[\tau \cdot  B_{R/\tau}(z) \setminus B_{r/\tau} (z) \cap   \P(z,\rho)  \subset \tau \cdot H\] can be written as a smooth section
   $\Gamma_\tau  \mm{ with } |\Gamma_\tau|_{C^5} < \delta$  of the normal bundle of
\[ B_R(0)\setminus B_r (0) \cap   \P(0,\rho)  \subset C_p^\tau.\]
Thus we can further reduce the consideration to the case where $\gamma_z \subset C_p^\tau$ with
\[B_R(0)\setminus B_r (0) \cap \gamma_z \subset B_R(0)\setminus B_r (0) \cap   \P(0,\rho)\]
for $1/c > \rho >0$.  Then we apply the skin uniformity of $S_C \setminus \sigma_C$, namely, we use [L1],Cor.4.7 and its refinement 4.11. From these results we see that any two points $a,b \in \{x \in S_C \setminus \Sigma_{S_C} \,|\, \delta_{\bp}(x) \ge a\} \subset S_C$ can be connected by an arc $\alpha_{a,b} \subset
\{x \in S_C \setminus \Sigma_{S_C} \,|\, \delta_{\bp}(x) \ge \eta(a) \cdot a\}$, for some $\eta(a) \in (0,1)$ of length $\le l(a)$, with $\eta(a)$ and $l(a)$ independent of $C$.\\

From \ref{pro} we get the following estimate
\[d_{\bp}(y, B_R(0)\setminus B_r (0) \cap \gamma_z) \le l(\rho)/(\eta(\rho) \cdot \rho), \mm{ for any } y \in B_R(0)\setminus B_r (0) \cap   \P(0,\rho). \]
Therefore, for sufficiently fine tangent cone approximations, we can choose $r:= 2 \cdot l(\rho)/(\eta(\rho) \cdot \rho)$ to
ensure that $\P(z,\rho) \cap B_{\ve}(z)$  equipped with $d_{\bp}$ is contained in $U_r(\gamma_z)$.\qed

The complementary result, showing that cylinders in $(H \setminus \Sigma, d_{\bp})$ are also contained in suitable pencils in $(H \setminus \Sigma, g_H)$, is left to the reader. We do not use this inclusion in our further discussion.\\

Now we reach the counterpart of the classical Fatou theorem on $B_1(0) \subset \R^2$. We recall from \ref{mrt} that we can write any positive solution of $L \,v=0$ on $H \setminus \Sigma$ for skin adapted operators in terms of the Martin
integral $u_{\mu}(x) =\int_{\widehat{\Sigma}} k(x;y) \, d \mu(y)$,
for some suitable finite Radon measure $\mu$.

\begin{proposition}\emph{\textbf{(Relative Fatou Theorem on $\mathbf{H \setminus \Sigma}$)}}\label{ftt} \, Let $H \in {\cal{H}}$ and $L$ some skin adapted operator on $H \setminus \Sigma$.
For any two finite Radon measures $\mu$, $\nu$ on $\widehat{\Sigma}$ we have: For $\nu$-almost any $z \in \widehat{\Sigma}$ and any $\rho >0$:
 \begin{equation}\label{rn}
u_{\mu}/u_{\nu}(x) \ra d \mu/d \nu(z),\,\mm{ for } x \ra z, \mm{ with } x \in \P(z,\rho),
\end{equation}  where $d \mu/d \nu$ denotes the Radon-Nikodym derivative of $\mu$ with respect to $\nu$.
\end{proposition}

To explain the statement we recall that for any two Radon measures $\mu$, $\nu$ on $\widehat{\Sigma}$, we have a (uniquely determined) Lebesgue decomposition
  $\mu=\mu_1 + \mu_2$ relative $\nu$ into some $\nu$-absolute continuous measure $\mu_1$ and a $\nu$-singular measure $\mu_2$.\\

That is, $\mu_1(E)= \int_E f \, d\nu, \mm{ for any measurable } E \subset \widehat{\Sigma},$
for some $\nu$-integrable function $f$ which is the Radon-Nikodym derivative $d \mu_1/d \nu$ and uniquely determined outside a set of $\nu$-measure zero. And, for
$\mu_2$, there is a set $F \subset \widehat{\Sigma}$, with $\nu(F)=0$ and $\mu_2(\widehat{\Sigma} \setminus F)=0$. Then \ref{ftt} reads as follows. For $\nu$-almost any $z \in \widehat{\Sigma}$:
\[u_{\mu}/u_{\nu}(x) \ra f(z) = d \mu_1/d \nu(z),\,\mm{ for } x \ra z, \mm{ with } x \in \P(z,\rho),\]
 for any given $\rho >0$.\\

\textbf{Proof} \quad From the unfolding correspondence in the proof of \ref{mbhs}  we get the following result for $(H \setminus \Sigma, d_{\bp^*})$ equipped with
$\delta_{\bp^*}^2 \cdot L$ from Ancona theory, in this case, [An4],p.100, Th.6.5 and the definitions on p.99, which, in turn, is an adaptation of the work of Gowrisankaran in [Go].\\

For any two finite Radon measures $\mu$, $\nu$ on $\p^0_M((H \setminus \Sigma, d_{\bp^*}),\delta_{\bp^*}^2 \cdot L)$ we get for almost any $z \in \p^0_M((H \setminus
\Sigma, d_{\bp^*}),\delta_{\bp^*}^2 \cdot L)$ and any distance tube $U_r(\gamma_z)$, for $r >0$, around a geodesic ray $\gamma_z \subset ((H \setminus \Sigma,
d_{\bp^*})$ from a basepoint $p \in H \setminus \Sigma$ to $z$, that is, $\gamma_z$ represents $z$:  for the solutions $u_{\mu}, u_{\nu}$ associated to $\mu$, $\nu$
\[u_{\mu}/u_{\nu}(x) \ra d \mu/d \nu(z),\,\mm{ for } x \ra z, \mm{ with } x \in U_r(\gamma_z).\]

However, \ref{cy} shows that for any $\rho >0$, and sufficiently small $\ve >0$ there is a some large $r >0$ so that $\P(z,\rho) \cap B_{\ve}(z)$ equipped with
$d_{\bp^*}$ is contained in $U_r(\gamma_z)$ and the result follows.  \qed

\begin{remark}  \quad   For measures with positive density along $\Sigma$ there actually is no general control over the tangential behavior of quotients $u_{\mu}/u_{\nu}$ of solutions.  Due to Littlewood and refined by Aikawa, cf.[Ai4], there are examples of  bounded harmonic functions $h$ on the unit disc $B_1(0) \subset \R^2$ so that the tangential limit, that is the limit along paths approaching $S^1$ tangentially, does \emph{not} exist in any $p \in S^1$.  It should be possible to derive similar results for skin adapted operators on $H \setminus \Sigma$.\qed
 \end{remark}

\subsubsection{Continuous Extension Problems to $\Sigma$}\label{nt}
\bigskip

We turn to a situation which is largely complementary to that in the Fatou theorem: we consider quotients of solutions $u_{\mu}, u_{\nu} >0$ on $H \setminus \Sigma$ whose associated Radon measures $\mu,\nu$ vanish on a common open subset $A \subset \Sigma$,  $\mu(A)=\nu(A)=0$.  Whereas, in this situation, the Fatou theorem makes no predictions, we show that, by the virtue of the BHP, the quotient $u_{\mu}/u_{\nu}$ remains well-controlled and this even holds in tangential directions. To this end we first notice

\begin{lemma} \label{cea}  For any $H \in {\cal{H}}$, any skin adapted operator $L$ on $H \setminus \Sigma$
and any open subset $A \subset \widehat{\Sigma}$, we have: $ \mu(A)=0 \, \Rightarrow \, u_{\mu} $ \emph{L}-vanishes along \emph{A}.\\
\end{lemma}

\textbf{Proof} \quad If $\mu(A)=0$, then the Martin integral and inequality (\ref{eqq2}) for minimal functions, in the BHP for Green's Functions, show that near each $z \in A$, there is a constant $c_z>0$, so that for the Green's functions $G(\cdot,p)$ for some basepoint $p \in H \setminus \Sigma$
\begin{equation}\label{lva} G(\cdot,p) \ge c_z \cdot u_{\mu}\mm{ near }z.\end{equation}
But $G(\cdot,p)$ $L$-vanishes along $\Sigma$ and thus $u_{\mu}$ $L$-vanishes along $A$. \qed

The latter result shows that solutions L-vanishing along $A$ are naturally related to a vanishing of the Radon measure they define. For any pair of such solutions we can now formulate one of the main results of this paper.

\begin{proposition}\emph{\textbf{(Continuous Extensions to $\mathbf{\Sigma}$)}}\label{cee}  For any $H \in {\cal{H}}$, any skin adapted operator $L$ on $H \setminus \Sigma$
and any two $u,v \in S^*_L(H \setminus \Sigma)$,  L-vanishing  along some common open subset $A \subset \widehat{\Sigma}$,   we have:\\

 The quotient $u/v$ on $H \setminus \Sigma$ admits a continuous extension to $H \setminus \Sigma \cup A \subset \widehat{H}$.\\

\end{proposition}

\textbf{Proof} \quad We modify a classical method due to Moser in [Mo],Ch.5, cf. [JK],7.9 and [Ai2],Th.2 for the case of boundary problems, originally employed to derive relative estimates for the oscillation of harmonic functions on concentric Euclidean balls. \\

In place of these concentric balls, we use a $\Phi$-neighborhood basis $\N_i(z) \subset H$ of (\ref{neb}), for any given $z \in \widehat{\Sigma}$. Then, we get Harnack estimates on the $\N_i$ from the BHP for skin adapted operators of \ref{mbhsq} and \ref{mbhsqg}. We set
\begin{equation}\label{mi}
\sup(k):= \sup_{\N_k} u/v  \mm{ and } \inf(k):= \inf_{\N_k} u/v
\end{equation}
 The oscillation of $u/v$ on $\N_k$ is written as  $osc(k)= \sup(k)-\inf(k)$.\\

We first note from \ref{mbhsq} resp.  \ref{mbhsqg} that $\sup(k_0) < \infty$, for any $k_0$ sufficiently large, so that the poles of $u$ and $v$ do not belong to $\N_{k_0-2}$.\\

For $k \ge k_0$, we consider the two solutions $\sup(k) \cdot v -u$ and  $u - \inf(k) \cdot v$. They are positive due Hopf's maximum principle and they also \emph{L}-vanish along $A$. Therefore the BHP \ref{mbhsq}  resp.  \ref{mbhsqg}  also applies to this pair of functions and shows for $\N_{k+1}$:
\[\sup_{\N_{k+1}} (\sup(k) \cdot v -u)/v  \le C^* \cdot  \inf_{\N_{k+1}} (\sup(k) \cdot v -u)/v\]
\[ \sup_{\N_{k+1}} (u - \inf(k) \cdot v)/v  \le C^* \cdot  \inf_{\N_{k+1}} (u - \inf(k) \cdot v)/v\]
From these inequalities we get:
{\small \[ \sup(k) - \inf(k+1) \le C^*  \cdot (\sup(k) - \sup(k+1)) \,\mm{ and  } \,\sup(k+1) - \inf(k) \le C^*  \cdot (\inf(k+1) - \inf(k))\]}
We add suitable multiples of these inequalities to get
\[ osc(k+1)= \sup(k+1) - \inf(k+1) \le \]
\[(C^*-1)/(C^*+1) \cdot (\sup(k) - \inf(k))=(C^*-1)/(C^*+1) \cdot osc(k) \]
In other words, for $a:=(C^*-1)/(C^*+1) <1$: \,
\begin{equation}\label{oscc}
osc(k) \le a^{k-k_0} \cdot osc(k_0) \ra 0, \mm{ for } k \ra \infty.
\end{equation}
Thus, $u/v$ remains bounded near $z$, directly from the BHP, and from (\ref{oscc}) it admits a continuous extension to $z$. \qed

\begin{remark}\textbf{ (H\"older and Equi-Continuity)}\label{cee1} \, For sufficiently regular Euclidean domains,
e.g. Lipschitz or uniform domains, cf. [JK],7.9 and [Ai2],Th.2, the latter argument even shows the \emph{H\"older continuity} of the extension of quotients of positive harmonic functions vanishing along (parts of) the boundary.  Here one can take profit from the underlying Euclidean structure and choose concentric Euclidean balls with the suitably adjusted radii.\\

It is conceivable that such arguments still apply to the unfolding $(H \setminus \Sigma, d_{\bp^*})$. This is more delicate since Gromov boundaries generally do not carry a canonical metric. But for $(H \setminus \Sigma, g_H)$ we only reach ordinary continuity since we have no metric control, relative $g_H$, over the way the family $\N_i(z)$ shrinks to $z$. \\

However, we observe from (\ref{oscc}) that the \emph{family of quotients} $u/v$, for $u,v \in S_L(H \setminus \Sigma)$,  \emph{L}-vanishing  along some common open subset $A \subset \widehat{\Sigma}$ is both \emph{uniformly bounded} and \emph{equi-continuous} in a neighborhood $W$ of any compact subset $K \subset A$. That is, due to the Arzel\`{a}-Ascoli theorem, any sequence of such quotients admits a $C^0$-convergent subsequence on $W$.  \qed
\end{remark}

For certain types of operators we can consider the absolute value of positive solutions on $H \setminus \Sigma$, then we also find the solvability of the Dirichlet problem

\begin{proposition}\emph{\textbf{(Dirichlet Problem for Skin Adapted Operators)}}\label{diri} \quad Let $H \in {\cal{H}}$ and $L$ some skin adapted operator on $H \setminus \Sigma$,
so that constant functions solve $L \, v=0$ and $G(x,p)\ra 0$,  for $x \ra \widehat{\Sigma}$, and given  $p \in H \setminus \Sigma$.\\

 Then, for any continuous function $f$ on $\widehat{\Sigma}$, there is a uniquely determined continuous function $F$ on $H$ so that
 \[F|_{\widehat{\Sigma}} \equiv f \mm{ and } L \, F =0.\]\
\end{proposition}

\textbf{Proof} \quad The assumptions and the Martin theory for skin adapted operators allow us to imitate the standard arguments well-known e.g. for harmonic functions on Euclidean domains. We first define the Radon measure $\mu^{\cs}$ on $\widehat{\Sigma}$ as the one associated to the constant function $1$, decompose $f$ into $f=f^+-f^-$, for $f^+:=\max\{f,0\}$, $f^-:=-\min\{f,0\}$ and set
\[F^\pm(x) =  \int_{\widehat{\Sigma}} f^\pm(y) \cdot k(x;y) \,d \mu^{\cs}(y)\]
Then $F^\pm \ge 0$ solves $L \, w=0$. We assert that $F^\pm$ extends continuously to $H$
and coincides with $f^\pm$ on $\widehat{\Sigma}$. This then implies the claims for $F:=F^+-F^-$.\\

We have to show that for any $z \in \widehat{\Sigma}$, $F^\pm(x)\ra f^\pm(z)$,  for $x \ra z$, $x \in H \setminus \Sigma$. To this end we notice from (\ref{lva}) that for any fixed neighborhood $U \subset \widehat{\Sigma}$ of $z$:
\[G(\cdot,p) \ge c_z \cdot \int_{\widehat{\Sigma} \setminus U} f^\pm(y) \cdot k(x;y) \,d \mu^{\cs}(y)\]
close to $z$. Thus, for any such $U$, we have $\int_{\widehat{\Sigma} \setminus U} f^\pm(y) \cdot k(x;y) \,d \mu^{\cs}(y) \ra 0$, for $x \ra z$.\\

In turn, since $f^\pm$ is continuous on $\widehat{\Sigma}$, we have $|f^\pm(y) - f^\pm(z)| < \ve$, for $y \in U(\ve)$, where $U(\ve) \subset \widehat{\Sigma}$ is a small enough neighborhood of $z$. Moreover, since $ \int_{\widehat{\Sigma}}  k(x;y) \,d \mu^{\cs}(y)=1$, by definition of $\mu^{\cs}$ and again  $\int_{\widehat{\Sigma} \setminus U(\ve)} k(x;y) \,d \mu^{\cs}(y) \ra 0$, for $x \ra z$,  we have for any $x \in H \setminus \Sigma$ close enough to $z$:
\[\left|\int_{U(\ve)} f^\pm(y) \cdot k(x;y) \,d \mu^{\cs}(y) - \int_{U(\ve)} f^\pm(z) \cdot k(x;y) \,d \mu^{\cs}(y)\right| \le 2 \cdot \ve, \mm{ for } y \in U(\ve).\]
We conclude that $F^\pm(x)  \ra f^\pm(z)$,  for $x \ra z$.\\

Finally, the maximum principle shows that $F \equiv 0$ is the only solution with $F|_{\widehat{\Sigma}} \equiv 0$ and thus we find the asserted unique solvability of the Dirichlet problem.
 \qed

\bigskip
\setcounter{section}{5}
\renewcommand{\thesubsection}{\thesection}
\subsection{Symmetric Operators}
\bigskip

In this chapter we particularly focus on eigenvalue problems. This makes it convenient to restrict to the case of symmetric operators.\\

The first observation for singular area minimizers, without counterpart on closed smooth manifolds, is the existence of a principal eigenvalue so that for any real value strictly
smaller we still find positive eigenfunction. Indeed this is case of most interest to us, since this is where the previously set-up Martin theory applies.\\

We use the Hardy inequality for skin transforms to show how this applies to classical elliptic operators. Many of them are actually Schr\"odinger operators and in this case we
 setup an inductive description of eigenfunctions.\\

To simplify the exposition, we make the mild regularity assumption that the coefficients of the given operator $L$ are at least $\alpha$-H\"older continuous, for some $\alpha \in (0,1)$. Hence, weak solutions of $L \, v=0$ all equations, including $\bp$-weighted versions, considered here, are $C^{2,\alpha}$-regular,cf. [BJS], Part II,,Ch.1.3. Since $\bp$ is locally Lipschitz regular, this also applies to solutions of weighted eigenvalue equations: $L \, v= \lambda \cdot \bp^2 \cdot v$, for some $\lambda \in \R$.\\

\subsubsection{Criticality and Principal Eigenvalues}\label{bgrb}

\bigskip

The eigenvalue theory for many $\bp$-adpated operators can be treated form that of some associated skin adapted operators we get from the simple shift of adding the term
$\lambda \cdot \bp^2 \cdot Id$, for some suitable $\lambda \in \R$. This is an important way to extend the range of problems the developed Martin theory can be applied to. Concrete examples are considered in the two sections that follow below.\\

\textbf{Shifted Skin Adaptedness} \, We start with an $\bp$-adpated symmetric operator $L$ on $H \setminus \Sigma$. Then, the $\bp$-weak coercivity can be expressed in terms of a variational integral for eigenvalues. This is well-known for the Laplacian on Euclidean domains, cf.[An1],Prop.1/Appendix. For completeness, we explain how this carries over to our case:

\begin{lemma} \emph{\textbf{(Hardy Inequalities)}} \label{hi} \, An $\bp$-adapted symmetric operator $L$ on $H \setminus \Sigma$ is $\bp$-weakly coercive if and only if the Hardy inequality
 {\small \begin{equation} \label{hadi}  \int_H  f  \cdot  L f  \,  dA \, \ge \, \tau \cdot \int_H \bp^2\cdot f^2
dA,\end{equation}}
holds for any smooth $f$, compactly supported in $H \setminus \Sigma$, and some positive constant $\tau = \tau(L,\bp,H)>0$.
\end{lemma}

 \textbf{Proof} \quad The vector space $V_L:=\{f \in L^2_{loc}(H \setminus \Sigma)\,|\,\nabla f, \bp \cdot f  \in L^2(H \setminus \Sigma)\}$
can be equipped with the bilinear form: $\langle u,v \rangle_{V_L}:= \int_{H \setminus \Sigma} u \cdot L v   +  \bp^2 \cdot u  \cdot v.$\\

When (\ref{hadi}) is satisfied, $\langle u,v \rangle_{V_L}$ becomes scalar product and that makes $V_L$ a Hilbert space. We define $V_L^0 \subset V_L$ as the $V_L$-closure of $C_0^\infty(H \setminus \Sigma)$. Now, we choose some $\ve \in \R$ and consider the following continuous quadratic form \[a(f,f):= \int_{H \setminus \Sigma} f \cdot L f   - \ve \cdot  \bp^2 \cdot f^2 \mm{ on } V_L^0.\]
For $\ve/ \tau \in (0,1)$, we have:
 \[a(f,f) \ge (1- \ve/ \tau) \int_{H \setminus \Sigma} u \cdot L v  \ge (1- \ve / \tau)/(1+ 1/\tau) \cdot |f|_{V_L}^2, \mm{ for any } f \in V_L^0.\]
Thus, $a(\cdot,\cdot)$ is a coercive bilinear form on $V_L^0$ and due to the symmetry of $L$ it is symmetric. Therefore we may apply the Lax-Milgram theorem to show that for any  $h \in C_0^\infty(H \setminus \Sigma)$
there is a unique $f \in V_L^0$ with $a(f,v)= \int_{H \setminus \Sigma} h \cdot v$, for any $v \in C_0^\infty(H \setminus \Sigma)$ and we infer that $L f + \ve \cdot \bp^2 \cdot f=h$. Now we choose some $h >0$, then \[a(\min\{0,f\},\min\{0,f\})=-a(f,v)=-\int_{H \setminus \Sigma} h \cdot v \le 0.\]
In other words, $f>0$ and satisfies $L f \ge \ve \cdot \bp^2 \cdot f$.\\

In turn, if such an $f>0$  exists, we consider the first Dirichlet eigenvalue $\lambda_D$ of  $L-\ve \cdot \bp^2 \cdot Id$ on any domain $D$ with $\overline{D} \subset H \setminus \Sigma$ and some positive eigenfunction $w$ with $(L-\ve \cdot \bp^2 \cdot Id) w = \lambda_D \cdot w$.\\

 We show that  $\lambda_D > 0$. Otherwise, we first observe that  $k \cdot f > w$ on $H \setminus \Sigma$, for sufficiently large $k >0$ and we take the infimum $k_0>0$ of all such $k$. Then there is some point $p \in D$ so that  $k \cdot f(p) - w(p)=0$ and $k \cdot f - w>0$  in some ball $B \subset D$ with $p \in \p B$. Moreover, $(L-\ve \cdot \bp^2 \cdot Id)(k \cdot f - w) >  0$ on $D$. Then the Hopf maximum principle can be applied to get a contradiction, cf. the proof of \ref{scal2} below for additional comments.\\

 Thus, $\lambda_D > 0$, and this means $\int_{H \setminus \Sigma}  v \cdot  (L-\ve \cdot \bp^2 \cdot Id) v \ge 0$ for any $v  \in C_0^\infty(D)$ and, since $D$ can be chosen arbitrarily large within $H \setminus \Sigma$, we get the Hardy inequality for $\tau = \ve$.
\qed

Thus for any symmetric skin adapted operator  $L$, there is a \emph{largest} $\lambda >0$ so that the Hardy inequality (\ref{hadi}) is satisfied. This characterization of skin
adaptedness suggests to consider also cases where $\lambda$ may be negative but finite: $\lambda > -\infty$.\\

The largest value for $\tau$ in (\ref{hadi}) can be viewed as an eigenvalue of the operator $\delta_{\bp}^2 \cdot L$. In \ref{scal2} below we will see that for \emph{singular} $H$ there are also
positive eigenfunctions for any $\lambda < \lambda^{\bp}_{L,H}$.
This is an important difference to the case of smooth compact manifolds  where the first eigenvalue $\lambda_1$ of such an elliptic
operator $L$ exists as the unique eigenvalue with a  positive eigenfunction $u_1$, which is unique up to multiples, cf.[C], Ch.VI. This detail
will be essential to  incorporate singular area minimizers in scalar curvature geometry in [L2].\\

We introduce the following terminology

\begin{definition} \label{saop} \quad  Let $L$ be an $\bp$-adpated symmetric operator $H \setminus \Sigma$ and assume there is some
\textbf{finite} $\tau \in \R$ so that the Hardy inequality (\ref{hadi}) holds, then $L$ is called a \textbf{shifted skin adapted} operator.\\

The largest value $\lambda^{\bp}_{L,H}$
for $\tau$, so that (\ref{hadi}) holds, is  called  the (generalized) \textbf{principal eigenvalue} of $\delta_{\bp}^2 \cdot L$. In particular, $L$ is skin adapted when $\lambda^{\bp}_{L,H} >0$.
\end{definition}

\textbf{Basic Spectral Theory} \, We get the following transparent description for the spectral theory of $\delta_{\bp}^2 \cdot L$ of a shifted skin adapted $L$, when either $\Sigma \n$, when $H$ is compact,
or if it is non-totally geodesic, when $H \subset \R^{n+1}$. In \ref{regu2} we have a brief look at the spectral theory for  the non-weighted $L$.

\begin{proposition}\emph{\textbf{(Criticality and Principal Eigenvalues)}} \label{scal2} \,  For any $H \in \cal{H}$  so that $H \setminus \Sigma$ is non-compact
and non-totally geodesic, $L$ any shifted skin adapted operator on $H \setminus \Sigma$, we set
\[L_\lambda:= L - \lambda \cdot \bp^2 \cdot Id, \mm{ for }\lambda \in\R.\] Then we have the following trichotomy for the spectral theory of $\delta_{\bp}^2 \cdot L$.
\begin{itemize}
    \item \emph{\textbf{Subcritical}} when $\lambda < \lambda^{\bp}_{L,H}$,  $L_\lambda$ is skin adapted.
    \item \emph{\textbf{Critical}} when $\lambda = \lambda^{\bp}_{L,H}$,  there is an, up to multiples, unique positive solution $\phi$, the \textbf{ground state}* of $L
        _{\lambda^{\bp}_{L,H}}, $  for the equation $L _{\lambda^{\bp}_{L,H}} \ v = 0$.

         The function $\phi$ can be described as the limit of first Dirichlet eigenfunctions for the operator $\delta_{\bp}^2 \cdot L$ on a sequence of smoothly bounded domains
        $\overline{D_m} \subset D_{m+1} \subset H \setminus \Sigma$, $m \ge 0$, with $\bigcup_m D_m = H \setminus \Sigma$.
    \item \emph{\textbf{Supercritical}} when $\lambda > \lambda^{\bp}_{L,H}$, $L_\lambda \, u = 0$ has no positive solution.
\end{itemize}
\end{proposition}
*When $L$ is a Schr\"odinger operator $\phi$ is the \emph{ground state} well-known from the study of quantum mechanical systems, cf.[GJ],3.3.

\begin{remark} \label{regu} \quad  Following this result and the conventions in [P],Ch.4, we also say the operator $L$ is \textbf{subcritical}, when $L$ admits a positive Green's function,
$L$ is  \textbf{critical} when it does \emph{not} admit a Green's function but a positive solution of $L \, v=0$ and
\textbf{supercritical} when the latter equation does \emph{not} admit any positive solutions.
\qed
\end{remark}

\textbf{Proof of \ref{scal2}} \quad The main case is the critical case, where $\lambda = \lambda^{\bp}_{L,H}$. The subcritical case is rather obvious. The supercritical one follows from the discussion of the critical case.\\

\textbf{Subcritical Case} ($\lambda < \lambda^{\bp}_{L,H}$)  \quad $L_{\lambda}$ clearly satisfies the adaptedness condition. The variational relation $\int_H  f  \cdot
L_{\lambda^{\bp}_{L,H}}  f  \,  dA \, \ge \, \lambda^{\bp}_{L,H} \cdot \int_H \bp^2\cdot f^2 dA$  shows that {\small \[\int_H  f \cdot L_{\lambda} f  \,  dA \, \ge
\, (\lambda^{\bp}_{L,H} - \lambda) \cdot \int_H \bp^2\cdot f^2 dA,\]}  and we conclude  that $L_{\lambda}$ is skin adapted.\\

\textbf{Critical Case} ($\lambda = \lambda^{\bp}_{L,H}$)  \quad We choose a sequence of smoothly bounded domains $\overline{D_m} \subset D_{m+1} \subset H
\setminus \Sigma$,
$m \ge 0$, with $\bigcup_m D_m = H \setminus \Sigma$.\\

 Now we consider the  uniquely determined first eigenvalue $\lambda_m$ and first
eigenfunctions $\phi_m$ on $\overline{D_m}$, for the operator $\delta_{\bp}^2 \cdot L$ for the Dirichlet problem, that is, with
\[L  \, \phi_m = \lambda_m \cdot \bp^2 \cdot \phi_m \mm{ and } \phi_m \equiv 0\mm{ on } \p D_m.\]

We may assume that $\phi_m > 0$ on $D_m$ and $\phi_m(p_0) = 1$, in some base point $p_0 \in \bigcap_m D_m$.
Also we formally extend the $\phi_m$ to $H \setminus \Sigma$ setting them equal to $0$ outside $D_m$.\\

The variational (or Rayleigh) characterization of $\lambda_m$ as \[\lambda_m = \inf \Big\{\int_H  f  \cdot  L f  \,  dA / \int_H \bp^2\cdot f^2 dA \,\Big|\, f \mm{ smooth,
supported in } D_m \subset H \setminus \Sigma\Big\}\]
 shows  $\lambda_m > \lambda_{m+1} > \lambda^{\bp}_{L,H} >0$. Also $\lim_{m \ra \infty}
\lambda_m = \lambda^{\bp}_{L,H}$, since the support $supp \, f$ of any smooth function with compact support on $H \setminus \Sigma$ is contained in $D_m$ once
we have chosen $m$ large enough.\\

Now we modify the operator $-L(u) = \sum_{i,j} a_{ij} \cdot \frac{\p^2 u}{\p x_i \p x_j} + \sum_i b_i \cdot \frac{\p u}{\p x_i} + c \cdot u$ and  choose one small ball $B \subset
D_1$. We slightly decrease the value of the function $c$ within $B$ to another $C^\beta$-function, but keep it fixed outside $B$.. We call the new function $c[m]$, that is, we
get an auxiliary, again skin adapted, operator
 \[-L[m](u) = \sum_{i,j} a_{ij} \cdot \frac{\p^2 u}{\p x_i \p x_j} + \sum_i b_i \cdot \frac{\p u}{\p x_i} + c[m] \cdot u.\]
In view of the Rayleigh characterization of the eigenvalue and the fact that $\lambda_m  > \lambda^{\bp}_{L,H}$,  we can choose $c[m]$, so that
\begin{itemize}
    \item the first eigenvalue $\lambda_m[m]$ of the operator $\delta_{\bp}^2 \cdot L[m]$ for the Dirichlet problem on $\overline{D_m}$, becomes
        $\lambda_m[m]=\lambda^{\bp}_{L,H}$,
    \item $c[m] \ra c$, for $m \ra \infty$, in  $C^{\beta}$-norm.
\end{itemize}
Then the new sequence of first eigenfunctions $\phi_m[m]$ with $\phi_m[m](p_0) = 1$ contains a compactly converging subsequence with limit $\phi >0$ on $H \setminus
\Sigma$.
 $\phi >0$ solves the equation $L  \, v = \lambda^{\bp}_{L,H} \cdot \bp^2 \cdot v$ on $H \setminus \Sigma$. Namely, due to $c[m] \ra c$, for $m \ra \infty$ we get Harnack inequalities, cf.[GT], 8.20,  for the $\phi_m[m]$ on domains $D \subset \overline{D} \subset H \setminus
\Sigma$ for sufficiently large
$m$ and with constants independent of $m$.\\

What remains is to show that  $\phi$ is the unique positive solution, up to multiples. Thus let $\psi>0$ be another solution of $L  \, v = \lambda^{\bp}_{L,H} \cdot
\bp^2 \cdot v$ on $H \setminus \Sigma$. Then we can find a constant $k
>0$: $k \cdot \phi \ge \phi_m[m]$ near $\overline{B}$, for any $m$.\\

 From this we can infer that $k \cdot \phi \ge\phi_m[m]$ on the whole of $H \setminus \Sigma$. Otherwise we can find a largest $k$, so that $k \cdot \phi \ge\phi_m[m]$, called $k_0$, so that there is some
  point $p \in  D_m \setminus \overline{B}$
so that $k_0 \cdot \phi(p) = \phi_m[m](p)$. But then $k_0 \cdot \phi -\phi_m[m]$ is a non-negative solution of $L  \, v = \lambda^{\bp}_{L,H} \cdot \bp^2 \cdot v$ on
$D_m
\setminus \overline{B}$ and we get an interior zero point while the boundary values remain positively lower bounded. \\

This contradicts the \emph{Hopf maximum principle} cf.[GT],Ch.3.2.  In this case the maximum principle applies without the typical assumptions for the sign of the zeroth order
coefficient. This is mentioned in a remark that follows the proof of [GT],Th.3.5. Actually, this follows from the strict inequality for the
outer normal derivative $\p (k_0 \cdot \phi -\phi_m[m])/\p n> 0$, [GT],p.24,  e.g. in extremal points of the zero set. It shows that $k_0 \cdot \phi -\phi_m[m] \equiv 0$ on the path component of $p$.\\

Thus we also get $k \cdot \psi \ge \phi$, for some suitable $k>0$. Now we choose the smallest such $k$, called $k_1$. Then, either $k_1 \cdot \psi \equiv \phi$ or $k_1 \cdot
\psi > \phi$. In the latter case, we write $F:=k_1 \cdot \psi - \phi >0$ and repeat the previous argument to find some $l>0$ with $l\cdot F \ge \phi$. But then we also have $k_1
\cdot \frac{l}{l+1} \cdot
\psi \ge \phi$ contradicting the definition of $k_1$. Thus the solution $\phi$ is uniquely determined up to positive multiples.\\

\textbf{Supercritical Case} ($\lambda > \lambda^{\bp}_{L,H}$)  \quad In this case $L_\lambda \, u = \lambda \cdot \bp^2 \cdot u$ has no positive solution. Otherwise, assume we had such a solution $u >0$.\\

Then, we can find a smoothly bounded domain $D \subset \overline{D} \subset H \setminus \Sigma$ so that the eigenvalue $\lambda_D$ for the first Dirichlet eigenfunction
$\phi_D$ for the operator
$\delta_{\bp}^2 \cdot L$ on $D$ equals $\lambda$.\\

We can find some constant $k>0$, so that $k \cdot u \ge \phi_D$. For the smallest $k$, with $k \cdot u \ge \phi_D$, we get an interior zero point while the boundary
values remain positively lower bounded. Again, this contradicts the Hopf maximum principle in the same fashion we have already seen before.\qed

\begin{remark} \textbf{(Non-Weighted Operators)} \label{regu2} \,  We also consider non-weighted eigenvalue
equations and operators for instance, on $S_C=\p B_1(0) \cap C$, cf.Ch.\ref{nee} below. The theory for $L - \lambda \cdot Id$ differs from that of $L_\lambda= L - \lambda \cdot \bp^2 \cdot Id$.\\

For a general shifted skin adapted $L$, the operator $L - \lambda \cdot Id$ is also shifted skin adapted and the non-weighted principal eigenvalue will be $>-\infty$. But $L$ need \emph{not} to be skin adapted for any $\lambda \in \R$ and
thus the Martin boundary of $(H \setminus \Sigma, L - \lambda \cdot Id)$ usually differs from $\widehat{\Sigma}$, for any $\lambda \in \R$. \qed
\end{remark}

\textbf{Minimal Growth Conditions} \quad The notion of  solutions $u >0$ \emph{L}-vanishing in some point  $p \in V \cap \widehat{\Sigma}$, saying that there is a supersolution $w >0$,   such that $u/w(x) \ra 0$, for $x \ra p$, $x \in H \setminus \Sigma$, suggests to view $u$ as a solution of minimal growth towards $p$. Here we make this idea more precise.\\

For domains in $\R^n$, and hence in manifolds, the concept of positive solutions of minimal growth is well-known, cf. [P],Ch.7.3 for a broad discussion of existence and representation results.  We adapt the notion as follows:

\begin{definition} \emph{\textbf{(Minimal Growth)}}\label{mini} \, For any $H \in \cal{H}$  so that $H \setminus \Sigma$ is non-compact
and non-totally geodesic, $L$ a shifted skin adapted operator on $H \setminus \Sigma$ we consider the subcritical operators $L_\lambda$, $\lambda < \lambda^{\bp}_{L,H}$ and solutions of $L_\lambda \, w=0$.\\

A solution $v>0$ on an open set $V \subset \widehat{H}$ with $\widehat{\Sigma} \cap V \n$ has \textbf{minimal growth} towards $\widehat{\Sigma} \cap V$, if for any other such solution with $u \ge v$ near $\p V$:\,$u \ge v$ on $V$.
\end{definition}

The following two results are counterparts to classical results on Euclidean domains, cf.[P],Th.3.7,p.297 and Th.3.9,p.302, for skin adapted operators on $H \setminus \Sigma$.

 \begin{lemma}  \label{hi} \, For $H \in \cal{H}$,  $L_\lambda$, $\lambda < \lambda^{\bp}_{L,H}$ as in \ref{mini} and any basepoint $p_0 \in H \setminus \Sigma$, the minimal Green's function $G(p_0,\cdot)$ of $L_\lambda$ has minimal growth towards $\widehat{\Sigma}$.
\end{lemma}

\textbf{Proof} \,  We choose a sequence of smoothly bounded domains $\overline{D_m} \subset D_{m+1} \subset H
\setminus \Sigma$, $m \ge 0$, with $\bigcup_m D_m = H \setminus \Sigma$ and $p_0 \in D_1$. From the proof of \ref{scal2}, we know that the first eigenvalue $\lambda_m$ of $\delta_{\bp}^2 \cdot L$ for the Dirichlet problem on $\overline{D_m}$ satisfies $\lambda_m >\lambda^{\bp}_{L,H}$.\\

Thus $L_\lambda$, for $\lambda < \lambda^{\bp}_{L,H}$, is again subcritical on each of the $D_m$ and we have a unique Green's function $G_m(p_0,\cdot)$ of $L_\lambda$  on $D_m$ classically vanishing along $\p D_m$. We get $G_m(p_0,\cdot) \nearrow G(p_0,\cdot)$, for $m \ra \infty$,  on each compact subset of $H \setminus \Sigma$. Then the same maximum principle argument we used in \ref{scal2} shows that $G(p_0,\cdot)$  has minimal growth towards $\widehat{\Sigma}$.\qed

Similarly we get the following useful estimates.

 \begin{proposition}  \emph{\textbf{(Minimal Growth Relations)} } \label{w1} \, Let $H \in \cal{H}$,  $L_\lambda$, $\lambda < \lambda^{\bp}_{L,H}$ be as in \ref{mini} with some basepoint $p_0 \in H \setminus \Sigma$. Let $v >0$ solve $L_\lambda \, w=0$ on $H \setminus \Sigma$, so that  $v(p_0)=1$. Then, we have the following estimates
\begin{enumerate}
\item For any neighborhood $U$ of $\widehat{\Sigma}$ with $p_0 \notin \overline{U}$, we have
\begin{equation}\label{g0}
G(p_0,\cdot) \le c \cdot v, \mm{ on } U
\end{equation}
for some $c(U,H,p_0) >0$, independent of $v$.
\item For any couple of open subsets $U \subset \! \subset V \subset \widehat{H}$, with $U \cap \widehat{\Sigma} \n$ and any solution $u >0$ L-vanishing along $V \cap \widehat{\Sigma}$, with $u(p_0)=1$ we have
        \begin{equation}\label{g01}
u \le c^* \cdot v, \mm{ on } U
\end{equation}
for some $c^*(U,V,H,p_0) >0$, independent of $u$ and $v$.
  \end{enumerate}
\end{proposition}

\textbf{Proof} \, For (i) we start with an individual $v$ and show that (\ref{g0}) holds for $c_v:=\max_{x \in \p U} G(p_0,x)/v(x)$, where we note that  $\p U$ is compact and $p_0 \notin \p U$.\\

  For the exhaustion $\overline{D_m} \subset D_{m+1} \subset H \setminus \Sigma$, we used in the previous proof, and from $G_m(p_0,\cdot) < G(p_0,\cdot)$, the maximum principle shows, for $m$ large enough so that $\p D_m \subset U$, we have: $G_m(p_0,\cdot) \le c_v \cdot v$  on $U$. Thus, from $G_m(p_0,\cdot) \nearrow G(p_0,\cdot)$, for $m \ra \infty$, we get  (\ref{g0}) for $v$ and $c_v$.\\

  To show that $c$ can be chosen independently of $v$, we use the standard Harnack inequality for positive solutions of $L_\lambda \, w=0$ on $H \setminus \Sigma$. Since $v(p_0)=1$, $H \setminus \Sigma$ is a connected manifold and $\p U$ is compact, we infer common bounds $0 <m<M$, for all such $v$, so that $m \le \min_{x \in \p U}  v(x) \le \max_{x \in \p U}  v(x) \le M$. Thus, for $k:=\max_{x \in \p U} G(p_0,x)$, and $c(U,H,p_0) :=k/m$ we get  (\ref{g0}) for arbitrary $v$.\\

  For (ii) we merely need to append the boundary Harnack inequality \ref{mbhsqg} applied to $G(p_0,\cdot)$ and $u$.  \qed

\subsubsection{Geometric and Natural Operators} \label{exp}
\bigskip

Now we consider some geometrically relevant examples of skin adapted operators. We start with an operator suggested right from the definition of skin transforms, henceforth
called the \emph{base operator} $L_{\bot}$. It will be used to minorize (the variational integrals for) other skin adapted operators.

\begin{lemma}\label{base} \emph{\textbf{(Base Operator)}}  \quad $L_{\bot}:=-\Delta + |A|^2$ is a skin adapted operator.
\end{lemma}

\textbf{Proof} \quad The claim virtually is a reformulation of the definition of skin transforms [L1],6.1 and 6.2 which asserts the wrapping condition $\bp \ge |A|$ and the validity
of a Hardy inequality. We only need to verify
that, relative to the charts $\psi_p$, we get the claimed regularity features for its coefficients:\\
\[-\Delta u + |A|^2 \cdot u = \sum_{i,j} a_{ij} \cdot \frac{\p^2 u}{\p x_i \p x_j} + \sum_i b_i \cdot \frac{\p u}{\p x_i} + c \cdot u.\]
We recall that the Laplacian written in local coordinates has the form {\small \[\Delta \, u = \frac{1}{\sqrt{det \, g}} \cdot \sum_{i,j} \frac{\p}{\p x_i} \Big(\sqrt{det \, g} \cdot
g^{ij} \cdot \frac{\p \, u}{\p x_j}\Big)\]} We also recall that $|A|^2$ is smooth. Hence the coefficients $a_{ij}$, $b_i$, $c$ are smooth.\\

Since the $\psi_p$ are the geodesic coordinates around $p$ we have
\[a_{ij} = \delta_{ij}, b_i = 0, c=|A|^2(p) \mm{ relative } \psi_p \mm{ in } 0 \in T_pH.\]
Now we can the unfolding correspondence of \ref{mbhs} to see how this transfers to  $(H \setminus \Sigma, d_{\bp^*})$ equipped with $\delta_{\bp^*}^2 \cdot L$ and use
$\bp \ge |A|$ to check that this operator is adapted weakly coercive. \qed

Now we use $L_{\bot} =-\Delta + |A|^2$ to find  other skin adapted operators. Of course, this can be accomplished by adding a smooth function $F:H \setminus \Sigma \ra
\R^{\ge 0}$ with $F \le c \cdot \bp^2$, so that $L_{\bot} + F$ satisfies the adpatedness condition. However, we are interested in operators naturally arising when we use and
study minimal hypersurfaces.\\

One source are operators related to actual geometric and physical variational problems where geometric properties of the ambient space  like its curvature or a given tensor $T$
expressing an physical constraint,  translate into a skin adaptedness of suitably chosen operators on the hypersurface. Here we observe, that skin adaptedness is not only a
matter of the formal shape of an operator. In some cases this depends on further global properties of $H \subset M$ which, in turn, are caused from properties of their ambient
space.

\begin{proposition}\emph{\textbf{(Curvature Constraints)}} \label{scal} \quad  Let $H \in \cal{H}$, $H^n \subset M^{n+1}$, so that $H \setminus \Sigma$ is non-compact
and non-totally geodesic. Then, we have
\begin{enumerate}
 \item If $scal_M \ge 0$, then the conformal Laplacian on $H$ \[L_H: = -\Delta +\frac{n-2}{4 (n-1)} \cdot scal_H\] is skin adapted. In general, $L_H$ is shifted skin adapted.
     \item More generally, let $S$ be any smooth function on $M$ and  $scal_M \ge S$, then the  $S$-conformal Laplacian on $H$
     \[L_{H,S}: = -\Delta +\frac{n-2}{4 (n-1)} \cdot (scal_H - S|_H)\] is skin adapted. Again, in general, $L_{H,S}$ is shifted skin adapted.
     \item The Laplacian $-\Delta_H$ is shifted skin adapted. When $H$ is compact, the principal eigenvalue $\lambda^{\bp}_{-\Delta,H}$ vanishes and the ground state is that of a constant function.
      In particular, $H \setminus \Sigma$ has the Liouville property, that is, all bounded harmonic functions are constant.
     \item The Jacobi field operator $J_H=-\Delta_H - |A|^2-Ric_M(\nu,\nu)$ is shifted skin adapted with  principal eigenvalue $\ge 0$.
\end{enumerate}
\end{proposition}

\textbf{Proof} \quad We start from the Gau\ss-Codazzi equations
\begin{equation} \label{miiq}  |A|^2 + 2 \cdot Ric_M(\nu,\nu)  =  scal_M - scal_H +(tr A_H)^2 , \end{equation}
where $tr A_H$ is the mean curvature of $H$,  $Ric_M(\nu,\nu)$ the Ricci curvature of $M$ for a normal vector $\nu$ to $H$, $scal_H$ and $scal_M$  the scalar
curvature of $H$ and $M$.\\

Since $H$ is minimal, $tr A_H = 0$. Now we observe that the terms $2 \cdot Ric(\nu,\nu)$ and $scal_M$ are remain bounded, whereas $|A|^2$ and $scal_H$ diverge when we
approach $\Sigma$ on $H$. Also  $S|_H$ remains bounded.\\

Therefore, $\bp \ge |A|$ shows that for some constant $k(\bp, S) \ge 1$: \[\bp \ge k \cdot \big|scal_H - S|_H \big|.\]

The argument for \ref{base} above also shows that the operators are adapted to $\bp$.\\

To see (i) and (ii), we prove the validity of the Hardy inequality for $L_{T,S}$. We start from the fact that $H$, as an area minimizing hypersurface, is also stable. That is,
the second variation of the area functional is non-negative.\\

In detail,  infinitesimal variations tangential to $H$ correspond to mere reparametrizations of $H$, we only consider normal fields $f \cdot \nu$, where $\nu$ is the outward
normal vector field over $H \setminus \Sigma$, $f$ is a smooth function on $H$ with $supp \: f \subset H \setminus \Sigma$. Then we have, using again $tr A_H = 0$, cf.  [G],
p.54 {\small\[Area'(f) =  \int_H tr A_H (z) \cdot f(z) \: dVol = 0,\]
\begin{equation} \label{secv}   Area''(f) = \int_{H}|\nabla_H f|^2 - \left( |A|^2 + Ric(\nu,\nu) \right) \cdot f^2 \: dVol \ge 0 \end{equation}}

The Gau\ss-Codazzi equations (\ref{miiq})  give the following equivalent formulation of $Area'' (f) \ge 0$: {\small
\begin{equation} \label{kwsy}  \int_H  f  \cdot  L_H f  \,  dA  = \int_H | \nabla f |^2 + \frac{n-2}{4 (n-1)} \cdot scal_H  \cdot  f^2 \, d A \ge \end{equation}
\[\int_H \frac{n}{2 (n-1)} \cdot  |  \nabla f |^2 + \frac{n- 2}{4 (n-1)}\cdot \left( | A |^2 + scal_M \right) \cdot  f^2\,  d A\]}

Now assume that $scal_M \ge S$, then : {\small
\begin{equation}   \int_H  f  \cdot  L_{H,S} f  \,  dA  = \int_H | \nabla f |^2 + \frac{n-2}{4 (n-1)} \cdot (scal_H - S|_H)  \cdot  f^2 \, d A \ge \end{equation}
\[\int_H \frac{n}{2 (n-1)} \cdot  |  \nabla f |^2 + \frac{n- 2}{4 (n-1)}\cdot \left( | A |^2 + scal_M - S|_H \right) \cdot  f^2\,  d A \ge\]
\[\frac{n- 2}{4 (n-1)}\cdot \int_H f  \cdot  L_{\bot} f   \, d A \ge \tau(\bp,H)  \cdot \frac{n- 2}{4 (n-1)}\cdot \int_H \bp^2 \cdot f^2   \, d A\]}
for $\tau(\bp,H) >0$. \\

For (iii), we note, for the Laplacian, that the condition $\lambda^{\bp}_{-\Delta,H} >-\infty$ is obvious since the variational integral is just the Dirichlet integral $\int_H | \nabla f |^2   \, d A \ge 0$.\\

When $H$ is compact, it is easy to see that the principal eigenvalue equals $0$. Here we use that the codimension of $\Sigma$ is $\ge 2$ and apply the coarea formula e.g. [GMS], Vol1, Ch.2.1.5, Th.3, p.103. \\

Every constant function $v$ solves $\Delta \, v=0$ and thus $v \equiv 1$ can be taken as the ground state for compact $H$. Also for any bounded harmonic function $u$,  $v:=u + \inf_{H \setminus \Sigma} u +1$ is a positive harmonic function. Thus $v$, and therefore $u$, are constant.\\

Finally, we check (iv), for $J_H$. Here we use the minimality of $H$ and get from (\ref{secv}),
\[\int_{H} f \cdot J_H f\, dVol  = \int_{H}|\nabla_H f|^2  -  \left( |A|^2 + Ric(\nu,\nu)  \right) \cdot f^2 \, dVol \ge 0,\]
for any smooth $f$ on $H$ with $supp \: f \subset H \setminus \Sigma$. \qed

 \textbf{Natural Operators} \, A common type of operators are \emph{naturally} assigned to any $H \in {\cal{H}}$, like the Jacobi field operator or the conformal Laplacian we have seen above. This means that given the incarnation $L(H)$ of the operator on $H$  the inherited operator under blow-ups is just the incarnation of this operator on the blow-up space.\\

These operators can be represented using natural coefficients relative e.g. geodesic or other naturally defined coordinates. Examples of coefficient functions
are expressions in terms of inherent curvatures or skin transforms.\\

Formally, we recall from [L1],Def.3.4, that an assignment $F$ of operators $H \mapsto F(H)$, for any $H \in {\cal{H}}$, is called a \emph{natural}, when $F(H)$ commutes with
convergence of the underlying spaces. We briefly call both, the assignment, defined on ${\cal{H}}$, but also the individual operators $F(H)$,
defined on $H \setminus \Sigma$, a \emph{natural operator}.\\

For natural operators the skin adaptedness survives blow-ups. Here we use the regularity theory of area minimizers to see that the approximation of $H$ by blow-up geometries leads to canonical local diffeomorphisms, the $\D$-maps of [L1], Ch.2.3. from smooth in regular regions of the limit geometry to $H$.

\begin{lemma}\emph{\textbf{(Inherited Skin Adaptedness)}} \label{inhnat} \quad  Let $L$ be a natural operator and assume $L(H)$  on some hypersurface $H \in {\cal{H}}$ is (shifted) skin adapted.
 Then $L(N)$ is (shifted) skin adapted  on any Euclidean hypersurface $N$ that arises as a blow-up from $H$ around a singular point in $\Sigma \subset H$ and
 we have: $\lambda^{\bp}_{L,H} \le \lambda^{\bp}_{L,N}$.
\end{lemma}

\textbf{Proof} \quad For any smooth $\phi$ with compact support   $K \subset N \setminus \Sigma_N$ we choose a large scaling factor $\Gamma \gg 1$ so that the $\D$-map
between $\Gamma \cdot H$ and $N$ in a neighborhood of $K$ is very close to an isometry in $C^5$-topology. Then the naturality of $L$ and $\bp$ says that {\small \[\int_H
\phi \circ \D^{-1} \cdot L(H) (\phi \circ \D^{-1}) \, dA/ \int_H \bp_H^2  \cdot  \phi^2\circ \D^{-1} dA \]}
is arbitrarily close to {\small \[\int_N  \phi   \cdot L(N) (\phi) \,  dA/ \int_N \bp_N^2  \cdot  \phi^2  dA\]} upon choosing $\Gamma$ large enough.\\

Thus the eigenvalue of $\delta_{\bp}^2 \cdot L(H)$ on $H$ is a lower estimate for the latter variational integral and, hence, we get the Hardy inequality for $L(N)$ on $N$. The
adaptedness property for $L(N)$ follows from the scaling invariance of the asserted estimates and the adaptedness on $H$, cf. the transformations in the proof of \ref{mbhs}.
\qed

\subsubsection{Schr\"odinger Operators and Scaling Actions} \label{nee}
\bigskip

For Schr\"odinger operators naturally associated to any $H  \in \cal{H}$ we get a transparent representation of some distinguished solutions. We estimate their radial growth in the  cone case and use this to support our understanding of the growth of the Martin kernel for compact $H$.

\begin{definition}\emph{\textbf{(Natural Schr\"odinger Operators)}} \label{spli} \,  A natural shifted skin adapted operator $L$
is called a natural Schr\"odinger operator, when $L(H)$ has the form \[L(H)(u)= -\Delta_H \, u + V_H(x) \cdot \, u \mm{ on } H \setminus \Sigma_H,\] for  some H\"older
continuous function $V_H(x)$, for any given $H \in \cal{H}$.\\
\end{definition}

This class of Schr\"odinger operators  covers the examples we discussed in Ch.\ref{exp}: $L=-\Delta_H$, the Jacobi field operator $J_H=-\Delta_H - |A|^2-Ric_M(\nu,\nu)$ and the $S$-conformal Laplacian on $H$, $L_{H,S} = -\Delta_H +\frac{n-2}{4 (n-1)}
\cdot (scal_H - S|_H)$, for any given $C^\infty$-function  $S$ on $M$.\\

  For an area minimizing cone $C \subset \cal{H}$ the naturality also means that $V_C(t \cdot x)=t^{-2} \cdot V_C(x)$, for any $x \in C \setminus \sigma_C$ and $t >0$.
That is, in this case we can write $V(x)=r^{-2} \cdot V^\times(\omega)$, for $x=(\omega,r) \in S_C \setminus \sigma \times \R^{> 0} = C \setminus \sigma$. Thus we get {\small
\begin{equation} \label{pol} L \, v= -\frac{\p^2 v}{\p r^2} - \frac{n-1}{r} \cdot \frac{\p v}{\p r} -
\frac{1}{r^2} \cdot \Delta_{S_C} v  +\frac{1}{r^2} \cdot  V^\times(\omega) \cdot v=: -\frac{\p^2 v}{\p r^2} - \frac{n-1}{r} \cdot \frac{\p v}{\p r} +\frac{1}{r^2} \cdot L^\times \, v.\end{equation}}

Thus we get an induced operator $L^\times$ on $S_C$, we will see in \ref{fix1} below that $L^\times(S_C)$ is also a natural Schr\"odinger operator.

\begin{remark} \label{ivn} \textbf{(Inheritance versus Naturality)} \, The definition \ref{spli} may leave some space of interpretation. The intended meaning is: the assignment
$H \mapsto L(H)$ on $\cal{H}$ is natural and, independently, on the given $H \in \cal{H}$ it is shifted skin adapted. Thus the condition does not assert it is shifted skin adapted for \emph{all} $H \in \cal{H}$. However, due to \ref{inhnat}, the naturality always implies that with $L(H)$,  $L(N)$, for a blow-up $N$ of $H$ around a singular point in $\Sigma \subset H$, is shifted skin adapted with $\lambda^{\bp}_{L,H} \le \lambda^{\bp}_{L,N}$.\qed
\end{remark}

\textbf{Scaling Actions} \, On an area minimizing cone equipped with a natural and skin adapted  operator $L$, both, the operator and solutions
of $L \, v =0$ reproduce under scalings of the cone up to constant multiples.\\

 More concretely, we write $L$ on $C = \R^{\ge 0} \times S_C$ in geodesic coordinates $x_1=r$,
$x_2,...x_n$, so that $x_2,...x_n$ locally parameterize $S_C=\p \cap B_1(0)$ {\small \begin{equation}  -L(u) =  \sum_{i,j} a_{ij} \cdot \frac{\p^2 u}{\p s_i \p s_j} + \sum_i b_i
\cdot \frac{\p u}{\p s_i} + c \cdot u.\end{equation}} When $L$ is natural this means that for any $\eta>0$:
\[a_{ij}(\eta \cdot x) = a_{ij}(x), b_i(\eta \cdot x) = \eta^{-1} \cdot b_i(x) \mm{ and } c(\eta \cdot x) = \eta^{-2} \cdot c(x)\]
and thus the chain rule shows that for any function $u(x)$ that solves $L \, v =0$ the scaled versions $u(\eta \cdot x)$ also solve this equation.\\

In particular, up to multiples, the Green's function and the set of minimal solutions of $L \, u = 0$ are reproduced under composition with the scaling map
\[ S_\eta: C \ra C, \mm{ given by } x \mapsto \eta \cdot x, \mm{ for } \eta \in (0,\infty).\]
That is we consider the map $u \mapsto u \circ S_\eta$ and regauge the values of the resulting functions $u \circ S_\eta$ in some base point $p \in C \setminus \sigma$ to $1$.
In this fashion we define a scaling action $S^*_\eta$ on the Martin boundary.

\begin{lemma}\emph{\textbf{(Attractors and Fixed Points of $\mathbf{S^*_\eta}$)}} \label{fix} \, Let $C$ be an area minimizing cone,
$L$ a natural skin adapted operator on $C$. Then we have
\begin{enumerate}
\item The scaling action $S^*_\eta$ on $\p_M (C,L)$ has exactly two fixed points \[[0]:= \{0\}/\sim \mm{ and } [1]:= \{1\}/\sim \mm{ on } \p_M (C,L) = [0,1]\times
    \Sigma_{S_C}/\sim.\] They correspond to the tip $0$ and the point at infinity $\infty$.
\item   For  $z \in \p_M (C,L) \setminus \{0,1\}$  we find
\[S^*_\eta(z) \ra [1], \mm{ for } \eta \ra \infty \, \mm{ and } \,  S^*_\eta(z) \ra [0], \mm{ for } \eta \ra 0.\]
\item More generally, we have for any given $z \in C \setminus \sigma$ {\small \[\frac{G(x,S_\eta(z))}{G(p,S_\eta(z))} \ra [1], \mm{ for } \eta \ra \infty \,\,\mm{ and }\,\, \frac{
    G(x,S_\eta(z))}{G(p,S_\eta(z))} \ra [0], \mm{ for } \eta \ra 0\]} locally uniformly for $x \in C \setminus \sigma$.
\end{enumerate}
\end{lemma}

{\bf Proof} \quad  This readily follows from the way the pole of $G(\cdot,S_\eta(z))$ shifts under these scaling operations.\qed

Now we give a description of the two fixed point solutions in $\p_M (C,L)$ to build inductive decomposition schemes for solutions on the original hypersurface $H$, despite the
mentioned fact that there is no proper transition of solutions in terms of Martin integrals.

\begin{proposition}\emph{\textbf{(Separation of Variables)}} \label{fix1} \quad Let $C$ be an area minimizing cone, $L$ a natural Schr\"odinger operator. Then we have for
 the  skin adapted operator $L_\lambda=L - \lambda \cdot \bp^2 \cdot Id$, for  $\lambda < \lambda^{\bp}_{L,C}$:
\begin{itemize}
\item The two fixed points $[0], [1] \in \p_M (C,L_\lambda)$ viewed as functions $\Psi_-=[0]$ and $\Psi_+=[1]$  on $C \setminus \sigma$ can be written as
\[\Psi_\pm(\omega,r) = \psi(\omega) \cdot r^{\alpha_\pm},\]
for $(\omega,r) \in S_C \setminus \sigma \times \R$, with {\small $\alpha_\pm = - \frac{n-2}{2} \pm \sqrt{ \Big( \frac{n-2}{2} \Big)^2 + \mu_{C,L^\times_\lambda}}$}.
 \item $\mu_{C,L^\times_\lambda} > - (\frac{n-2}{2})^2$ is the non-weighted principal eigenvalue and  $\psi(\omega) >0$  the ground state of an associated natural
     Schr\"odinger operator $L_\lambda^\times$:
 \[L_\lambda^\times(v)(\omega)= - \Delta_{S_C} v(\omega) + \big(V^\times(\omega) - \lambda \cdot (\bp^\times)^2(\omega) \big)\cdot v(\omega),\]
    defined on $S_C \setminus \Sigma_{S_C}$, where $\bp(x)=r^{-1} \cdot \bp^\times(\omega)$, for $x=(\omega,r) \in C \setminus \sigma$. That is,  we have
    \[L_\lambda^\times \, \psi =  \mu_{C,L^\times_\lambda} \cdot  \psi, \mm{ on } S_C \setminus \Sigma_{S_C}.\]
\end{itemize}
\end{proposition}

{\bf Proof} \quad We proceed in three steps: we show  $\Psi_\pm$ can written as a product $\Psi_\pm(\omega,r) = \psi_\pm(\omega) \cdot r^{\alpha_\pm}$. Then we determine the $\psi_\pm$ and $\alpha_\pm$ for some inner approximation of by regular subcones of $C$ and finally show that the resulting values converge to the $\psi_\pm$ and $\alpha_\pm$ on $C$.\\

\textbf{Product Shape} \, We first restrict the two fixed points $\Psi_\pm \in \p_M (C,L_\lambda)$ to a regular ray $\Gamma_v = \R^{>0} \cdot v \in C \setminus \sigma$, for some $v
\in S_C$, and consider the map $\Psi_\pm|_{\Gamma_v}: \R^{>0} \cdot v \ra \R^{>0}$. We view such a restriction as a function $f:\R^{>0} \ra \R^{>0}$.\\

We know that $\Psi_\pm$ reproduces under scalings, up to a constant: for any $\eta >0$, there is constant $c_\eta>0$ so that
\[f(\eta \cdot x)=c_\eta \cdot f(x), \mm{ for all } x \in \R^{>0}.\]

From this it follows that $f$ is a monomial: $f(x)=a \cdot x^b$, for some constants $a(v) > 0, b(v) \in\R$.  This argument applies to any regular ray $\Gamma_v$.\\

Next we consider the Harnack inequality for $L_\lambda$ on a ball $B_{2 \cdot R}(v) \subset C \setminus \sigma$, for some $R>0$. We get, for any positive solution $u$ of
$L_\lambda \, v =0$, the Harnack inequality
\begin{equation}\label{h1}
 \sup_{B_R(v)}\, u \le c(L_\lambda,v,R)  \cdot \inf_{B_R(v)}\, u,
\end{equation}  for some constant independent of $u$. Now the point is that the scaling symmetry of $C$ and the naturality of $L_\lambda$ imply that
 the same constant $c$ can still be used in the Harnack inequality after scalings around the tip $0$.
\begin{equation}\label{h2}
 \sup_{B_{s \cdot R}(s \cdot v)}  u \le c(L_\lambda,v,R)  \cdot \inf_{B_{s \cdot R}(s \cdot v)} u, \mm{ holds for any } s >0.
\end{equation}

This implies that for all rays $\Gamma_w$ passing through $B_R(v)$ the exponent $b(w)$ equals $b(v)$. Since $S_C \setminus \Sigma_{S_C}$ is connected, this shows that $b(v)$ is constant on $S_C \setminus \Sigma_{S_C}$.\\

Thus $\Psi_\pm$ can written as $\Psi_\pm(\omega,r) = \psi_\pm(\omega) \cdot r^{\alpha_\pm},$ for $(\omega,r) \in \p B_1(0) \cap C \setminus
\sigma \times \R^{> 0}$ and some  $C^{2,\beta}$- regular function $\psi_\pm$ on $\p B_1(0) \cap C \setminus \sigma \times \R^{> 0}$, for some $\beta>0$, cf.\ref{regu}.1.\\

When we reinsert $\Psi_\pm(\omega,r) = \psi_\pm(\omega) \cdot r^{\alpha_\pm}$ into the equation $L_\lambda \, w =0$, written in polar coordinates (\ref{pol}), we observe,
from a separation of variables, that the $\psi_\pm$ solve the following equations on  $S_C \setminus \Sigma_{S_C}$:  {\small \begin{equation}\label{sph}  L_\lambda^\times \,
v := - \Delta_{S_C} v(\omega) + \big(V^\times(\omega) - \lambda \cdot (\bp^\times)^2(\omega) \big)\cdot v(\omega) = (\alpha_\pm^2 + (n-2)\cdot  \alpha_\pm) \cdot v,
\end{equation}}
$L_\lambda^\times$ is again a natural Schr\"odinger operator and adapted to the skin transform $\bp_C|_{S_C}$.\\

\textbf{Inner Regular Approximation} \, We use an approximation by Dirichlet eigenvalue problems to show
\begin{itemize}
    \item  {\small $\alpha_-^2 + (n-2)\cdot  \alpha_- = \alpha_+^2 + (n-2)\cdot  \alpha_+ > - ((n-2)/2)^2$}
    \item  {\small $\mu_{C,L^\times_\lambda} := \alpha_\pm^2 + (n-2)\cdot  \alpha_\pm $}
    is the non-weighted  \emph{principal eigenvalue} of $L_\lambda^\times$. Also $L_\lambda^\times$ is shifted skin adapted.
    \item $\psi_-=\psi_+$  and $\psi:=\psi_\pm$ is the corresponding \emph{ground state} of $L_\lambda^\times$.
\end{itemize}

To this end we choose smoothly bounded domains $G_i \subset \overline{G_i} \subset S_C \setminus \Sigma_{S_C}$, $G_i \subset G_{i+1}$ and $\bigcup_i G_i = S_C
\setminus \Sigma_{S_C}$. We consider the positive solutions of $L_\lambda \, v=0$ on the cone $C(G_i) \subset C$ over $G_i$ with vanishing boundary value
along $\p C(G_i) \setminus \{0\}$.\\

Now we apply \ref{mgr}, case (i),  to  $L_\lambda$ on $C(G_i)$. Note that  we have some $a_{C(G_i)}>0$ so that:\, $dist(z,\Sigma) \le a_{C(G_i)} \cdot \delta_{\bp}(z) \le a_{C(G_i)} \cdot L_{\bp} \cdot  dist(z,\Sigma)$. Also one readily checks the adaptedness conditions.\\

Therefore, we get that the Martin boundary of $L_\lambda$
equals $(\p C(G_i) \setminus \{0\}) \cup \{[0],[1]\}$ and the two fixed point solutions $\Psi_\pm[i]$ of $L_\lambda \,v$ vanishing along $\p C(G_i) \setminus \{0\}$, alternatively
labelled $[0][i],[1][i]$, are again of the form $\Psi_\pm[i](\omega,r) = \psi_\pm[i](\omega) \cdot r^{\alpha[i]_\pm}$.\\

Again, we insert $\Psi_\pm[i](\omega,r)$ into the equation $L_\lambda \, w =0$, written in polar coordinates (\ref{pol}), we find, that the $\psi_\pm[i]$ solve \[
L_\lambda^\times \, v = - \Delta_{S_C} v(\omega) + \big(V^\times(\omega) - \lambda \cdot (\bp^\times)^2(\omega) \big)\cdot v(\omega) = (\alpha_\pm[i]^2 + (n-2)\cdot
\alpha_\pm[i]) \cdot v,\] For $G_i$ we apply the spectral theory for bounded domains and observe that the \emph{positive} eigenfunctions $\psi_\pm[i]$ must equal the uniquely determined first Dirichlet eigenfunction $\psi[i]$ for the first eigenvalue $\mu[i]$ of $L_\lambda^\times$. Thus we have
\[\psi[i]=\psi_-[i]=\psi_+[i] \,\mm{ and } \, \mu[i]=\alpha_-[i]^2 + (n-2)\cdot \alpha_-[i]=\alpha_+[i]^2 + (n-2)\cdot  \alpha_+[i],\]
hence $\mu[i]\ge- ((n-2)/2)^2$ and we have $\mu[i]\ge \mu[i+1]$, as is seen from the variational characterization of these eigenvalues, since the space of admissible test functions on $G_i$ is a subset of the corresponding function space over $G_{i+1}$.\\

When suitably normalized, the first Dirichlet eigenfunction $\psi[i]>0$ of $L_\lambda^\times$ on the $G_i$, for $i \ra \infty$, uniquely converge $C^3$-compactly on $S_C
\setminus \Sigma_{S_C}$ to an eigenfunction
  $\psi^*>0$ with eigenvalue $\mu^* = \lim_{i \ra \infty} \mu[i]\ge- ((n-2)/2)^2$. We also define the limits $\alpha_\pm^* := \lim_{i \ra \infty} \alpha_\pm[i]$.\\

 From this we observe as in \ref{scal2} that $\psi^*>0$ is the non-weighted ground state of $L_\lambda^\times$ on $S_C \setminus \Sigma_{S_C}$ for the eigenvalue $\mu^* >-\infty$. We observe that
 $L_\lambda^\times$, and hence $L^\times$, are shifted skin adapted. Namely, $\bp^\times > c_{S_C} >0$ and thus with the non-weighted principal eigenvalue $\mu^*$ of $L_\lambda^\times$, the principal eigenvalue of $\delta^2_{\bp^\times} \cdot L_\lambda^\times$ remains finite.\\

 Therefore, $L_\lambda^\times$ is a natural Schr\"odinger operator.\\

\textbf{Comparison Arguments} \, Now we compare these solutions with $[0], [1] \in \p_M (C,L_\lambda)$. From the fact that there are no positive solutions for $L_\lambda^\times \,  v = \mu \cdot v$, when
$\mu <\mu^*$, we see that
\begin{equation}\label{aq}
\alpha_\pm^2 + (n-2)\cdot  \alpha_\pm\ge\mu^*\ge- ((n-2)/2)^2
\end{equation}
In turn, the solution $\Psi_+=\psi(\omega) \cdot r^{\alpha_+}$ $L$-vanishes along $\widehat{\sigma_C} \setminus \{\infty\}$, in particular in $0$, and this implies the estimate
\begin{equation}\label{aqq}
\alpha_+^* \ge \alpha_+\mm{ and thus  }\mu^* \ge  \alpha_\pm^2 + (n-2)\cdot  \alpha_\pm
\end{equation}
From (\ref{aq}) and  (\ref{aqq})  we conclude \[\mu^* =  \alpha_\pm^2 + (n-2)\cdot  \alpha_\pm \mm{ and therefore }
\alpha_\pm=\alpha^*_\pm.\] Since $[0] \neq [1] \in \p_M (C,L_\lambda)$, we further infer that \[\alpha_- <  -(n-2)/2 < \alpha_+ \mm{ and } \alpha_\pm^2 + (n-2)\cdot
\alpha_\pm > - ((n-2)/2)^2.\] Also, the $\psi_\pm$ belong to the eigenvalue $\alpha_\pm^2 + (n-2)\cdot \alpha_\pm=\mu^*$. Since the ground state is uniquely determined,
we have $\psi:=\psi_-=\psi_+=\psi^*$  is the ground state of $L_\lambda^\times$. \qed

\begin{remark} \label{recc} \textbf{(Symmetry versus  Uniformity)} \quad  At first sight, the previous result seems to be plausible as a consequence of the cone symmetry alone. However, it also substantially rely its (skin) uniformity. To illustrate this point, we discuss an example  due to Ioffe and Pinsky [IP]. They studied the following domains in $\R^n$, $n \ge 3$: for some smooth
$f:\R^{\ge 0} \ra \R^{>0}$ and $F:\R^{n-1}\ra \R$ with $F(x):=f(|x|)$  choose the  domain $D_f:=\{(x_1,..,x_n) \in \R^n \,|\, |x_1|< F(x_2,..,x_n) \}$.\\

[IP] contains an analysis of Martin boundary of $D_f$ for the Laplacian, for instance, for the following two standard types of functions:
\[f_1(t)= c_1 \cdot t + c_2, \, f_2(t)=c_1 \cdot (t+1)^\gamma + c_2 , \mm{ for constants } c_i >0, \gamma \in (0,1).\]
Note that $D_{f_1}$ contains the open cone $D_{c_1 \cdot t} \subset \R^{n}$, while  $D_{f_2}$ compresses to the hyperplane $\R^{n-1} \times \{0\}$ when we scale $D_{f_2}$
by constants $s \ra 0$.\\

In both cases $\p  D_{f_i}$ is a Lipschitz boundary and thus one has BHPs, with the possible exception for the point at infinity, and [IP],Th.1, actually says:
\begin{itemize}
    \item the \emph{Martin boundary  at infinity} of $D_{f_1}$ equals precisely one point,
    \item the \emph{Martin boundary  at infinity}  of $D_{f_2}$ equals $S^{n-1}$.
\end{itemize}
We may retrieve the first result, for $D_{f_1}$, from the fact that it is a uniform domain, as one readily checks. $D_{f_2}$ is not uniform as can be seen from
the mentioned scaling compression to a hyperplane and, with hindsight, from the fact that its Martin boundary at infinity contains more than one point. \qed
\end{remark}

 Along the lines of \ref{fix1}, we also find all fixed points of the scaling action.

\begin{corollary}\emph{\textbf{(Non-Minimal Fixed Points)}} \label{nmi} \, Let $C$ be an area minimizing cone, $L$ a natural Schr\"odinger operator. Then we have for
 the  skin adapted operator $L_\lambda=L - \lambda \cdot \bp^2 \cdot Id$, for  $\lambda < \lambda^{\bp}_{L,C}$, and any $\alpha \in (\alpha_-,\alpha_+)$:
 \begin{itemize}
\item There are \textbf{non-minimal} solutions of $L_\lambda \, u = 0$ reproducing under the scaling action $C \setminus \sigma$:
\begin{equation}\label{nmieq}
\psi(\omega) \cdot r^{\alpha}, \mm{  where }
\end{equation}
\[\psi > 0 \mm{ is a solution of } L_\lambda^\times \, v  = (\alpha^2 + (n-2)\cdot  \alpha) \cdot v, \mm{ on } S_C \setminus \Sigma_{S_C}.\]
\item The collection of these solutions and the minimal solutions $\Psi_\pm$ of \ref{fix1} make up the complete set of all fixed point solutions (up to multiples).
\item The operator $L_{\lambda,\alpha}^\times := L_\lambda^\times - (\alpha^2 + (n-2)\cdot  \alpha)$ is skin adapted. \item Therefore, $\p_M
    (S_C,L_{\lambda,\alpha}^\times) \cong \Sigma_{S_C}$ and any positive non-minimal fixed point solutions of $L_\lambda \, u = 0$ can written
    as\begin{equation}\label{nmm} \Psi_{\alpha,\widetilde{\mu}}(\omega,r) := \psi_{\alpha,\widetilde{\mu}}(\omega) \cdot r^{\alpha}:=\int_{\Sigma_{S_C}}
    \widetilde{k}(x;y) \, d \widetilde{\mu}(y)
\cdot r^\alpha
\end{equation}
where $\widetilde{k}(x;y)$ denotes the Martin kernel of $L_{\lambda,\alpha}^\times$ on $S_C \setminus \sigma$ and $\widetilde{\mu}$  is a (unique) finite Radon measure
on $\Sigma_{S_C}$ associated to $\Psi_{\alpha,\widetilde{\mu}}$.\\
\end{itemize}
\end{corollary}

{\bf Proof} \quad Most of this follows readily from the proof of \ref{fix1}. What remains is to see that for any $\alpha \in (\alpha_-,\alpha_+)$, $\alpha_\pm=
\alpha_\pm(\lambda)=- \frac{n-2}{2} \pm \sqrt{ \Big( \frac{n-2}{2} \Big)^2 + \mu_{C,L^\times_\lambda}}$, the operator is
$L_{\lambda,\alpha}^\times := L_\lambda^\times - (\alpha^2 + (n-2)\cdot  \alpha)$ is skin adapted.\\

For this, we use $\lambda < \lambda_C$.  For any $\ve >0$, we find a small $\delta>0$ with $\lambda+\delta < \lambda_C$ so that
\[|\mu_{C,L^\times_\lambda}-\mu_{C,L^\times_{\lambda+\delta}}| \le \ve.\]
We infer, that for any $\alpha \in (\alpha_-(\lambda),\alpha_+(\lambda))$ there are sufficiently small $\ve$ and $\delta$ so that the non-weighted principal eigenvalue of
$L_{\lambda+\delta}^\times - (\alpha^2 + (n-2)\cdot  \alpha)$  and, thus, that of $\delta_{\bp^\times}^2 \cdot (L_{\lambda+\delta}^\times - (\alpha^2 + (n-2)\cdot \alpha))$,
is non-negative.\\

 This shows that any  $\alpha \in (\alpha_-(\lambda),\alpha_+(\lambda))$, there is some $\delta(\alpha) >0$ so that the principal eigenvalue of $\delta_{\bp^\times}^2 \cdot
L_{\lambda,\alpha}^\times$ is $\ge \delta(\alpha)$, that is, $L_{\lambda,\alpha}^\times$ is skin adapted. \qed

Since any positive solution of $L \, w=0$ for natural skin adapted $L$ can be written as a convex combination of minimal solutions, the Martin integral,
the minimal solutions  $\Psi_\pm$ still play their role of attractors for the scaling action relative general solutions. From this we also see that minimal growth of solutions  $L\, w=0$ is a natural property conserved under convergence of area minimizing hypersurfaces:
\begin{proposition} \emph{\textbf{(Recovery Processes)}}  \label{rp} Let $L$ be a natural Schr\"odinger  operator, $H \in {\cal{SH}}$ and $C \in {\cal{SC}}$. Then we have:
\begin{enumerate}
  \item  For any solution $u >0$ on $C \setminus \sigma_C$ of $L(C) \, w=0$, we get a locally uniform convergence:
\[S^*_\eta(u+\Psi_-) \ra \Psi_- \mm{ for } \eta \ra \infty \, \mm{ and } \,  S^*_\eta(u+\Psi_+) \ra \Psi_+, \mm{ for } \eta \ra 0.\]
where we think of the functions $S^*_\eta(u+\Psi_\pm)$ and the limit $\Psi_\pm$ being (re)normalized to value $1$ in some common basepoint $p \in C \setminus \sigma_C$.
  \item For any solution $u >0$ on $H \setminus \Sigma_H$ of $L(H) \, w=0$, $L$-vanishing along $V \cap \widehat{\Sigma}$, for some open $V \subset H$ and $y \in V \cap \widehat{\Sigma}$, the induced* solutions on  tangent cones $C$ in $y$ $L$-vanish along $\sigma_C$.
\end{enumerate}
\end{proposition}
*A detailed discussion of the definition and properties of induced solutions is given in [L2],Ch.2.1 and 2.2.\\

{\bf Proof} \quad \textbf{(i)}  follows from \ref{fix}(ii) viewing $u$ as a Martin integral and the relation $\alpha_- <  -(n-2)/2 < \alpha_+$ for the extremal growth rates of $\Psi_\pm$, which characterizes these minimal functions. Observe that scaling of $\R^+$ by some $\eta >0$ and $\alpha \in (\alpha_-, \alpha_+)$  gives the locally uniform convergence:
\[(\eta \cdot r)^{\alpha} /(\eta \cdot r)^{\alpha_-} \ra 0,\mm{ for } \eta \ra 0\mm{ and } (\eta \cdot r)^{\alpha} /(\eta \cdot r)^{\alpha_+} \ra 0,\mm{ for } \eta \ra \infty,\]
from  $(\eta \cdot r)^{\alpha}  = \eta^{\alpha} \cdot r^{\alpha}$.\\

\textbf{(ii)}\, We consider the minimal growth relations \ref{w1}(ii) and observe that for a solution $u >0$ that $L$-vanishes along $V \cap \widehat{\Sigma}$, we can infer that for the minimal function $k(\cdot ; y)$, associated to the Dirac measure at $y$,  with $y \in  U \cap \widehat{\Sigma}$ with  $U \subset \! \subset V$, there is a constant $c_y >0$ so that
 \begin{equation}\label{iii}
u \le c_y \cdot k(\cdot ; y), \mm{ on } U
\end{equation}
Under increasing scalings around $y$, and suitable (re)normalizations, the solutions $v_1=u$ and $v_2=k(\cdot ; y)$ of $L(H) \, w=0$ induce positive solutions $v^*_i$ of $L(C) \, w=0$ on each of the tangent cones $C$ at $y$. Under this blow-up process this inequality, with suitably adjusted constant, is inherited to any ball $B_r(0) \cap C$, $r>0$.\\

In turn, on $C$, one verifies that $v^*_2$ is again the limit of Green's functions $G_C(x,p_n)$, $p_n \ra 0$. This follows from the general characterization of $G_C$, cf.[An3],Ch.1. But this general argument does not imply that the $G_C$ are \emph{minimal} Green's functions. Thus, at this stage, the Radon measure associated to $v^*_2$ could be the Dirac measure at $0 \in C$, but accompanied by some  non-trivial density outside $0$. However, for a natural Schr\"odinger  operator it is precisely the Dirac measure at $0$.\\

Namely, we observe from part (i) that any such $v^*_2$ converges to $k(\cdot ; 0)$ for increasingly large scalings around $0$. But $v^*_2$ may already be viewed as the result of such a blow-up process and thus it already was $k(\cdot ; 0)$.\\

Now we recall that $v^*_2=k(\cdot ; 0)$, $L$-vanishes towards $\widehat{\sigma_C} \setminus \{0\}$. This shows that $v^*_1$  $L$-vanishes towards $\sigma_C$.   \qed

\subsubsection{Conformal Laplacians} \label{ard}
\bigskip

Now we focus on a subclass of natural Schr\"odinger operators one typically encounters in applications in scalar curvature geometry and derive some (more) explicit growth estimates.

\begin{proposition}\emph{\textbf{(Eigenvalue Estimates for $J_C$ and $L_C$)}} \label{evee} \, Let $C$ be a singular area minimizing cone. Then we get the following estimates for the Jacobi field operator $J_C$ and the conformal Laplacian $L_C$:
\begin{itemize}
  \item The principal eigenvalue $\lambda_{J_C,C}$ of $J_C$ is non-negative  and we have
 \[ \mu_{C,(J_C)^\times_\lambda} \ge - \left(\frac{n-2}{2}\right)^2, \mm{ for any } \lambda \le \lambda_{J_C,C}.\]
  \item   There are constants $\Lambda_n > \lambda_n >0$ depending only on $n$ so that:
For any $\lambda \le \Lambda_n$, $(L_C)_\lambda$ is skin adapted and we have \[\mu_{C,(L_C)^\times_\lambda} \ge - 1/4 \cdot \left(\frac{n-2}{2}\right)^2,\mm{ for any }\lambda \le \lambda_n, \mm{ and thus }\]
\[\alpha_+ \ge - (1- \sqrt{3/4}) \cdot \frac{n-2}{2}, \,\, \alpha_- \le - (1+ \sqrt{3/4}) \cdot \frac{n-2}{2},\mm{ for }\lambda \le \lambda_n.\]
Moreover, for $\lambda \in (0,\lambda_n]$, we have $\mu_{C,(L_C)^\times_\lambda} <-\eta_\lambda <0$ and hence:
 {\small \[-\vartheta_\lambda > \alpha_+ \ge - (1- \sqrt{3/4}) \cdot \frac{n-2}{2} >-\frac{n-2}{2} >- (1+ \sqrt{3/4}) \cdot \frac{n-2}{2} \ge \alpha_- > \vartheta_\lambda -(n-2).\]}
 for some constants $\eta_\lambda, \vartheta_\lambda>0$, depending only on $\lambda$ and $n$.
\end{itemize}
\end{proposition}

{\bf Proof} \quad From \ref{scal}(iv) we know that $\lambda_{J_C,C} \ge 0$. For $\lambda <\lambda_{J_C,C}$, $J_C$ is skin adapted and thus we have the lower bound $\mu_{C,(J_C)^\times_\lambda} >- \left(\frac{n-2}{2}\right)^2$. From \ref{fix1}, we also find a common upper bound for $\mu_{C,(J_C)^\times_\lambda}$, for all $\lambda$ with $\lambda<  \lambda_{J_C,C} \le \lambda + \ve$, for any $\ve >0$ given.\\

Next we recall that the ground state $v>0$ is uniquely determined.  It follow from standard elliptic theory applied to the given set of eigenvalue equations that $v$ is the unique $C^3$-limit of solutions $v_n >0$ of $(J_C)_{(\lambda_{J_C,C}-1/n)}\, w=0$, $n \in \Z^{>0}$, (note that in this case we have locally uniform bounds on the coefficients) and hence we find that $\mu_{C,(J_C)^\times_{\lambda_{J_C,C}}} \ge - \left(\frac{n-2}{2}\right)^2$. \\

For $L_C$ we first notice that there is common $\Lambda_n>0$ for all area minimizing cones in $\R^{n+1}$ so that $(L_C)_{\lambda}$ is skin adapted for all $\lambda \le \Lambda_n$. For each of the $C$ this is just \ref{scal}(i). The uniform estimate then follows from the compactness of the space of singular cones.\\

Now we estimate the variational integral for $\mu_{C,(L_C)^\times_\lambda}$ for any smooth function $f$ compactly supported in $S_C \setminus \Sigma_{S_C} =\p B_1(0) \cap C \setminus \sigma_C$.  To this end, we use
 \begin{itemize}
   \item  the estimate  $\mu_{C,(J_C)^\times_0} \ge - \left(\frac{n-2}{2}\right)^2$,
   \item the Hardy inequality (\ref{hh}) for the skin transform $\bp|_{S_C}$ where we can choose the constant $\tau(\bp,C)=\tau_n >0$ uniformly for all cones (this follows from the naturality of $\bp$ as discussed in [L1],Prop.3.4 for an explicit sample of $\bp$),
   \item the identity $scal_C\equiv -|A|^2$, valid for minimal hypersurfaces in any Ricci flat space, from the Gau\ss-Codazzi equation (\ref{miiq}):
 \end{itemize}
\begin{equation}   \int_{S_C}  f  \cdot  (L_C)^\times_{\lambda} f  \,  dA  = \int_{S_C} | \nabla f |^2 + \left(\frac{n-2}{4 (n-1)} \cdot scal_C|_{S_C} - \lambda \cdot \bp|_{S_C}^2\right)  \cdot  f^2 \, d A \ge \end{equation}
\[\int_{S_C} | \nabla f |^2 - \frac{n-2}{4 (n-1)} \cdot |A|^2 \cdot f^2 dA - \max\{0,\lambda/\tau_n\} \cdot \int_{S_C}|\nabla f|^2  + |A|^2 \cdot f^2 dA.\]
For $\lambda/\tau_n < \frac{1}{4 (n-1)}$, and $|f|_{L^2}=1$, this can be lower estimated by
\[\ge 1/4 \cdot  \int_{S_C} | \nabla f |^2 -  |A|^2 \cdot f^2 dA \ge  1/4 \cdot  \mu_{C,(J_C)^\times_0} \ge - 1/4 \cdot \left(\frac{n-2}{2}\right)^2.\]
Hence, for $\lambda_n:=1/2 \cdot \min \{\Lambda_n, \frac{\tau_n}{4 (n-1)}\}$, we have
\[\mu_{C,(L_C)^\times_\lambda} \ge - 1/4 \cdot \left(\frac{n-2}{2}\right)^2,\mm{ for }\lambda \le \lambda_n.\]
Finally, for $\lambda \in (0,\lambda_n]$, we observe, since the codimension of $\Sigma_{S_C}$ is  $>2$ and $scal_C \le 0$, we can find a function $f \in C_0^\infty(S_C \setminus \Sigma_{S_C})$ equal to $1$ outside a small enough neighborhood of $\Sigma_{S_C}$ and   equal to $0$ very close to $\Sigma_{S_C}$ so that
\[ \int_{S_C} | \nabla f|^2 + \left(\frac{n-2}{4 (n-1)} \cdot scal_C|_{S_C} - \lambda \cdot \bp|_{S_C}^2\right)  \cdot  f^2 \, d A < -1/2 \cdot \int_{S_C} \lambda \cdot \bp|_{S_C}^2  \cdot  f^2 \, d A  < 0.\]
Due to the naturality of skin transforms and the \emph{compactness} of the space of all \emph{singular} area minimizing cones $C \subset \R^{n+1}$ in flat norm topology and in compact $C^5$-topology outside the singular sets, we know that there common positive lower bound on $\int_{S_C}\bp|_{S_C}^2 \, d A$ for all such $C$.\\

This means,  there are constants $\eta_\lambda, \vartheta_\lambda>0$, depending only on $\lambda>0$ and $n$, so that
\[\mu_{C,(L_C)^\times_\lambda} <-\eta_\lambda<0 \mm{ and, therefore, } 0>-\vartheta_\lambda>\alpha_+ \mm{ and }\alpha_- > \vartheta_\lambda -(n-2)>  -(n-2).\]
for any singular area minimizing cone $C$. \qed

In inductive cone reduction arguments we use iterated blow-ups around a series of singular points. This operation suggests to consider also some variants of these operators.\\

Iterative blow-ups are processes where we first blow-up around some $p_0 \in \Sigma_H$ and get a tangent cone $C_1$. While scaling around $0 \in \sigma_{C_1}$ just reproduces $C_1$, blow-ups around singular points $p_1\neq 0 \in \sigma_{C_1}$ generate some area minimizing cones $C_2$ that can be written as a Riemannian product $C_2=\R \times C_3$, for a lower dimensional area minimizing cone $C_3$.\\

This way we encounter cones $C^n\subset \R^{n+1}$ which can written as a Riemannian product $\R^{n-k} \times C^k$, where $C^k \subset \R^{k+1}$ is a lower dimensional area minimizing cone.\\

 In this case we observe that the minimal function $\Psi_+(\omega,r)$ in Martin boundary for $(L_{C^n})_\lambda$ on $C^n$ shares the $\R^{n-k}$-translation symmetry with the underlying space $\R^{n-k} \times C^k$ since  $\Psi_+(\omega,r)$  is uniquely determined. From this, we observe that $\Psi_+(\omega,r)|_{\{0\} \times C^k}$ satisfies the equation \[(L_{C^k,n})_\lambda \Psi_+(\omega,r)|_{\{0\} \times C^k}=0.\]
  where \[L_{C^k,n} :=-\Delta - \frac{n-2}{4 (n-1)} \cdot |A|^2,\,  \mm{ for }
 n \ge k \mm{ and } (L_{C^k,n})_\lambda  = L_{C^k,n}-\lambda  \cdot  \bp^2\] Thus $L_{C^k,n}$ is a dimensionally shifted version of the conformal Laplacian on the  cone $C^k$. $\Delta$, $|A|^2$ and $\bp^2$, are the entities intrinsically  defined on $C^k$. The dimensional shift comes from using $\frac{n-2}{4 (n-1)}$ in place of $\frac{k-2}{4 (k-1)}$. The next two results describe the analysis of these operators.

 \begin{proposition}\emph{\textbf{(Dimensionally Shifted $L_C$)}} \, Let $C^k \subset \R^{k+1}$ be a singular area minimizing cone and $n \ge k$. Then we have
\begin{itemize}
  \item $L_{C,n}$ is skin adapted and  the principal eigenvalue can be uniformly lower estimated for all $k$-dimensional cones by a positive constant $\Lambda^*_k >0$ independent of
 $n \ge k$.
  \item   There is a constant $\lambda^*_k \in (0,\Lambda^*_k)$, depending only on $k$, so that  for any
 $n \ge k$:
{\small  \[\mu_{C,(L_{C,n})^\times_\lambda} \ge - 1/3 \cdot \left(\frac{k-2}{2}\right)^2,\mm{ for any }\lambda \le \lambda^*_k \mm{ and hence }\]
\[\alpha_+ \ge - (1- \sqrt{2/3}) \cdot \frac{k-2}{2}, \,\, \alpha_- \le - (1+ \sqrt{2/3}) \cdot \frac{k-2}{2}.\]}
Moreover, for $\lambda \in (0,\lambda^*_k]$  we have $\mu^*_{C,(L_{C,n})^\times_\lambda} <-\eta^*_\lambda <0$ and hence
{\small  \[-\vartheta^*_\lambda >  \alpha_+ \ge - (1- \sqrt{2/3}) \cdot \frac{k-2}{2}>  - \frac{k-2}{2} >- (1+ \sqrt{2/3}) \cdot \frac{k-2}{2} \ge \alpha_- > \vartheta^*_\lambda -(k-2),\]}
 for some constants $\eta^*_\lambda, \vartheta^*_\lambda>0$, depending only on $\lambda$ and $k$.
\end{itemize}
\end{proposition}

{\bf Proof} \quad  To check the skin adpatedness assertions we note that  $\frac{n-2}{4 (n-1)} < 1/4$. Thus, we have from (\ref{kwsy}):
 {\small
\begin{equation}    \int_C | \nabla f |^2 - \frac{n-2}{4 (n-1)} \cdot |A|^2  \cdot  f^2 \, d A  \ge \int_C | \nabla f |^2 - 1/4 \cdot |A|^2  \cdot  f^2 \, d A \ge \end{equation}
\[  \ge \int_C \frac{k}{2 (k-1)} \cdot  |  \nabla f |^2 + \frac{k- 3}{4 (k-1)}\cdot   | A |^2   \cdot  f^2\,  d A\]
and for some $\tau^*_k >0$ depending only on the dimension we get the estimate:
\[ \ge\frac{k-3}{4 (k-1)}\cdot \int_C f  \cdot  L_{\bot} f   \, d A \ge \tau^*_k  \cdot \frac{k- 3}{4 (k-1)}\cdot \int_C \bp^2 \cdot f^2   \, d A=:\Lambda^*_k \cdot \int_C \bp^2 \cdot f^2   \, d A.\]}\\
For the eigenvalue estimate on $S_C$ we write for any $\lambda < \Lambda^*_k$:
\begin{equation}   \int_{S_C}  f  \cdot (L_{C,n})_\lambda f  \,  dA  = \int_{S_C} | \nabla f |^2 - \left(\frac{n-2}{4 (n-1)} \cdot | A |^2 +\lambda \cdot \bp^2\right)  \cdot  f^2 \, d A \ge \end{equation}
\[\int_{S_C} | \nabla f |^2 - 1/4 \cdot |A|^2 \cdot f^2 dA - \max\{0,\lambda/\tau_k\}  \cdot \int_{S_C}|\nabla f|^2  + |A|^2 \cdot f^2 dA.\]
For $\lambda/\tau_k < 1/16$, and $|f|_{L^2}=1$, this can be lower estimated by
\[\ge 1/3 \cdot  \int_{S_C} | \nabla f |^2 -  |A|^2 \cdot f^2 dA \ge  1/3 \cdot  \mu_{C,(J_C)^\times_0} \ge - 1/3 \cdot \left(\frac{k-2}{2}\right)^2.\]
Hence, for $\lambda^*_k:=1/2 \cdot \min \{\Lambda^*_k, 1/16\}$, we have
$\mu_{C,(L_{C,n})^\times_\lambda} \ge - 1/3 \cdot \left(\frac{k-2}{2}\right)^2,$ for $\lambda \le \lambda^*_k$.\\

 For $\lambda \in (0,\lambda^*_k]$, we observe again, as in \ref{evee}, that there are $f \in C_0^\infty(S_C \setminus \Sigma_{S_C})$ with
\[ \int_{S_C} | \nabla f|^2 + \left(\frac{n-2}{4 (n-1)} \cdot scal_C|_{S_C} - \lambda \cdot \bp|_{S_C}^2\right)  \cdot  f^2 \, d A < 0.\]
Similarly as before, there are constants $\eta^*_\lambda, \vartheta^*_\lambda>0$, depending only on $\lambda>0$ and $k$, so that
\[\mu^*_{C,(L_{C,n})^\times_\lambda} <-\eta^*_\lambda<0 \mm{ and, therefore, } 0>-\vartheta^*_\lambda>\alpha_+ \mm{ and }\alpha_- > \vartheta^*_\lambda -(k-2)>  -(k-2).\]
for any singular area minimizing $k$-dimensional cone $C^k$.  \qed

Complementary to the previous discussions, where we mostly focussed on the radial growth rate, we now describe some global properties of the spherical component $\psi_C(\omega)$ of $\Psi_\pm[n,k](\omega,r) = \psi_C(\omega) \cdot r^{\alpha_\pm}$
defined on $S_C \setminus \Sigma_{S_C}$. In turn, these results ensure uniform estimates for the radial growth, near the singular set.

\begin{proposition}\emph{\textbf{(Global Harnack Estimates for $\psi_C$)}}   \label{skk} \, For any $k$ with $n \ge k \ge 7$, we consider the two solutions $\Psi_\pm[n,k](\omega,r) = \psi_C(\omega) \cdot r^{\alpha_\pm}$ of  $(L_{C^k,n})_\lambda \, w=0$, for some $\lambda \in (0,\lambda^*_k]$,  on a cone $C^k \in {\cal{SC}}_k$, corresponding to $0$ resp. $\infty$ in the Martin boundary $\widehat{\sigma_C}$.\\

Then, there are constants $a_{n,k,\lambda} >0$, depending only on $n,k,\lambda$, and $b_{n,k,\lambda,\rho}$ additionally depending on $\rho>0$, but not on the cone $C^k$:
\begin{enumerate}
  \item $|\psi_C|_{L^p(S_C \setminus \Sigma_{S_C})}<\infty$, for $ p< \frac {k-1} {k-3}$.
  \item $|\psi_C|_{L^1(S_C \setminus \Sigma_{S_C})} \le a_{n,k,\lambda} \cdot \inf_{\omega \in S_C \setminus \Sigma_{S_C}}  \psi_C(\omega)$.
  \item $\sup_{\omega \in\E(\rho)}\psi_C(\omega) \le  b_{n,k,\lambda,\rho} \cdot|\psi_C|_{L^1(S_C \setminus \Sigma_{S_C})}$, for $\E(\rho)=\{x \in S_C \,|\, \bp(x) \le \rho^{-1}\}$.
\end{enumerate}
Similarly, $|v|_{L^q(B_1(0) \cap C^k)}<\infty$, for $q< \frac {k} {k-2}$, for any solution $v>0$ of  $(L_{C^k,n})_\lambda \, w=0$
\end{proposition}

\textbf{Proof} \quad For $0 <\rho^2:= \lambda/\eta^*_\lambda < -\lambda/\mu^*_{C,(L_{C,n})^\times_\lambda}$ and any $C \in {\cal{SH}}^{\R}_n$ we have
 \[\frac{n-2}{4 (n-1)} \cdot |A|_C^2(\omega,1) + \lambda \cdot (\bp^\times)^2(\omega)  +\mu^*_{C,(L_{C,n})^\times_\lambda}>0 \mm{ on } \I(\rho).\]
Thus,  we observe, $\psi_C(\omega)>0$ is a superharmonic function on $\I(\rho)$:
  $- \Delta_{S_C} \psi_C(\omega) \ge$
   \[- \Delta_{S_C} \psi_C(\omega) - \Big(\frac{n-2}{4 (n-1)} \cdot |A|_C^2(\omega,1) + \lambda \cdot (\bp^\times)^2(\omega) \Big)\cdot \psi_C(\omega) - \mu^*_{C,(L_{C,n})^\times_\lambda} \cdot \psi_C(\omega) = 0\]
   Although $S_C \subset \p B_1(0)$ is not a global area minimizer, it is an almost minimizer and it shares the regularity theory with proper area minimizers and we can also locally apply the Bombieri-Giusti Harnack inequality [BG],Th.6,p.39 in the following form: For any superharmonic $w > 0$ defined on the regular portion of on a extrinsically measured ball $B_R(x) \cap S_C$ of  sufficiently small radius $R>0$ we have
\[0< \left\{\frac {1}{Vol_{n-1}(S_C \cap B_r(x))}\int w^p \right\}^{1/p} \leqq C  \cdot  \inf_{B_r(x)}w\]
for $ r \leqq \beta_n \cdot R $ and $ p< \frac {k-1} {k-3}$ for some constants $C = C(S_C,p) , \beta_n > 0$. We apply this to a finite cover $B_r(p_j)$, $j=1,...m$ of $\Sigma_{S_C}$ by sufficiently small balls so that $B_R(p_j) \subset \I(\rho)$.\\

Since the complement of these open balls in $S_C \setminus \Sigma_{S_C}$ is compact, we see that $\inf_{S_C \setminus \Sigma_{S_C}}  \psi_C(\omega) >0$ and
$|\psi_C|_{L^p(S_C \setminus \Sigma_{S_C})} < \infty$ and we also obtain the (trivial) estimate $\sup_{\omega \in\E(\rho)}\psi_C(\omega) \le  b_{C,\rho} \cdot|\psi_C|_{L^1(S_C \setminus \Sigma_{S_C})}$, for some suitably large $b_{C,n,k,\lambda,\rho}>0$, for each individual cone $C$. \\

 From the, up to multiples, uniqueness of $\psi_C$, the compactness of ${\cal{SC}}_k$, the naturality of $|A|$ and $\bp$ and the standard elliptic theory for $(L_{C,n})_\lambda^\times$, we infer the existence of some common $a_{n,k,\lambda} >0$ for all  $C^k \in {\cal{SC}}_k$, so that \[|\psi_C|_{L^1(S_C \setminus \Sigma_{S_C})} \le a_{n,k,\lambda} \cdot \inf_{\omega \in S_C \setminus \Sigma_{S_C}}  \psi_C(\omega),\]
 for any $C^k \in {\cal{SC}}_k$. This way we also get a common $b_{n,k,\lambda,\rho}>0$ for all cones  $C^k \in {\cal{SC}}_k$ so that  $\sup_{\omega \in\E(\rho)}\psi_C(\omega) \le b_{n,k,\lambda,\rho} \cdot|\psi_C|_{L^1(S_C \setminus \Sigma_{S_C})}$.\\

 The assertion that $|v|_{L^q(B_1(0) \cap C^k)}<\infty$, for $q< \frac {k} {k-2}$, for any solution $v>0$ of $(L_{C^k,n})_\lambda \, w=0$ follows completely similarly from the Bombieri-Giusti Harnack inequality. \qed

\end{document}